\documentclass[a4paper]{amsart}
\pdfoutput=1
\usepackage{mathtools}
\usepackage{amssymb}
\usepackage{enumerate}
\usepackage{xcolor}
\usepackage{graphicx}
\usepackage{tikz}
\definecolor{accent}{cmyk}{0.94,0.11,0,0.15}
\usetikzlibrary{arrows.meta}
\usepackage{pdflscape}
\usepackage{afterpage}
\makeatletter
\newcommand{\centerfloatpage}{\setlength{\@fptop}{0pt plus 1fil}\setlength{\@fpbot}{0pt plus 1fil}}
\makeatother
\usepackage{microtype}
\usepackage[margin=1in]{geometry}
\usepackage[utf8]{inputenc}
\usepackage[pagebackref, bookmarks=true, colorlinks=true, linkcolor=blue, citecolor=blue, urlcolor=blue]{hyperref}
\usepackage{cleveref}

\usepackage{aliascnt}
\numberwithin{equation}{section}
\newtheorem{theorem}{Theorem}[section]
\crefname{theorem}{Theorem}{Theorems}\Crefname{theorem}{Theorem}{Theorems}
\crefname{section}{Section}{Sections}\crefname{subsection}{Section}{Sections}
\crefname{figure}{Figure}{Figures}\Crefname{figure}{Figure}{Figures}
\newcommand{\newsharedtheorem}[2]{%
  \newaliascnt{#1}{theorem}\newtheorem{#1}[#1]{#2}\aliascntresetthe{#1}%
  \crefname{#1}{#2}{#2s}\Crefname{#1}{#2}{#2s}}
\newsharedtheorem{proposition}{Proposition}
\newsharedtheorem{corollary}{Corollary}
\newsharedtheorem{lemma}{Lemma}
\theoremstyle{definition}
\newsharedtheorem{definition}{Definition}
\newsharedtheorem{example}{Example}
\theoremstyle{remark}
\newsharedtheorem{remark}{Remark}

\newcommand{\bbC}{\mathbb{C}}
\newcommand{\mfh}{\mathfrak{h}}
\newcommand{\mfg}{\mathfrak{g}}
\newcommand{\mfn}{\mathfrak{n}}
\newcommand{\mfk}{\mathfrak{k}}
\newcommand{\mfp}{\mathfrak{p}}
\newcommand{\mcG}{\mathcal{G}}
\newcommand{\mcH}{\mathcal{H}}
\newcommand{\mcN}{\mathcal{N}}
\newcommand{\mcB}{\mathcal{B}}
\newcommand{\mcC}{\mathcal{C}}
\newcommand{\mcS}{\mathcal{S}}
\newcommand{\mcL}{\mathcal{L}}
\newcommand{\mcI}{\mathcal{I}}
\newcommand{\oA}{\overline{\mathcal{A}}}
\newcommand{\bfi}{\mathbf{i}}
\newcommand{\tbfi}{\widetilde{\boldsymbol{\imath}}}
\newcommand{\bfj}{\mathbf{j}}
\newcommand{\bx}{\mathbf{x}}
\newcommand{\by}{\mathbf{y}}
\newcommand{\tB}{\widetilde{B}}
\newcommand{\tG}{\widetilde{\mcG}}
\newcommand{\tg}{\widetilde{\gamma}}
\newcommand{\vt}{\vartheta}
\newcommand{\Vt}{\Theta}
\newcommand{\Kil}[2]{\langle #1, #2\rangle}
\newcommand{\rr}{\mathrm{r}}
\newcommand{\cc}{\mathrm{c}}
\newcommand{\ex}{\mathrm{ex}}
\newcommand{\fr}{\mathrm{fr}}
\newcommand{\pr}{\mathrm{pr}}
\newcommand{\len}{\operatorname{len}}
\newcommand{\Cas}{\mathsf t}
\newcommand{\toric}{\Pi}
\newcommand{\Sp}{Sp_{2n}}
\newcommand{\SL}{SL_{2n}}
\newcommand{\newword}[1]{\textbf{\emph{#1}}}

\begin{document}
\title{Exotic cluster structures on $Sp_{2n}$} \author{Idan Eisner} \address{Department of Computer Science, Tel-Hai University of the Galilee, Kiryat Shmona, Israel} \email{eisnerida@telhai.ac.il} \author{Salvatore Stella} \address{Universit\`a degli Studi dell'Aquila, Dipartimento di Ingegneria e Scienze dell'Informazione e Matematica, Via Vetoio, 67100 L'Aquila, Italy} \email{salvatore.stella@univaq.it} \thanks{Salvatore Stella was partially supported by PRIN 20223FEA2E and INdAM/GNSAGA} \date{\today}

\begin{abstract}
  We modify the seeds of Berenstein, Fomin and Zelevinsky attached to reduced words for~$(w_0,w_0)$, for any simply connected simple complex algebraic group and any aperiodic Belavin--Drinfeld pair.
  In every type, moves of words act on the resulting seeds by mutations.
  After a change of coordinates introduced by Gekhtman, Shapiro and Vainshtein, the seeds are, for $SL_N$, seeds of the cluster structures of these authors, and for $Sp_{2n}$ they define regular complete cluster structures with the expected properties.
  This proves the aperiodic version of their conjecture for $Sp_{2n}$, with its uniqueness part in a refined form.
\end{abstract}
\maketitle

\section{Introduction}
\label{sec: intro}

Let $\mcG$ be a simply connected simple complex algebraic group of rank $r$.
Belavin and Drinfeld \cite{BD82} classified the quasitriangular solutions of the classical Yang--Baxter equation on $\mfg=\mathrm{Lie}(\mcG)$ in terms of combinatorial data, which we encode as \emph{BD maps} $\gamma_\bullet:\Gamma_\bullet\to\Gamma$ between subsets of simple roots.
A pair~$\gamma=(\gamma_{\rr},\gamma_{\cc})$ of BD maps, together with some Cartan data $s\in\mcS(\gamma)$, defines a Poisson bracket~$\{\cdot,\cdot\}_{\gamma,s}$ on $\mcG$, where $\gamma_{\rr}$ and $\gamma_{\cc}$ govern the behavior under left and right multiplication, respectively.
For $\gamma=(0,0)$ and $s=0$ this is the standard Poisson--Lie structure, which is \emph{compatible} with the cluster structure of Berenstein, Fomin and Zelevinsky \cite{BFZ05}, i.e.\ log-canonical with respect to every extended cluster.
Gekhtman, Shapiro and Vainshtein conjectured \cite[Conjecture~3.2]{GSV11} that every Poisson--Lie structure of this kind is compatible with a regular complete cluster structure on $\mcG$.
In the setting of pairs \cite{GSV19}, every \emph{aperiodic} pair $\gamma$ (\cref{def: aperiodic}) should admit a cluster structure $\mcC_\gamma$ on $\mcG$ with the following properties.
\begin{enumerate}[(i)]
  \item The number of frozen variables is $2r-|\Gamma_{\rr}|-|\Gamma_{\cc}|$, and the exchange matrices have full rank.
  \item $\mcC_\gamma$ is regular and its upper cluster algebra is $\bbC[\mcG]$.
  \item The global toric action on $\mcC_\gamma$ is generated by the action of $\mcH_{\gamma_{\rr}}\times\mcH_{\gamma_{\cc}}$ by left and right multiplication.
    Here $\mcH_{\gamma_\bullet}$ is the identity component of the subgroup of the maximal torus on which the characters $\alpha$ and $\gamma_\bullet(\alpha)$ agree for every $\alpha\in\Gamma_\bullet$.
  \item $\mcC_\gamma$ is compatible with $\{\cdot,\cdot\}_{\gamma,s}$ for every $s\in\mcS(\gamma)$.
  \item The Poisson brackets compatible with $\mcC_\gamma$ are exactly the brackets $q\{\cdot,\cdot\}_{\gamma,s}+\toric_\Xi$, where $q\in\bbC$, $s\in\mcS(\gamma)$, and $\toric_\Xi$ is a bivector field built from the infinitesimal action of $\mcH_{\gamma_{\rr}}\times\mcH_{\gamma_{\cc}}$.
\end{enumerate}
Property (v) extends \cite[Conjecture~3.2(v)]{GSV11}, which concerns Poisson--Lie brackets, to all Poisson brackets.
For Poisson--Lie brackets, (v) gives the conclusion of \cite[Conjecture~3.2(v)]{GSV11} when $q\neq0$, while for $q=0$ some toric brackets are also Poisson--Lie (\cref{rmk: v}).
Property (v) follows from a general criterion: for every $\mcG$ and every BD pair $\gamma$, aperiodic or not, a cluster structure on $\mcG$ with (i) and (iii), compatible with $\{\cdot,\cdot\}_{\gamma,s}$ with non-zero constants $\kappa_\ell$ in \eqref{eq: compatibility}, and with block-indecomposable exchange matrices, satisfies (v) (\cref{thm: classification}).
For $\mcG=SL_N$, cluster structures with (i)--(iv) are constructed in \cite[Theorem~3.3]{GSV19} for aperiodic oriented pairs, and cluster structures with (i), (ii) and (iv) in \cite{GSV23} for all aperiodic pairs; the latter also have (iii) and (v) (\cref{thm: A}(c), \cref{cor: classification}).
In the latter construction, the initial seed is described by a quiver, and its cluster variables are obtained from standard ones through a rational map $h_\gamma$, defined in \cite{GSV23} for every $\mcG$.
In this paper we give seeds for every $\mcG$ and every double reduced word, and, using the same rational map $h_\gamma$, prove the conjecture for $\mcG=Sp_{2n}$.

\subsection*{Main results}
For a double reduced word $\bfi$, i.e.\ a reduced word for $(w_0,w_0)$ in $W\times W$, let $\Sigma(\bfi)$ be the seed of \cite[Theorem~2.10]{BFZ05}, whose cluster variables are generalized minors indexed by the chambers of~$\bfi$.
The seed $\Sigma^\Phi_\gamma(\bfi)$ of \cref{def: BD word seed} is obtained from $\Sigma(\bfi)$ by three operations.
We thaw one end chamber for each root in $\Gamma_{\rr}\sqcup\Gamma_{\cc}$.
We add to the exchange matrix the entries of some virtual letters and dipoles determined by $\gamma$.
Finally, we multiply each minor by a monomial in the minors of the end chambers.

\begin{theorem}
  \label{thm: A}
  Let $\gamma$ be an aperiodic BD pair of $\mcG$.
  \begin{enumerate}[(a)]
    \item The seeds $\Sigma^\Phi_\gamma(\bfi)$ of all double reduced words are mutation equivalent with moves on words corresponding to (sequences of) mutations (\cref{thm: moves,cor: mutation equivalent}).
    \item With the thawed chambers frozen, $\Sigma^\Phi_\gamma(\bfi)$ and $\Sigma(\bfi)$ are quasi-isomorphic in the sense of \cite{Fra16} via the identity of $\bbC(\mcG)$ (\cref{cor: Fraser}).
    \item If $\mcG=SL_N$, the seeds $\Sigma^\Phi_\gamma(\bfi)$ composed with $h_\gamma^{-1}$ are seeds of the cluster structure of \cite{GSV23}, and its global toric action is generated by $\mcH_{\gamma_{\rr}}\times\mcH_{\gamma_{\cc}}$ (\cref{sec: type A}).
  \end{enumerate}
\end{theorem}

Whenever $h_\gamma$ is birational, (a) holds for the \emph{BD word seeds} $\Sigma_\gamma(\bfi)$, obtained by composing the functions of $\Sigma^\Phi_\gamma(\bfi)$ with $h_\gamma^{-1}$.
This is the case for $SL_N$ \cite{GSV23} and, as we show, for $Sp_{2n}$.

\begin{theorem}
  \label{thm: main}
  Let $\gamma$ be an aperiodic BD pair of type $C_n$.
  Then $h_\gamma$ is birational, and the BD word seeds $\Sigma_\gamma(\bfi)$ are seeds of a single cluster structure $\mcC_\gamma$ on $\Sp$ (\cref{thm: all words}), which satisfies (i)--(v) (\cref{thm: structure,thm: complete,cor: classification}).
\end{theorem}

As $-w_0$ acts trivially on the simple roots of type $C_n$, a Poisson--Lie pair $\gamma_{\cc}=\gamma_{\rr}^{-1}\neq0$ is never aperiodic, so the nontrivial Poisson--Lie classes of $Sp_{2n}$ are not covered by \cref{thm: main} (\cref{rmk: aperiodic not PL}).

The group $\Sp$ is the fixed point set of an involution $\vt$ of $\SL$, and $\gamma$ lifts to a $\vt$-invariant aperiodic pair of type $A_{2n-1}$.
The BD word seeds of lifted words are $\vt$-symmetric seeds of the cluster structures of \cite{GSV23}, and their folds are the BD word seeds of $\Sp$ (\cref{thm: fold}).
We deduce the properties of $\mcC_\gamma$ from those of the ambient cluster structure, using that $\bbC[\Sp]$ is factorial for completeness.

\subsection*{Organization}
\Cref{sec: conventions} fixes conventions and recalls the results of \cite{GSV19,GSV23}.
\Cref{sec: word seeds} defines the BD word seeds for every $\mcG$ and proves \cref{thm: A}(a) and~(b), and \cref{sec: type A} proves \cref{thm: A}(c).
\Cref{sec: Sp} describes $\Sp$ inside $\SL$, BD pairs of type $C_n$ and their lifts, and the restriction of brackets.
\Cref{sec: folding} treats the folding of symmetric seeds, \cref{sec: symmetric} shows that the BD word seeds of lifted words are symmetric with the BD word seeds of $\Sp$ as folds, and \cref{sec: main} proves \cref{thm: main} except for (v).
\Cref{sec: brackets} classifies the Poisson brackets compatible with a cluster structure on any $\mcG$, for any BD pair (\cref{thm: classification}), and deduces (v) for $SL_N$ and $\Sp$.

\subsection*{Acknowledgments}
We are grateful to M. Gekhtman, M. Shapiro and A. Vainshtein for suggesting this problem to us and for the many illuminating conversations on the topic over the long years it took us to complete it.

\section{Conventions and the results of Gekhtman, Shapiro and Vainshtein}
\label{sec: conventions}
For integers $a,b$ we write $[a]_+:=\max(a,0)$ and, if $a\le b$, $[a,b]:=\{a,\dots,b\}$.
For a statement $\mathsf P$, $[\mathsf P]$ is~$1$ if $\mathsf P$ holds and $0$ otherwise.

\subsection{Cluster algebras}
\label{sec: cluster conventions}
We use the conventions of \cite{BFZ05}; in particular exchange matrices are \emph{tall}.
A seed is a pair $(\bx,\tB)$, where $\bx=(x_k)_{k\in\mcI}$ is an extended cluster of algebraically independent generators of a field.
Its finite index set $\mcI=\mcI_\ex\sqcup\mcI_\fr$ is partitioned into exchangeable and frozen indices.
The integer matrix $\tB=(b_{k\ell})_{k\in\mcI,\ell\in\mcI_\ex}$ has skew-symmetrizable principal part $B=(b_{k\ell})_{k,\ell\in\mcI_\ex}$.
We call $B$ \newword{block-indecomposable} if no permutation of $\mcI_\ex$ makes it block diagonal, i.e.\ if the graph on $\mcI_\ex$ with edges $\{k,\ell\}$, $b_{k\ell}\neq0$, is connected; as $\mu_k$ preserves block decompositions, this is a mutation invariant.
The mutation in direction $k\in\mcI_\ex$ replaces $x_k$ by $x'_k$ with
\begin{equation}
  \label{eq: exchange}
  x_kx'_k=\prod_{\ell\in\mcI}x_\ell^{[b_{\ell k}]_+}+\prod_{\ell\in\mcI}x_\ell^{[-b_{\ell k}]_+},
\end{equation}
and $\tB$ by $\mu_k(\tB)$, with entries $-b_{\ell m}$ if $\ell=k$ or $m=k$, and $b_{\ell m}+[b_{\ell k}]_+[b_{km}]_+-[-b_{\ell k}]_+[-b_{km}]_+$ otherwise.
We write $\mu_k(\bx)$ for the new extended cluster and $\mu_k(\bx,\tB)=(\mu_k(\bx),\mu_k(\tB))$ for the new seed.
The papers \cite{GSV10,GSV19,GSV23} use \emph{wide} matrices, which are the negative transposes of ours.
Both conventions give the same exchange relations, and mutation commutes with $\tB\mapsto-\tB^T$, so they define the same seeds.
For an extended cluster $\by$ let $\mcL(\by)$ be the ring of Laurent polynomials in the exchangeable variables of $\by$ with coefficients polynomials in the frozen ones.
Frozen variables are never inverted.
The \newword{upper cluster algebra} $\oA:=\oA(\bx)$ is the intersection of the rings $\mcL(\by)$ over the extended clusters $\by$ mutation equivalent to $\bx$; it only depends on the cluster structure defined by $(\bx,\tB)$.
If $\tB$ has full rank, $\oA(\bx)$ equals the \newword{upper bound} $\mcL(\bx)\cap\bigcap_{k\in\mcI_\ex}\mcL(\mu_k(\bx))$ \cite[Theorem~3.11]{GSV18dd}.
This result is stated for generalized cluster structures and for a ground ring in which frozen variables need not be inverted.
Ordinary cluster structures are the case of trivial strings; see also \cite[Theorem~3.21 and Section~3.4]{GSV10}.
A cluster structure on an irreducible variety $X$ is a mutation class of seeds in $\bbC(X)$.
It is \newword{regular} if all cluster variables are regular functions on $X$, and \newword{complete} if moreover $\oA=\bbC[X]$.

A Poisson bracket $\{\cdot,\cdot\}$ on the field containing the extended cluster $\bx$ is \newword{log-canonical} with respect to $\bx$ if~$\{x_k,x_\ell\}=\omega_{k\ell}x_kx_\ell$ with $\omega_{k\ell}\in\bbC$, and \newword{compatible} with a cluster structure if it is log-canonical with respect to every extended cluster.
Assume that the bracket is log-canonical with respect to $\bx$, with coefficient matrix $\Omega=(\omega_{k\ell})$, and that $\tB$ has full rank.
Then it is compatible if and only if
\begin{equation}
  \label{eq: compatibility}
  \sum_{m\in\mcI}b_{m\ell}\,\omega_{mk}=\delta_{k\ell}\,\kappa_\ell
  \qquad(k\in\mcI,\ \ell\in\mcI_\ex)
\end{equation}
for some $\kappa_\ell\in\bbC$ \cite[Proposition~2.2]{GSV11}, see also \cite[Theorem~4.5]{GSV10}, stated there for the wide matrix~$-\tB^T$.
Some $\kappa_\ell$ may vanish, for example for the zero bracket.
Conversely, \eqref{eq: compatibility} with all $\kappa_\ell\neq0$ implies both compatibility and full rank of $\tB$ \cite[Proposition~2.5]{GSV18dd}.
Since $\Omega$ is skew-symmetric, so is~$\tB^T\Omega\tB$, whose entry at $(\ell,m)\in\mcI_\ex^2$ is $\kappa_\ell b_{\ell m}$ by \eqref{eq: compatibility}; the $\kappa_\ell$ need not be positive.
An integer matrix~$P=(P_{kj})$ with rows indexed by $\mcI$ and columns indexed by $[1,q]$ defines an action of the torus~$(\bbC^*)^q$ on the field generated by $\bx$, the \newword{local toric action} $x_k\mapsto x_k\prod_jt_j^{P_{kj}}$ for $(t_1,\dots,t_q)\in(\bbC^*)^q$.
Its \newword{rank} is the rank of~$P$, i.e.\ the dimension of the image of $(\bbC^*)^q$ in the group of automorphisms of the field.
It is a \newword{global toric action} if every cluster variable of every seed mutation equivalent to $(\bx,\tB)$ is homogeneous for it, so that it rescales all cluster variables.
This holds if and only if $\tB^TP=0$, i.e.\ the two monomials of every exchange relation have the same weight \cite[Lemma~5.3]{GSV10}.
If $\tB$ has full rank, $\ker\tB^T$ has dimension~$|\mcI_\fr|$, so the rank of a global toric action is at most $|\mcI_\fr|$ \cite[Section~2.2]{GSV11}.

\subsection{Lie theory}
\label{sec: Lie}
Let $\mcG$ be a simply connected simple complex algebraic group of rank $r$, with Lie algebra~$\mfg$.
We fix a maximal torus $\mcH$, with Lie algebra $\mfh$, and a pinning: simple roots $\Gamma=\{\alpha_i\}_{i\in[1,r]}$, simple coroots $\alpha_i^\vee$, and Chevalley generators $E_i,F_i$ with $[E_i,F_i]=\alpha_i^\vee$.
We write $\lambda(\zeta)$ for the pairing of $\lambda\in\mfh^*$ and $\zeta\in\mfh$, $\varpi_i$ for the fundamental weights, and $a_{ij}=\alpha_j(\alpha_i^\vee)$ for the Cartan matrix entries.
We follow the conventions of \cite[Section~1.1]{Kac90} for the Cartan matrix and the numbering of \cite[Chapter~VI, Plate~III]{Bou02}, so that for $C_n$ the root $\alpha_n$ is long, $a_{n-1,n}=-2$ and $a_{n,n-1}=-1$.

\emph{Subgroups.} Let $\mfn_\pm$ be the nilpotent subalgebras spanned by the positive, resp.\ negative, root vectors,~$\mcN_\pm=\exp\mfn_\pm$ the corresponding maximal unipotent subgroups, and $\mcB_\pm=\mcH\mcN_\pm$ the opposite Borel subgroups.
For $K\subseteq\Gamma$, let $\mcG^K$ be the Levi subgroup generated by $\mcH$ and the one-parameter subgroups $\exp(tE_i)$, $\exp(tF_i)$, $\alpha_i\in K$; we write $\mcN_\pm^K=\mcN_\pm\cap\mcG^K$ and $\mcB_\pm^K=\mcB_\pm\cap\mcG^K$.

\emph{The Weyl group.} Let $W=N_\mcG(\mcH)/\mcH$ be the Weyl group, generated by the simple reflections $s_i$, with length function $\len$ and longest element $w_0$.
A \newword{reduced word} for $w\in W$ is a sequence $(j_1,\dots,j_m)$ with~$w=s_{j_1}\cdots s_{j_m}$ and $m=\len(w)$.
For $w\in W$ let $\overline{w}\in\mcG$ be the representative of \cite[Section~1.4]{FZ99}: $\overline{s_i}:=\exp(-E_i)\exp(F_i)\exp(-E_i)$, and $\overline{uv}=\overline{u}\,\overline{v}$ whenever $\len(uv)=\len(u)+\len(v)$.
For $K\subseteq\Gamma$, the Weyl group of $\mcG^K$ is the parabolic subgroup $W_K$ generated by the $s_i$, $\alpha_i\in K$, and $w_0^K$ is its longest element.
The involution $i\mapsto i^*$ of $[1,r]$ is defined by $w_0(\alpha_i)=-\alpha_{i^*}$; then $w_0(\alpha_i^\vee)=-\alpha_{i^*}^\vee$ and $w_0\varpi_i=-\varpi_{i^*}$.
If~$(j_1,\dots,j_m)$ is a reduced word for $w\in W$, the roots $s_{j_1}\cdots s_{j_{p-1}}(\alpha_{j_p})$, $p\in[1,m]$, are positive and pairwise distinct, and they are the positive roots $\beta$ with $w^{-1}\beta<0$ \cite[Chapter~VI, Section~1.6, Corollary~2 of Proposition~17]{Bou02}, applied to $w^{-1}$.
For $w=w_0$ every positive root occurs exactly once.

\emph{The invariant form.} Fix a non-degenerate invariant symmetric bilinear form $\Kil{\cdot}{\cdot}$ on $\mfg$.
Let $\Cas\in\mfg\otimes\mfg$ be the Casimir element, $\Cas=\Cas_0+\Cas_{-+}+\Cas_{+-}$ with $\Cas_0\in\mfh\otimes\mfh$, $\Cas_{-+}\in\mfn_-\otimes\mfn_+$, $\Cas_{+-}=\Cas_{-+}^{21}$, where $\Psi^{21}$ denotes~$\Psi\in\mfg\otimes\mfg$ with its two tensor factors swapped.
We write $\xi\wedge\eta:=\xi\otimes\eta-\eta\otimes\xi$, and $\lambda\otimes1:\mfh\otimes\mfh\to\mfh$ for $\zeta\otimes\zeta'\mapsto\lambda(\zeta)\zeta'$.

\emph{Generalized minors.} Multiplication $\mcN_-\times\mcH\times\mcN_+\to\mcG$ is an isomorphism onto an open dense subset $\mcN_-\mcH\mcN_+$, the \newword{big cell}.
For $U$ in it we write $U=U_-U_0U_+$, $U_\pm\in\mcN_\pm$, $U_0\in\mcH$, for its \newword{Gauss factorization}, which is unique.
Following \cite{FZ99}, let $\Delta_{\varpi_i,\varpi_i}$ be the regular function on $\mcG$ given by~$\Delta_{\varpi_i,\varpi_i}(U)=U_0^{\varpi_i}$ on the big cell.
For $u,v\in W$ let the \newword{generalized minor} $\Delta_{u\varpi_i,v\varpi_i}$ be the function $U\mapsto\Delta_{\varpi_i,\varpi_i}(\overline{u}^{\,-1}U\overline{v})$.
For $SL_N$ with the standard pinning these are the minors of the matrix $U$ (\cref{sec: type A}).
It depends only on the weights $u\varpi_i$ and $v\varpi_i$, not on the choice of $u$ and $v$ \cite[Section~2.3]{FZ99}, which justifies the notation.
Let $U\mapsto U^T$ be the anti-involution of $\mcG$ whose differential swaps $E_i$ and $F_i$ and fixes $\mfh$; for $SL_N$ it is the transpose.
Then $\overline{w}^{\,T}=\overline{w}^{\,-1}$, so $\Delta_{u\varpi_i,v\varpi_i}(U^T)=\Delta_{v\varpi_i,u\varpi_i}(U)$.
A function $f$ on $\mcG$ is \newword{$\mcH$-homogeneous of bidegree} $(\lambda^\rr,\lambda^\cc)$ if $f(y_{\rr}Uy_{\cc})=y_{\rr}^{\lambda^\rr}y_{\cc}^{\lambda^\cc}f(U)$ for $y_{\rr},y_{\cc}\in\mcH$; for instance~$\Delta_{u\varpi_i,v\varpi_i}$ has bidegree $(u\varpi_i,v\varpi_i)$.

\subsection{Belavin--Drinfeld data and Poisson brackets}
\label{sec: BD}
A \newword{BD map} is an injective map $\gamma_\bullet:\Gamma_\bullet\subseteq\Gamma\to\Gamma$ which is an isometry for $\Kil{\cdot}{\cdot}$ and is nilpotent.
Nilpotency means that for every $\alpha\in\Gamma_\bullet$ some $\gamma_\bullet^m(\alpha)$ is not in $\Gamma_\bullet$.
Thus this is a BD triple $(\Gamma_\bullet,\gamma_\bullet(\Gamma_\bullet),\gamma_\bullet)$ of \cite{BD82}.
The empty map is denoted $0$.
Single BD maps always appear as the components $\gamma_{\rr},\gamma_{\cc}$ of a BD pair (defined below), so we write them with an index~$\bullet\in\{\rr,\cc\}$ and attach the same index to the objects built from them.
An unindexed $\gamma$ always denotes a pair.
We identify $\Gamma$ with $[1,r]$ and write $\gamma_\bullet(i)=j$ for $\gamma_\bullet(\alpha_i)=\alpha_j$.
The assignment $E_i\mapsto E_j$, $F_i\mapsto F_j$ for $\gamma_\bullet(i)=j$ extends to an isomorphism $\mfg_{\Gamma_\bullet}\to\mfg_{\gamma_\bullet(\Gamma_\bullet)}$.
Here, for $K\subseteq\Gamma$, we write $\mfg_K$ for the semisimple subalgebra generated by the $E_i,F_i$, $\alpha_i\in K$.
It is the derived algebra of $\mathrm{Lie}(\mcG^K)$.
Precomposing with the orthogonal projection onto $\mfg_{\Gamma_\bullet}$ gives an endomorphism $\gamma_\bullet$ of $\mfg$, nilpotent on $\mfn_\pm$ and commuting with the transpose.
Set $C_{\gamma_\bullet}:=\frac12+\sum_{k\ge1}\gamma_\bullet^k$ on $\mfn_+\oplus\mfn_-$, and
\[
  \mfh_{\gamma_\bullet}:=\{\zeta\in\mfh\mid\alpha(\zeta)=\gamma_\bullet(\alpha)(\zeta)\ \forall\alpha\in\Gamma_\bullet\},
  \qquad
  \mcH_{\gamma_\bullet}:=\{y\in\mcH\mid y^\alpha=y^{\gamma_\bullet(\alpha)}\ \forall\alpha\in\Gamma_\bullet\}^{\circ},
\]
where $\circ$ denotes the identity component; thus $\mcH_{\gamma_\bullet}=\exp\mfh_{\gamma_\bullet}$ is a torus and $\dim\mfh_{\gamma_\bullet}=|\Gamma|-|\Gamma_\bullet|$.
We take the identity component because the subgroup defined by the character equations can be disconnected.
For example, let $\mcG=SL_3$ and $\gamma_\bullet=\{1\to2\}$.
For $y=\mathrm{diag}(t_1,t_2,t_3)$ the only equation is $t_1/t_2=t_2/t_3$, that is $t_1t_3=t_2^2$.
Together with $t_1t_2t_3=1$, this gives $t_2^3=1$.
So the subgroup is~$\{\mathrm{diag}(t,\omega,\omega^2t^{-1})\mid t\in\bbC^*,\ \omega^3=1\}$, which has three components, and $\mcH_{\gamma_\bullet}=\{\mathrm{diag}(t,1,t^{-1})\mid t\in\bbC^*\}$.
In terms of characters, $\alpha_1-\alpha_2=3(\varpi_1-\varpi_2)$ is not primitive in the weight lattice.
Equivalently, the subgroup contains the center $\{\omega I\mid\omega^3=1\}$ of $SL_3$, on which every root is trivial, and it is the direct product of $\mcH_{\gamma_\bullet}$ and the center.

Let $\mcS(\gamma_\bullet)\subseteq\mfh\wedge\mfh$ be the set of solutions $s_\bullet$ of
\begin{equation}
  \label{eq: system s}
  2\big((\alpha-\gamma_\bullet(\alpha))\otimes1\big)s_\bullet=\big((\alpha+\gamma_\bullet(\alpha))\otimes1\big)\Cas_0
  \qquad\forall\alpha\in\Gamma_\bullet;
\end{equation}
by \cite{BD82} it is a non-empty affine space over $\mfh_{\gamma_\bullet}\wedge\mfh_{\gamma_\bullet}$.
For $s_\bullet\in\mcS(\gamma_\bullet)$,
\[
  r_{\gamma_\bullet,s_\bullet}:=\tfrac12\Cas+s_\bullet+(1\otimes C_{\gamma_\bullet})\Cas_{-+}-(C_{\gamma_\bullet}\otimes1)\Cas_{+-}
\]
is a classical $R$-matrix \cite{BD82}.

For $\xi\in\mfg$ and a function $f$ on $\mcG$ let
\[
  \partial^\rr_\xi f(U):=\tfrac{d}{dt}f(\exp(t\xi)U)|_{t=0},
  \qquad
  \partial^\cc_\xi f(U):=\tfrac{d}{dt}f(U\exp(t\xi))|_{t=0},
\]
and for $\Psi=\sum_u\xi_u\otimes\eta_u\in\mfg\otimes\mfg$ let $\Psi(f,g)^\rr:=\sum_u\partial^\rr_{\xi_u}f\,\partial^\rr_{\eta_u}g$, and similarly $\Psi(f,g)^\cc:=\sum_u\partial^\cc_{\xi_u}f\,\partial^\cc_{\eta_u}g$.

A \newword{BD pair} is a pair $\gamma=(\gamma_{\rr},\gamma_{\cc})$ of BD maps $\gamma_{\rr}:\Gamma_{\rr}\to\Gamma$, $\gamma_{\cc}:\Gamma_{\cc}\to\Gamma$.
We set $\mcS(\gamma):=\mcS(\gamma_{\rr})\times\mcS(\gamma_{\cc})$.
For $s=(s_{\rr},s_{\cc})\in\mcS(\gamma)$ let $r_\rr:=r_{\gamma_{\rr},s_{\rr}}$ and $r_\cc:=\big((T\otimes T)r_{\gamma_{\cc},s_{\cc}}\big)^{21}$, where $T$ is the transpose, and define
\begin{equation}
  \label{eq: bracket}
  \{f,g\}_{\gamma,s}:=r_\cc(f,g)^\cc-r_\rr(f,g)^\rr.
\end{equation}
In the notation of \cite{GSV19,GSV23}, $\gamma_{\rr}=\gamma^\rr$ and $\gamma_{\cc}=(\gamma^\cc)^{-1}$, and \eqref{eq: bracket} is the bracket of \cite[Section~2.2]{GSV23} for the pair $(r_\rr,r_\cc)$, which is Poisson \cite[Section~2.2]{GSV23}, \cite[Section~2.2]{GSV19}.
We set~$\mfh_\gamma:=\mfh_{\gamma_{\rr}}\oplus\mfh_{\gamma_{\cc}}$ and $\mcH_\gamma:=\mcH_{\gamma_{\rr}}\times\mcH_{\gamma_{\cc}}$, acting on $\mcG$ by $(y_{\rr},y_{\cc})\cdot U=y_{\rr}Uy_{\cc}$.

\subsection{The map \texorpdfstring{$h_\gamma$}{h}}
\label{sec: h}
Let $\gamma$ be a BD pair for $\mcG$, and put $\sigma_{\rr}=w_0^{\Gamma_{\rr}}w_0$ and $\sigma_{\cc}=w_0^{\Gamma_{\cc}}w_0$.
A generic $U\in\mcG$ admits factorizations
\begin{equation*}
  U=V_{\rr}M_{\rr},\quad V_{\rr}\in\mcN_-^{\Gamma_{\rr}},\ M_{\rr}\in\mcB_+\overline{\sigma_{\rr}}\mcB_+
  \qquad\text{and}\qquad
  U=M_{\cc}V_{\cc},\quad V_{\cc}\in\mcN_+^{\Gamma_{\cc}},\ M_{\cc}\in\mcB_-\overline{\sigma_{\cc}}^{\,-1}\mcB_-,
\end{equation*}
which determine $V_{\rr}$ and $V_{\cc}$.
Write $\lfloor A\rfloor_+$ (resp.\ $\lceil A\rceil_-$) for the factor in $\mcN_+$ (resp.\ $\mcN_-$) of $A\in\mcN_+\mcB_-$ (resp.~$A\in\mcB_+\mcN_-$), and put
\begin{equation*}
  \overline{V}_{\rr}=\lfloor V_{\rr}\overline{w_0^{\Gamma_{\rr}}}\rfloor_+\in\mcN_+^{\Gamma_{\rr}}
  \qquad\text{and}\qquad
  \overline{V}_{\cc}=\lceil \overline{w_0^{\Gamma_{\cc}}}V_{\cc}\rceil_-\in\mcN_-^{\Gamma_{\cc}}.
\end{equation*}
Let $\boldsymbol\gamma_{\rr}:\mcN_+\to\mcN_+$ and $\boldsymbol\gamma_{\cc}:\mcN_-\to\mcN_-$ be the homomorphisms integrating the restrictions of the endomorphisms $\gamma_{\rr}$ to $\mfn_+$ and $\gamma_{\cc}$ to $\mfn_-$; their images lie in $\mcN_+^{\gamma_{\rr}(\Gamma_{\rr})}$ and $\mcN_-^{\gamma_{\cc}(\Gamma_{\cc})}$.
The map $h_\gamma:\mcG\dashrightarrow\mcG$ is
\begin{equation}
  \label{eq: h}
  h_\gamma(U)=H_{\rr}(U)\,U\,H_{\cc}(U),
  \qquad
  H_{\rr}(U)=\cdots\boldsymbol\gamma_{\rr}^2(\overline{V}_{\rr})\,\boldsymbol\gamma_{\rr}(\overline{V}_{\rr}),
  \qquad
  H_{\cc}(U)=\boldsymbol\gamma_{\cc}(\overline{V}_{\cc})\,\boldsymbol\gamma_{\cc}^2(\overline{V}_{\cc})\cdots,
\end{equation}
where the products are finite because $\gamma_{\rr}$ and $\gamma_{\cc}$ are nilpotent.
This is the map of \cite[Section~3.1]{GSV23}, defined there for every simple $\mcG$.
For $\Sp$ we study the map $h_\gamma$ in \cref{prop: h theta}.

\begin{lemma}
  \label{lem: h transpose}
  For every BD pair $\gamma$ of $\mcG$, $h_{(\gamma_{\rr},\gamma_{\cc})}(U)=h_{(\gamma_{\cc},\gamma_{\rr})}(U^T)^T$.
\end{lemma}
\begin{proof}
  The anti-involution $T$ fixes $\mcH$, swaps $\mcN_+^K$ and $\mcN_-^K$, and intertwines the homomorphisms integrating~$\gamma_{\rr}$ and $\gamma_{\cc}$.
  Also, $\overline{w}^{\,T}=\overline{w}^{\,-1}\in\mcH\overline{w}$ for an involution $w$.
  Hence replacing $\overline{w}$ by $\overline{w}^{\,-1}$, for $w=w_0^{\Gamma_{\rr}},w_0^{\Gamma_{\cc}}$, only changes the Borel factors in the definitions of $\overline{V}_{\rr},\overline{V}_{\cc}$.
  Transposing the factorizations of $U^T$ for~$(\gamma_{\cc},\gamma_{\rr})$ therefore gives those of $U$ for $(\gamma_{\rr},\gamma_{\cc})$.
  Thus $H_{\cc}(U)$, resp.\ $H_{\rr}(U)$, is the transpose of $H_{\rr}(U^T)$, resp.\ $H_{\cc}(U^T)$, computed for $(\gamma_{\cc},\gamma_{\rr})$.
\end{proof}

\begin{definition}[{\cite[Section~4.1]{GSV23}}]
  \label{def: aperiodic}
  Take disjoint copies $\{a^\rr\mid a\in\Gamma_{\rr}\}$ and $\{b^\cc\mid b\in\Gamma_{\cc}\}$ of $\Gamma_{\rr}$ and $\Gamma_{\cc}$.
  The \newword{alternating map} of a BD pair $\gamma$ is the partial injective map on their union sending $a^\rr$ to $(\gamma_{\rr}(a)^*)^\cc$ if $\gamma_{\rr}(a)^*\in\Gamma_{\cc}$, and $b^\cc$ to $(\gamma_{\cc}(b)^*)^\rr$ if $\gamma_{\cc}(b)^*\in\Gamma_{\rr}$.
  The pair $\gamma$ is \newword{aperiodic} if the alternating map has no periodic points.
\end{definition}

For $SL_N$ the periodic points of the alternating map correspond to the alternating cycles of the BD graph of \cite[Section~4.1]{GSV23}, so \cref{def: aperiodic} agrees with the definition there.

\subsection{The cluster structures on \texorpdfstring{$SL_N$}{SLN}}
\label{sec: GSV}
\begin{theorem}[{\cite[Section~1, Theorem~4.4 and (4.6)]{GSV23}}]
  \label{thm: GSV h}
  If $\gamma$ is an aperiodic BD pair of type $A_{N-1}$, the map $h_\gamma$ is birational.
\end{theorem}

Birationality of $h_\gamma$ is proved in \cite{GSV23} only in type $A$.
For $\Sp$ we obtain it in \cref{thm: fold}(c).

\begin{theorem}[{\cite[Theorem~3.3]{GSV19}, \cite[Theorems~4.3, 4.11, 4.14, 4.17]{GSV23}}]
  \label{thm: GSV}
  Let $\gamma$ be an aperiodic BD pair of type $A_{N-1}$.
  The seed of \cite[Theorem~4.11]{GSV23} defines a cluster structure $\mcC_\gamma$ on $SL_N$ with the properties (i)--(iv) of the introduction:
  \begin{enumerate}[(i)]
    \item the number of frozen variables is $2(N-1)-|\Gamma_{\rr}|-|\Gamma_{\cc}|$ and the exchange matrices have full rank;
    \item $\mcC_\gamma$ is regular and $\oA=\bbC[SL_N]$;
    \item the global toric action on $\mcC_\gamma$ is generated by the action of $\mcH_\gamma$;
    \item $\mcC_\gamma$ is compatible with $\{\cdot,\cdot\}_{\gamma,s}$ for every $s\in\mcS(\gamma)$; in every extended cluster its coefficient matrix satisfies \eqref{eq: compatibility} with all $\kappa_\ell\neq0$.
  \end{enumerate}
\end{theorem}

For (ii), \cite[Theorem~4.14]{GSV23} gives regularity of the adjacent cluster variables, and \cite[Theorem~4.17]{GSV23} gives Laurent expressions of all matrix entries in terms of the initial cluster and of each adjacent cluster; the latter is proved by induction on $|\Gamma_{\rr}|+|\Gamma_{\cc}|$, the case of size $1$ being \cite[Section~7.3]{GSV19}.
Regularity and completeness follow as in \cite[Sections~4.4--4.5]{GSV23}.
For $\gamma=(0,0)$, completeness of the standard cluster structure on $SL_N$ is \cite{Oya25} (and \cite[Theorem~B.1]{QY25} in all types).
For (iv), the proof of \cite[Theorem~4.11]{GSV23} shows that at the initial seed the coefficient matrix satisfies \eqref{eq: compatibility}, and that all~$\kappa_\ell$ are equal to some $\kappa\neq0$.
The coefficient matrices in adjacent clusters are related as the matrices~$\Lambda$ of \cite[(3.4)]{BZ05}, and the proof of \cite[Proposition~3.4]{BZ05} shows that $\tB^T\Omega$ does not change, so the same~$\kappa_\ell$ work in every extended cluster.
Part (iii) is proved in \cite[Theorem~3.3]{GSV19} for the class of oriented pairs considered there, and is not stated in \cite{GSV23}.
We prove it for arbitrary aperiodic pairs after \cref{prop: base case}, using only (i), (ii), and the initial cluster described there.

\section{The BD word seeds}
\label{sec: word seeds}
In this section we define the BD word seeds for every $\mcG$ and prove that moves of words act on them by mutations.

\subsection{Words, slots and chambers}
\label{sec: words}
We describe the seeds of \cite[Section~2.2]{BFZ05} in a form in which the two ends of a word play the same role.
We need this because the BD word seeds are obtained by adding letters at both ends of a word.
The translation to the labels of \cite{BFZ05} is given in \cref{prop: BFZ}.
Readers can follow along using \cref{fig: Sp6 word seed}, ignoring for the moment all the parts in color.

A \newword{word} is a family $\bfi=(i_j)_{j\in J}$ of \newword{letters} $i_j\in\pm[1,r]$, indexed by a finite non-empty interval $J$ of integers, the \newword{positions}.
The \newword{level} of $i_j$ is $|i_j|$, and $\varepsilon(i_j)$ is its sign.
The \newword{slots} of $\bfi$ are the integers $z$ such that $z$ or $z+1$ is a position.
The letter $i_j$ lies between the slots $j-1$ and $j$, and the first and the last slot are the \newword{end slots}.
For each $x\in[1,r]$, the letters at level $x$ cut the slots into intervals, the \newword{chambers} of $\bfi$ at level $x$.
We write $c_x(z)$ for the one at level $x$ containing the slot $z$, and $C(\bfi)$ for the set of all chambers.
A chamber is \newword{$\bfi$-exchangeable} if it contains no end slot, and an \newword{end chamber} otherwise.
So if $\bfi$ has $n$ letters at level $x$, there are $n+1$ chambers at level $x$; the \newword{initial} one contains the first slot and the \newword{final} one the last slot (they are equal if $n=0$).
If $x=|i_j|$, the letter $i_j$ separates $c_x(j-1)$ from $c_x(j)$.

For distinct $x,y\in[1,r]$, let $M_{xy}$ be the set of positions of the letters at levels $x$ and $y$, and for $j\in M_{xy}$ put $\varphi_{xy}(j):=\varepsilon(i_j)$ if $|i_j|=y$ and $\varphi_{xy}(j):=-\varepsilon(i_j)$ if $|i_j|=x$.
Let $\tB(\bfi)=(b_{k\ell})$ be the matrix with rows labeled by $C(\bfi)$ and columns by the $\bfi$-exchangeable chambers, where
\begin{equation}
  \label{eq: B entries}
  b_{k\ell}=
  \begin{cases}
    \varepsilon(i_j) & \text{if } (k,\ell)=\big(c_x(j-1),c_x(j)\big) \text{ with } x=|i_j|,\\
    -\varepsilon(i_j) & \text{if } (k,\ell)=\big(c_x(j),c_x(j-1)\big) \text{ with } x=|i_j|,\\
    -\frac12a_{xy}\big(\varphi_{xy}(p')-\varphi_{xy}(p)\big) & \text{if } (k,\ell)=\big(c_x(p),c_y(p)\big) \text{ with } p<p' \text{ consecutive in } M_{xy},\\
    0 & \text{otherwise.}
  \end{cases}
\end{equation}
The cases are exclusive.
The first two concern chambers at equal levels and the third at distinct levels.
In the third case, distinct pairs of consecutive positions give distinct pairs of chambers: if $p_1<p'_1\le p_2<p'_2$, the letter at $p'_1$ separates the chambers of $p_1$ and $p_2$ at its level.
So the entries of $\tB(\bfi)$ lie in $\{0,\pm1\}$ at equal levels and in $\{0,\pm a_{xy}\}$ at levels $x\neq y$.
As $\varphi_{yx}=-\varphi_{xy}$, if $d_xa_{xy}=d_ya_{yx}$ with $d_x>0$, then $d_xb_{k\ell}=-d_yb_{\ell k}$ for exchangeable chambers $k,\ell$ at levels $x,y$.
So the principal part of $\tB(\bfi)$ is skew-symmetrizable, and skew-symmetric if the Cartan matrix is symmetric.
Only the order of the letters in a word matters: translating all positions by the same integer changes neither the chambers, nor their exchangeability, nor $\tB(\bfi)$.
We use this freely.

\begin{proposition}
  \label{prop: BFZ}
  Let $J=[1,L]$, and for $k\in J$ label the chamber $c_{|i_k|}(k)$ by $k$, and, for $x\in[1,r]$, the initial chamber $c_x(0)$ by $-x$.
  Then $\tB(\bfi)$ is the matrix of \cite[Definition~2.3]{BFZ05}.
\end{proposition}
\begin{proof}
  With these labels, the chamber $k$ is $\bfi$-exchangeable if and only if $k$ is $\bfi$-exchangeable in the sense of \cite[Section~2.2]{BFZ05}.
  For every label $k$ the next label at the same level is $k^+$, the smallest index~$j>k$ with~$|i_j|=|i_k|$, or $L+1$ if there is none.
  Here $i_{-x}:=-x$ as in \cite{BFZ05}, so that $|i_k|$ is the level of the chamber $k$ for every label $k$.
  By \cite[Remark~2.4]{BFZ05}, for $\ell$ exchangeable, $p=\max(k,\ell)$ and~$p'=\min(k^+,\ell^+)$, the entry of the matrix of \cite{BFZ05} is
  \begin{equation}
    \label{eq: BFZ entries}
    b_{k\ell}=
    \begin{cases}
      -\operatorname{sgn}(k-\ell)\,\varepsilon(i_p) & \text{if } p=p',\\
      -\operatorname{sgn}(k-\ell)\,\varepsilon(i_p)\,a_{|i_k|,|i_\ell|} & \text{if } p<p' \text{ and } \varepsilon(i_p)\varepsilon(i_{p'})(k-\ell)(k^+-\ell^+)>0,\\
      0 & \text{otherwise,}
    \end{cases}
  \end{equation}
  where $p$ and $p'$ lie in $[1,L]$, so that $\varepsilon(i_p)$ and $\varepsilon(i_{p'})$ are signs of letters.
  If $p=p'$, then $k$ and $\ell$ are at the same level, $p$ is the letter separating them, and this is one of the first two cases of \eqref{eq: B entries}.
  Otherwise, if $k$ and $\ell$ are at the same level, or at distinct levels $x,y$ with $a_{xy}=0$, both formulas give $0$.
  Let $k$ and $\ell$ be at distinct levels $x,y$ with $a_{xy}\neq0$.
  If $p>p'$, the chambers have no common slot, and both formulas give~$0$.
  If $p<p'$, then $p$ and $p'$ are consecutive positions of letters at the levels $x$ and $y$.
  The slot $p$ lies in both chambers, so $k=c_x(p)$ and $\ell=c_y(p)$.
  Conversely, every such pair $p<p'$ arises in this way.
  We check that the two formulas agree in each of the four cases for $p\in\{k,\ell\}$ and $p'\in\{k^+,\ell^+\}$.
  For instance, if $p=k$ and $p'=\ell^+$, then $(k-\ell)(k^+-\ell^+)>0$, so $b_{k\ell}=-a_{xy}\varepsilon(i_p)$ if $\varepsilon(i_p)\varepsilon(i_{p'})>0$ and $0$ otherwise, while $\varphi_{xy}(p)=-\varepsilon(i_p)$ and $\varphi_{xy}(p')=\varepsilon(i_{p'})$.
  If $p=k$ and $p'=k^+$, then $b_{k\ell}=-a_{xy}\varepsilon(i_p)$ if $\varepsilon(i_p)\varepsilon(i_{p'})<0$ and $0$ otherwise, while $\varphi_{xy}(p)=-\varepsilon(i_p)$ and $\varphi_{xy}(p')=-\varepsilon(i_{p'})$.
  The cases $p=\ell$ follow in the same way with the signs of $\varphi_{xy}$ reversed.
\end{proof}

A \newword{commutation} swaps two consecutive letters $i_j,i_{j+1}$ at different levels which have opposite signs or satisfy $a_{|i_j|,|i_{j+1}|}=0$.
A \newword{sign swap} replaces two consecutive letters $i,-i$ by $-i,i$, or conversely.
A \newword{braid move} replaces $m$ consecutive letters $i,j,i,\dots$ by $j,i,j,\dots$, where $i$ and $j$ have the same sign,~$a_{|i|,|j|}\neq0$.
The number $m$, the \newword{length} of the move, is the order of $s_{|i|}s_{|j|}$ in $W$, i.e.\ $m=3,4,6$ for $a_{|i|,|j|}a_{|j|,|i|}=1,2,3$.
These are the \newword{moves}.

A move replaces the letters at the positions in an interval $[p,p']$ by letters with the same set of levels.
If~$\bfi'$ is obtained from $\bfi$ by a move, a chamber $c_x(z)$ of $\bfi$ with $z\notin[p,p'-1]$ \newword{corresponds} to the chamber~$c_x(z)$ of $\bfi'$.
The other chambers, of $\bfi$ and of $\bfi'$, are \newword{replaced} by the move.
A commutation replaces no chamber.
A sign swap, resp.\ a braid move of length $3$, replaces one chamber of $\bfi$ and one of $\bfi'$, with the slot $p$, resp.\ the slots $p,p+1$, and we let these correspond to each other.

\begin{lemma}
  \label{lem: commutation}
  Let $\bfi'$ be obtained from $\bfi$ by a commutation.
  Then the correspondence of chambers is a bijection $C(\bfi)\to C(\bfi')$ which preserves levels and exchangeability, and $\tB(\bfi')$ is obtained from $\tB(\bfi)$ by relabeling rows and columns along it.
\end{lemma}
\begin{proof}
  Let the commutation swap the letters at positions $j,j+1$, with levels $x\neq y$.
  The correspondence is a bijection since no chamber is replaced, and it preserves levels and end chambers, hence exchangeability.

  By \eqref{eq: B entries}, the entries between chambers at given levels depend only on the signed sequence of letters at those levels.
  The commutation leaves this sequence unchanged at each single level and at every pair of levels other than $\{x,y\}$.
  Thus it remains to check the entries between levels $x$ and $y$.

  If $a_{xy}=0$, these entries vanish.
  Otherwise the swapped letters have opposite signs, therefore $\varphi_{xy}(j)=\varphi_{xy}(j+1)$, and the commutation leaves both $M_{xy}$ and $\varphi_{xy}$ unchanged.
  Every consecutive pair~$p<p'$ in $M_{xy}$ with $p\neq j$ therefore gives the same entry at corresponding chambers, since the slot $p$ lies outside the move.
  The remaining pair $j<j+1$ contributes zero in both words.
  Hence all entries agree along the correspondence.
\end{proof}

\subsection{Double reduced words}
\label{sec: double words}
A double reduced word for $(u,v)\in W\times W$ is a shuffle of a reduced word for $u$ in the alphabet $-[1,r]$ and one for $v$ in the alphabet $[1,r]$ \cite[Section~2.2]{BFZ05}.
We only use~$(u,v)=(w_0,w_0)$, and call a double reduced word $\bfi=(i_1,\dots,i_L)$ for $(w_0,w_0)$ simply a \newword{double reduced word}.
So $J=[1,L]$, the slots are $[0,L]$, $L=2\len(w_0)$, and every level carries letters of both signs.
Following \cite[(2.9)--(2.11)]{BFZ05}, for a slot $z$ let $u_{\le z}$ be the product of the $s_{|i_j|}$ over the negative letters $i_j$ with $j\le z$, in increasing order of $j$.
Let $v_{>z}$ be the product of the $s_{|i_j|}$ over the positive letters $i_j$ with $j>z$, in decreasing order of $j$.
The \newword{row weight} and the \newword{column weight} of the chamber $k=c_x(z)$ are $\lambda^\rr_k:=u_{\le z}\varpi_x$ and $\lambda^\cc_k:=v_{>z}\varpi_x$, which do not depend on $z\in k$, as the letters between two slots of~$k$ are at levels $y\neq x$ and $s_y\varpi_x=\varpi_x$.
Its \newword{minor} is $\Delta(k;\bfi):=\Delta_{\lambda^\rr_k,\lambda^\cc_k}$.
By \cite[Theorem~2.10]{BFZ05} and \cref{prop: BFZ}, $\Sigma(\bfi):=\big((\Delta(k;\bfi))_{k\in C(\bfi)},\tB(\bfi)\big)$ is a seed in $\bbC(\mcG)$; in particular the $\Delta(k;\bfi)$ are algebraically independent and generate $\bbC(\mcG)$.

The \newword{root} of a negative letter $i_j$ is $u_{\le j-1}(\alpha_{|i_j|})$, and that of a positive letter $i_j$ is $v_{>j}(\alpha_{|i_j|})$.
The negative letters, in increasing order of position, form a reduced word for $w_0$, and so do the positive letters in decreasing order of position.
So by \cref{sec: Lie} the roots of the negative letters are the positive roots, each occurring once, and similarly for the positive letters.
In particular, every simple root is the root of exactly one negative and exactly one positive letter.

Moves transform double reduced words into double reduced words, and any two double reduced words are related by moves, by \cite[Chapter~IV, Section~1.5, Proposition~5]{Bou02} applied to $W\times W$ and to the monoid presented by the braid relations.

\begin{theorem}[{\cite[Theorem~3.5]{SSVZ00}, \cite[Propositions~3.7 and~3.25]{SW21}}]
  \label{thm: SSVZ}
  Let $\bfi'$ be obtained from $\bfi$ by a move.
  Along the correspondence of chambers:
  \begin{enumerate}
    \item For a commutation, $\Sigma(\bfi')=\Sigma(\bfi)$.
    \item For a sign swap or a braid move of length $3$, $\Sigma(\bfi')=\mu_k(\Sigma(\bfi))$, where $k$ is the replaced chamber.
    \item For a braid move of length $4$, $\Sigma(\bfi')=\mu_k\mu_\ell\mu_k(\Sigma(\bfi))$, where $k$ and $\ell$ are the replaced chambers at the levels of the long and of the short root (equivalently $\mu_\ell\mu_k\mu_\ell(\Sigma(\bfi))$); the replaced chambers of $\bfi'$ correspond to them by level.
    \item For a braid move of length $6$ at the positions $[p,p+5]$, with $\alpha_x$ short for $x=|i_p|$, and $y=|i_{p+1}|$, $\Sigma(\bfi')$ is obtained from $\Sigma(\bfi)$ by the mutations at
      \[
        \ell_2,\ell_1,k_2,k_1,\ell_2,k_2,\ell_2,\ell_1,k_1,\ell_2,
      \]
      in this order, where $k_1=c_x(p)$, $k_2=c_x(p+2)$, $\ell_1=c_y(p+1)$, $\ell_2=c_y(p+3)$, corresponding to $c_x(p+1)$, $c_x(p+3)$, $c_y(p)$, $c_y(p+2)$ in $\bfi'$.
      If $\alpha_{|i_p|}$ is long, the inverse move is of this form, and $\Sigma(\bfi')$ is obtained from $\Sigma(\bfi)$ by the same mutations in the reverse order.
  \end{enumerate}
  The statement on the matrices also holds for non-reduced words.
\end{theorem}
For simply-laced $\mcG$, where the moves have length at most $3$, the statement on the matrices is \cite[Theorem~3.5]{SSVZ00}, see \cite[Remark~2.14]{BFZ05}, where the non-simply-laced case is left open.
In general the statement follows from \cite[Proposition~3.7]{SW21}.
The strings of the string diagram of $\bfi$ are our chambers, the closed ones being the $\bfi$-exchangeable chambers.
The Cartan matrix is ours \cite[Appendix~A.1]{SW21}, and the exchange matrix of \cite[Section~3.1]{SW21} is $\epsilon_{\ell k}=b_{k\ell}$.
For a braid move of length $6$, the strings $a,b,c,d$ of \cite[Proposition~3.7(2)]{SW21} are $k_1,k_2,\ell_1,\ell_2$.
In \cite[Proposition~3.7]{SW21} the braid move of length $4$ mutates first at the short root, which gives the same seed, as the alternating mutations at two chambers~$k,\ell$ with $b_{k\ell}b_{\ell k}=-2$ have period $6$.
The statement on the functions is \cite[Proposition~3.25]{SW21}, through the realization of the decorated double Bott--Samelson cell \cite[Section~3.3]{SW21} sending $U\in\mcG^{w_0,w_0}$ to the chains of decorated flags $U^T\overline{u_{\le z}}\,\mcN_+$ and $\mcN_-\overline{v_{>z}}^{\,-1}$, $z$ a slot, after deleting repetitions.
This is not the map of \cite[Section~1.4]{SW21}.
These chains represent the same configuration under the canonical identifications for different reduced words, by uniqueness of compatible refinements \cite[Lemmas~2.10--2.11]{SW21}.
The coordinate \cite[Definition~3.20]{SW21} of $k=c_x(z)$ is $\Delta_{\varpi_x,\varpi_x}\big(\overline{v_{>z}}^{\,-1}U^T\overline{u_{\le z}}\big)=\Delta(k;\bfi)(U)$, by \cref{sec: Lie}.

\begin{remark}
  As any two double reduced words are related by moves, \cref{thm: SSVZ} shows that the seeds~$\Sigma(\bfi)$ of all double reduced words are related by sequences of mutations, up to the correspondence of chambers.
  The quantum version of this claim is \cite[Theorem~A]{QY25}.
\end{remark}

\subsection{The BD word seeds}
\label{sec: BD word seeds}
The colored parts of \cref{fig: Sp6 word seed} illustrate the following definition.
Fix an aperiodic BD pair $\gamma$ of $\mcG$ and a double reduced word $\bfi$, i.e.\ a reduced word for $(w_0,w_0)$, and call the elements of $\Gamma_{\rr}$ and of $\Gamma_{\cc}$ \newword{row roots} and \newword{column roots}.

\begin{definition}
  \label{def: BD word seed}
  \begin{enumerate}[(a)]
    \item For a row root $a$, its \newword{inversion letter} is the negative letter with root $\alpha_{\gamma_{\rr}(a)}$, and~$p^\rr_a$ denotes its position; its \newword{thawed chamber} is $t^\rr_a:=c_{a^*}(L)$.
      For a column root $b$, its inversion letter is the positive letter with root $\alpha_{\gamma_{\cc}(b)}$, and $p^\cc_b$ denotes its position; its thawed chamber is $t^\cc_b:=c_{b^*}(0)$.
    \item Order the row roots as $a_1,\dots,a_{|\Gamma_{\rr}|}$ and the column roots as $b_1,\dots,b_{|\Gamma_{\cc}|}$ so that $p^\rr_{a_1}>p^\rr_{a_2}>\cdots$ and~$p^\cc_{b_1}>p^\cc_{b_2}>\cdots$.
      The \newword{extended word} $\bfi^\gamma$ is obtained from $\bfi$ by adding \newword{virtual letters} $b_1^*,\dots,b_{|\Gamma_{\cc}|}^*$ at positions $1-|\Gamma_{\cc}|,\dots,0$ and $-a_1^*,\dots,-a_{|\Gamma_{\rr}|}^*$ at positions $L+1,\dots,L+|\Gamma_{\rr}|$:
      \[
        \bfi^\gamma=(b_1^*,\dots,b_{|\Gamma_{\cc}|}^*,\,i_1,\dots,i_L,\,-a_1^*,\dots,-a_{|\Gamma_{\rr}|}^*).
      \]
      For $z\in[0,L]$ we identify the chamber $c_x(z)$ of $\bfi$ with the chamber $c_x(z)$ of $\bfi^\gamma$.
    \item Let $\mcI:=C(\bfi)$, let $\mcI_\ex$ consist of the $\bfi$-exchangeable chambers and of the thawed chambers, and let~$\mcI_\fr:=\mcI\setminus\mcI_\ex$.
      Let $\tB(\bfi^\gamma)_\mcI$ be the submatrix of $\tB(\bfi^\gamma)$ with rows $\mcI$ and columns $\mcI_\ex$.
      For a row or column root $a$, with inversion letter $i_p$ at level $x$ and thawed chamber $t$, let $D_a$ be the matrix with rows $\mcI$ and columns $\mcI_\ex$ whose only nonzero entries are
      \[
        d_{c_x(p-1),t}=-d_{t,c_x(p-1)}=-\varepsilon(i_p),
        \qquad
        d_{c_x(p),t}=-d_{t,c_x(p)}=\varepsilon(i_p),
      \]
      an entry being omitted if its column is not in $\mcI_\ex$.
      The \newword{dipole matrix} is $D_\gamma(\bfi):=\sum_aD_a$, and~$\tB_\gamma(\bfi):=\tB(\bfi^\gamma)_\mcI+D_\gamma(\bfi)$.
    \item For $k\in\mcI$, a row root $a$ and a column root $b$ let
      \[
        e^\rr_a(k):=\big[\lambda^\rr_k(\alpha^\vee_{\gamma_{\rr}(a)})\big]_+,
        \qquad
        e^\cc_b(k):=\big[\lambda^\cc_k(\alpha^\vee_{\gamma_{\cc}(b)})\big]_+,
      \]
      and define recursively
      \begin{equation}
        \label{eq: Phi}
        \Phi_\gamma(k;\bfi):=\Delta(k;\bfi)\prod_{a\in\Gamma_{\rr}}\Phi_\gamma(t^\rr_a;\bfi)^{e^\rr_a(k)}\prod_{b\in\Gamma_{\cc}}\Phi_\gamma(t^\cc_b;\bfi)^{e^\cc_b(k)}.
      \end{equation}
    \item If $h_\gamma$ is birational, the \newword{BD word seed} is $\Sigma_\gamma(\bfi):=\big((x_\gamma(k;\bfi))_{k\in\mcI},\tB_\gamma(\bfi)\big)$, $x_\gamma(k;\bfi):=\Phi_\gamma(k;\bfi)\circ h_\gamma^{-1}$.
  \end{enumerate}
\end{definition}

\begin{lemma}
  \label{lem: tails}
  Let $x\in[1,r]$, $a\in\Gamma_{\rr}$ and $b\in\Gamma_{\cc}$.
  The weights of the final chamber $c_x(L)$ are $(-\varpi_{x^*},\varpi_x)$ and those of the initial chamber $c_x(0)$ are $(\varpi_x,-\varpi_{x^*})$, and
  \[
    e^\rr_a(c_x(L))=e^\cc_b(c_x(0))=0,
    \qquad
    e^\cc_b(c_x(L))=\delta_{\gamma_{\cc}(b),x},
    \qquad
    e^\rr_a(c_x(0))=\delta_{\gamma_{\rr}(a),x}.
  \]
  So in \eqref{eq: Phi} the function $\Phi_\gamma(c_x(L);\bfi)$ involves exactly the $\Phi_\gamma(t^\cc_b;\bfi)$ with $\gamma_{\cc}(b)=x$, and $\Phi_\gamma(c_x(0);\bfi)$ exactly the $\Phi_\gamma(t^\rr_a;\bfi)$ with $\gamma_{\rr}(a)=x$, each with exponent $1$.
  As $\gamma$ is aperiodic, the recursion \eqref{eq: Phi} terminates.
  In particular the functions $\Phi_\gamma(k;\bfi)$ of the end chambers $k$ are monomials in the minors of end chambers, and do not depend on $\bfi$.
\end{lemma}
\begin{proof}

  The weights of $c_x(L)$ are computed at the slot $L$, with $u_{\le L}=w_0$ and $v_{>L}=e$, and they are~$w_0\varpi_x=-\varpi_{x^*}$ and $\varpi_x$.
  Those of $c_x(0)$ are computed at the slot $0$, with $u_{\le0}=e$ and $v_{>0}=w_0$, and are $\varpi_x$ and $-\varpi_{x^*}$.
  Fundamental weights pair non-negatively with simple coroots, and $\varpi_x(\alpha^\vee_y)=\delta_{xy}$.
  This gives the four formulas.
  As $t^\rr_a=c_{a^*}(L)$, $\Phi_\gamma(t^\rr_a;\bfi)$ involves $\Phi_\gamma(t^\cc_b;\bfi)$ if and only if $\gamma_{\cc}(b)=a^*$, i.e.\ the alternating map, which sends $b^\cc$ to $(\gamma_{\cc}(b)^*)^\rr$, sends $b^\cc$ to $a^\rr$.
  Similarly $\Phi_\gamma(t^\cc_b;\bfi)$ involves $\Phi_\gamma(t^\rr_a;\bfi)$ if and only if it sends $a^\rr$ to $b^\cc$.
  Thus the recursion from a thawed chamber follows the alternating map backwards and terminates since this map is injective and has no periodic points.
  As the functions of end chambers only involve those of thawed chambers, which are end chambers, they are monomials in minors of end chambers.
  These minors do not depend on $\bfi$, since the weights of end chambers do not.
\end{proof}

\cref{def: BD word seed} makes sense: the thawed chambers are distinct, as those of row roots are final and those of column roots initial, and every level carries at least two letters.
They are $\bfi^\gamma$-exchangeable, as the virtual letter at their level separates them from the end slot of $\bfi^\gamma$, while the other end chambers of $\bfi$ remain end chambers of $\bfi^\gamma$, as there is no virtual letter beyond them at their level.
So $\mcI_\ex$ is the set of~$\bfi^\gamma$-exchangeable chambers of $\bfi$, and $\tB(\bfi^\gamma)_\mcI$ is the restriction of $\tB(\bfi^\gamma)$ to $\mcI$.
The dipole entries are off the diagonal, i.e.\ $t^\rr_a$ is neither of the two chambers separated by the inversion letter of $a$.
Indeed, paired with~$\alpha^\vee_{\gamma_{\rr}(a)}$, the row weights of these chambers give $1$ and $-1$.
In contrast, the row weight $-\varpi_a$ of $t^\rr_a$ gives $0$, as $\gamma_{\rr}(a)\neq a$ by nilpotency.
For a column root the column weights give $-1$, $1$ and $0$.
The exponents $e^\rr_a(k)$, $e^\cc_b(k)$ in \eqref{eq: Phi} are $0$ or $1$ for $SL_N$ (\cref{sec: type A}), but exponents equal to $2$ occur in type $C_n$ (\cref{ex: Sp6}).

By \cref{lem: commutation}, only the relative order of virtual letters at adjacent levels matters.
For adjacent row roots $a,a'$, the virtual letters of $a$ and $a'$ are in the order opposite to that of $p^\rr_a$ and $p^\rr_{a'}$, and similarly for column roots.

We now check that the pair $\Sigma^\Phi_\gamma(\bfi):=\big((\Phi_\gamma(k;\bfi))_{k\in\mcI},\tB_\gamma(\bfi)\big)$ is a seed in $\bbC(\mcG)$, for the two reasons given below.
Then $\Sigma_\gamma(\bfi)$ is also a seed if $h_\gamma$ is birational, as it is obtained from $\Sigma^\Phi_\gamma(\bfi)$ by composing all functions with $h_\gamma^{-1}$.
First, the principal part of $D_\gamma(\bfi)$ is skew-symmetric.
Its entries join chambers at the level $x$ of an inversion letter $i_p$ and at the level $a^*$, with $d_x=d_{a^*}$.
Indeed, the root $\alpha_{\gamma_\bullet(a)}$ of $i_p$ is conjugate to $\alpha_x$ and $\gamma_\bullet$ and $*$ preserve root lengths; so the principal part of $\tB_\gamma(\bfi)$ is skew-symmetrizable with the symmetrizer of $\tB(\bfi^\gamma)$ (\cref{sec: words}).
Second, by \cref{lem: tails} the thawed chambers can be ordered so that each $\Phi_\gamma(t;\bfi)$ involves only earlier ones.
Putting the other chambers after them, each $\Phi_\gamma(k;\bfi)$ is~$\Delta(k;\bfi)$ times a monomial in earlier minors, so the $\Phi_\gamma(k;\bfi)$ are algebraically independent generators of~$\bbC(\mcG)$, as the $\Delta(k;\bfi)$ are (\cref{sec: double words}).

The map $h_\gamma$ is birational for $SL_N$ by \cref{thm: GSV h} and for $\Sp$ by \cref{thm: fold}(c).
It is also birational for $\gamma=(0,0)$ in every type.
In this case $h_\gamma$ is the identity and there are no thawed chambers, so $\Sigma_{(0,0)}(\bfi)$ is the seed $\Sigma(\bfi)$ of \cite[Theorem~2.10]{BFZ05}.

\begin{example}
  \label{ex: Sp6}
  In type $C_n$, $w_0=-1$, so $*$ is the identity, and the only long simple root is $\alpha_n$, so a BD map, being an isometry and nilpotent, only involves $\alpha_1,\dots,\alpha_{n-1}$.
  Thus $Sp_4$ only has $\gamma=(0,0)$, and $Sp_6$ is the smallest case with non-trivial BD pairs.
  Let $\mcG=Sp_6$, with $\alpha_3$ long, $\gamma=(\{1\to2\},0)$, and let $\bfi$ be the word of \cref{fig: Sp6 word seed},
  \[
    \bfi=(-3,-2,-3,-1,-2,-3,-1,-2,-1,\,1,2,1,3,2,1,3,2,3).
  \]
  The negative letter $i_7=-1$ has the root $\alpha_2=\alpha_{\gamma_{\rr}(1)}$, so it is the inversion letter.
  The thawed chamber is~$t^\rr_1=c_1(18)$, with $\Phi_\gamma(t^\rr_1;\bfi)=\Delta_{-\varpi_1,\varpi_1}$, and $\bfi^\gamma$ is $\bfi$ followed by the virtual letter $-1$.
  The non-zero entries of the column $t^\rr_1$ of $\tB_\gamma(\bfi)$ are
  \[
    b_{c_1(6),t^\rr_1}=b_{c_1(14),t^\rr_1}=1,
    \qquad
    b_{c_1(7),t^\rr_1}=b_{c_2(14),t^\rr_1}=-1,
  \]
  those at $c_1(6)$ and $c_1(7)$ from the dipole and the others from the virtual letter.
  The other entries of $\tB_\gamma(\bfi)$ not in $\tB(\bfi)$ are the dipole entries $b_{t^\rr_1,c_1(6)}=-1$ and $b_{t^\rr_1,c_1(7)}=1$.

\afterpage{%
  \begin{landscape}
    \centerfloatpage
    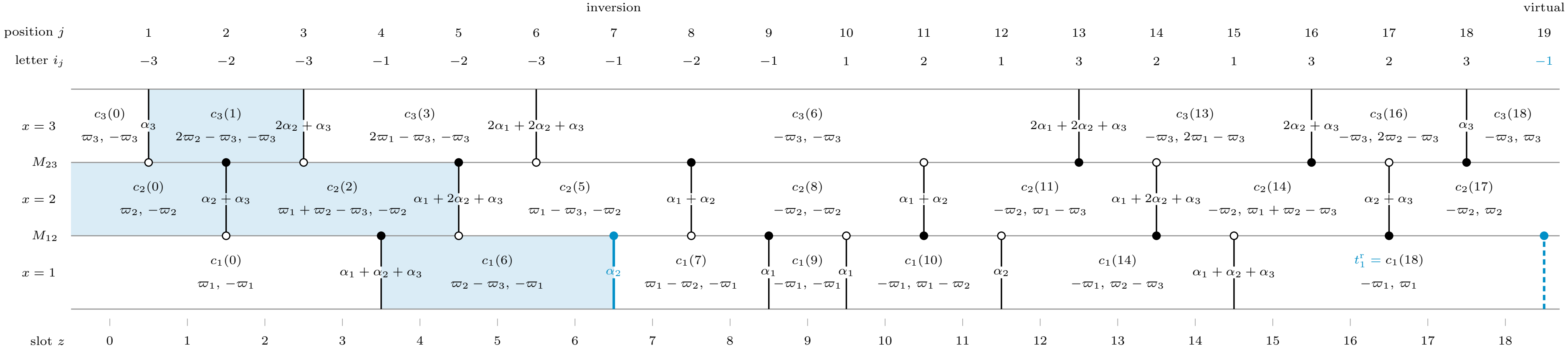
\begin{figure}[p]
      \centering
      \begingroup
        \colorlet{ink}{black}
        \colorlet{dipole}{accent}
        \colorlet{virtual}{accent}
        \colorlet{muted}{black}
        \resizebox{\linewidth}{!}{%
          \begin{tikzpicture}[x=1.18cm,y=1.45cm]
            \path[use as bounding box] (-1.1,.35) rectangle (24.96,5.25);
            \fill[dipole!12] (0,2) rectangle (2.6,3);
            \fill[dipole!12] (1.3,3) rectangle (3.9,4);
            \fill[dipole!12] (2.6,2) rectangle (6.5,3);
            \fill[dipole!12] (5.2,1) rectangle (9.1,2);
            \foreach \y in {1,2,3,4} \draw[ink!42,line width=.6pt] (0,\y)--(24.96,\y);
            \node[anchor=east,font=\scriptsize,text=black,inner sep=2pt] at (-0.15,2) {$M_{12}$};
            \node[anchor=east,font=\scriptsize,text=black,inner sep=2pt] at (-0.15,3) {$M_{23}$};
            \draw[ink,line width=0.8pt] (1.3,3)--(1.3,3.42) (1.3,3.58)--(1.3,4);
            \draw[draw=ink,fill=white,line width=.65pt] (1.3,3) circle[radius=2.1pt];
            \node[font=\scriptsize,text=ink,inner sep=1.2pt] at (1.3,3.5) {$\alpha_3$};
            \node[font=\scriptsize,text=ink] at (1.3,4.38) {$-3$};
            \node[font=\scriptsize,text=muted] at (1.3,4.77) {$1$};
            \draw[ink,line width=0.8pt] (2.6,2)--(2.6,2.42) (2.6,2.58)--(2.6,3);
            \draw[draw=ink,fill=white,line width=.65pt] (2.6,2) circle[radius=2.1pt];
            \draw[draw=ink,fill=ink,line width=.65pt] (2.6,3) circle[radius=2.1pt];
            \node[font=\scriptsize,text=ink,inner sep=1.2pt] at (2.6,2.5) {$\alpha_2+\alpha_3$};
            \node[font=\scriptsize,text=ink] at (2.6,4.38) {$-2$};
            \node[font=\scriptsize,text=muted] at (2.6,4.77) {$2$};
            \draw[ink,line width=0.8pt] (3.9,3)--(3.9,3.42) (3.9,3.58)--(3.9,4);
            \draw[draw=ink,fill=white,line width=.65pt] (3.9,3) circle[radius=2.1pt];
            \node[font=\scriptsize,text=ink,inner sep=1.2pt] at (3.9,3.5) {$2\alpha_2+\alpha_3$};
            \node[font=\scriptsize,text=ink] at (3.9,4.38) {$-3$};
            \node[font=\scriptsize,text=muted] at (3.9,4.77) {$3$};
            \draw[ink,line width=0.8pt] (5.2,1)--(5.2,1.42) (5.2,1.58)--(5.2,2);
            \draw[draw=ink,fill=ink,line width=.65pt] (5.2,2) circle[radius=2.1pt];
            \node[font=\scriptsize,text=ink,inner sep=1.2pt] at (5.2,1.5) {$\alpha_1+\alpha_2+\alpha_3$};
            \node[font=\scriptsize,text=ink] at (5.2,4.38) {$-1$};
            \node[font=\scriptsize,text=muted] at (5.2,4.77) {$4$};
            \draw[ink,line width=0.8pt] (6.5,2)--(6.5,2.42) (6.5,2.58)--(6.5,3);
            \draw[draw=ink,fill=white,line width=.65pt] (6.5,2) circle[radius=2.1pt];
            \draw[draw=ink,fill=ink,line width=.65pt] (6.5,3) circle[radius=2.1pt];
            \node[font=\scriptsize,text=ink,inner sep=1.2pt] at (6.5,2.5) {$\alpha_1+2\alpha_2+\alpha_3$};
            \node[font=\scriptsize,text=ink] at (6.5,4.38) {$-2$};
            \node[font=\scriptsize,text=muted] at (6.5,4.77) {$5$};
            \draw[ink,line width=0.8pt] (7.8,3)--(7.8,3.42) (7.8,3.58)--(7.8,4);
            \draw[draw=ink,fill=white,line width=.65pt] (7.8,3) circle[radius=2.1pt];
            \node[font=\scriptsize,text=ink,inner sep=1.2pt] at (7.8,3.5) {$2\alpha_1+2\alpha_2+\alpha_3$};
            \node[font=\scriptsize,text=ink] at (7.8,4.38) {$-3$};
            \node[font=\scriptsize,text=muted] at (7.8,4.77) {$6$};
            \draw[dipole,line width=1.4pt] (9.1,1)--(9.1,1.42) (9.1,1.58)--(9.1,2);
            \draw[draw=dipole,fill=dipole,line width=.65pt] (9.1,2) circle[radius=2.1pt];
            \node[font=\scriptsize,text=dipole,inner sep=1.2pt] at (9.1,1.5) {$\alpha_2$};
            \node[font=\scriptsize,text=black] at (9.1,4.38) {$-1$};
            \node[font=\scriptsize,text=muted] at (9.1,4.77) {$7$};
            \draw[ink,line width=0.8pt] (10.4,2)--(10.4,2.42) (10.4,2.58)--(10.4,3);
            \draw[draw=ink,fill=white,line width=.65pt] (10.4,2) circle[radius=2.1pt];
            \draw[draw=ink,fill=ink,line width=.65pt] (10.4,3) circle[radius=2.1pt];
            \node[font=\scriptsize,text=ink,inner sep=1.2pt] at (10.4,2.5) {$\alpha_1+\alpha_2$};
            \node[font=\scriptsize,text=ink] at (10.4,4.38) {$-2$};
            \node[font=\scriptsize,text=muted] at (10.4,4.77) {$8$};
            \draw[ink,line width=0.8pt] (11.7,1)--(11.7,1.42) (11.7,1.58)--(11.7,2);
            \draw[draw=ink,fill=ink,line width=.65pt] (11.7,2) circle[radius=2.1pt];
            \node[font=\scriptsize,text=ink,inner sep=1.2pt] at (11.7,1.5) {$\alpha_1$};
            \node[font=\scriptsize,text=ink] at (11.7,4.38) {$-1$};
            \node[font=\scriptsize,text=muted] at (11.7,4.77) {$9$};
            \draw[ink,line width=0.8pt] (13,1)--(13,1.42) (13,1.58)--(13,2);
            \draw[draw=ink,fill=white,line width=.65pt] (13,2) circle[radius=2.1pt];
            \node[font=\scriptsize,text=ink,inner sep=1.2pt] at (13,1.5) {$\alpha_1$};
            \node[font=\scriptsize,text=ink] at (13,4.38) {$1$};
            \node[font=\scriptsize,text=muted] at (13,4.77) {$10$};
            \draw[ink,line width=0.8pt] (14.3,2)--(14.3,2.42) (14.3,2.58)--(14.3,3);
            \draw[draw=ink,fill=white,line width=.65pt] (14.3,3) circle[radius=2.1pt];
            \draw[draw=ink,fill=ink,line width=.65pt] (14.3,2) circle[radius=2.1pt];
            \node[font=\scriptsize,text=ink,inner sep=1.2pt] at (14.3,2.5) {$\alpha_1+\alpha_2$};
            \node[font=\scriptsize,text=ink] at (14.3,4.38) {$2$};
            \node[font=\scriptsize,text=muted] at (14.3,4.77) {$11$};
            \draw[ink,line width=0.8pt] (15.6,1)--(15.6,1.42) (15.6,1.58)--(15.6,2);
            \draw[draw=ink,fill=white,line width=.65pt] (15.6,2) circle[radius=2.1pt];
            \node[font=\scriptsize,text=ink,inner sep=1.2pt] at (15.6,1.5) {$\alpha_2$};
            \node[font=\scriptsize,text=ink] at (15.6,4.38) {$1$};
            \node[font=\scriptsize,text=muted] at (15.6,4.77) {$12$};
            \draw[ink,line width=0.8pt] (16.9,3)--(16.9,3.42) (16.9,3.58)--(16.9,4);
            \draw[draw=ink,fill=ink,line width=.65pt] (16.9,3) circle[radius=2.1pt];
            \node[font=\scriptsize,text=ink,inner sep=1.2pt] at (16.9,3.5) {$2\alpha_1+2\alpha_2+\alpha_3$};
            \node[font=\scriptsize,text=ink] at (16.9,4.38) {$3$};
            \node[font=\scriptsize,text=muted] at (16.9,4.77) {$13$};
            \draw[ink,line width=0.8pt] (18.2,2)--(18.2,2.42) (18.2,2.58)--(18.2,3);
            \draw[draw=ink,fill=white,line width=.65pt] (18.2,3) circle[radius=2.1pt];
            \draw[draw=ink,fill=ink,line width=.65pt] (18.2,2) circle[radius=2.1pt];
            \node[font=\scriptsize,text=ink,inner sep=1.2pt] at (18.2,2.5) {$\alpha_1+2\alpha_2+\alpha_3$};
            \node[font=\scriptsize,text=ink] at (18.2,4.38) {$2$};
            \node[font=\scriptsize,text=muted] at (18.2,4.77) {$14$};
            \draw[ink,line width=0.8pt] (19.5,1)--(19.5,1.42) (19.5,1.58)--(19.5,2);
            \draw[draw=ink,fill=white,line width=.65pt] (19.5,2) circle[radius=2.1pt];
            \node[font=\scriptsize,text=ink,inner sep=1.2pt] at (19.5,1.5) {$\alpha_1+\alpha_2+\alpha_3$};
            \node[font=\scriptsize,text=ink] at (19.5,4.38) {$1$};
            \node[font=\scriptsize,text=muted] at (19.5,4.77) {$15$};
            \draw[ink,line width=0.8pt] (20.8,3)--(20.8,3.42) (20.8,3.58)--(20.8,4);
            \draw[draw=ink,fill=ink,line width=.65pt] (20.8,3) circle[radius=2.1pt];
            \node[font=\scriptsize,text=ink,inner sep=1.2pt] at (20.8,3.5) {$2\alpha_2+\alpha_3$};
            \node[font=\scriptsize,text=ink] at (20.8,4.38) {$3$};
            \node[font=\scriptsize,text=muted] at (20.8,4.77) {$16$};
            \draw[ink,line width=0.8pt] (22.1,2)--(22.1,2.42) (22.1,2.58)--(22.1,3);
            \draw[draw=ink,fill=white,line width=.65pt] (22.1,3) circle[radius=2.1pt];
            \draw[draw=ink,fill=ink,line width=.65pt] (22.1,2) circle[radius=2.1pt];
            \node[font=\scriptsize,text=ink,inner sep=1.2pt] at (22.1,2.5) {$\alpha_2+\alpha_3$};
            \node[font=\scriptsize,text=ink] at (22.1,4.38) {$2$};
            \node[font=\scriptsize,text=muted] at (22.1,4.77) {$17$};
            \draw[ink,line width=0.8pt] (23.4,3)--(23.4,3.42) (23.4,3.58)--(23.4,4);
            \draw[draw=ink,fill=ink,line width=.65pt] (23.4,3) circle[radius=2.1pt];
            \node[font=\scriptsize,text=ink,inner sep=1.2pt] at (23.4,3.5) {$\alpha_3$};
            \node[font=\scriptsize,text=ink] at (23.4,4.38) {$3$};
            \node[font=\scriptsize,text=muted] at (23.4,4.77) {$18$};
            \draw[virtual,line width=1.4pt,densely dashed] (24.7,1)--(24.7,2);
            \draw[draw=virtual,fill=virtual,line width=.65pt] (24.7,2) circle[radius=2.1pt];
            \node[font=\scriptsize,text=accent] at (24.7,4.38) {$-1$};
            \node[font=\scriptsize,text=muted] at (24.7,4.77) {$19$};

            \node[anchor=east,font=\scriptsize,text=muted] at (0,4.77) {position $j$};
            \node[anchor=east,font=\scriptsize,text=muted] at (0,4.38) {letter $i_j$};
            \node[font=\scriptsize,text=black] at (9.1,5.12) {inversion};
            \node[font=\scriptsize,text=black] at (24.7,5.12) {virtual};
            \foreach \x in {1,2,3} \node[anchor=east,font=\scriptsize] at (-0.15,\x+.5) {$x=\x$};
            \foreach \s/\slotx in {0/0.65,1/1.95,2/3.25,3/4.55,4/5.85,5/7.15,6/8.45,7/9.75,8/11.05,9/12.35,10/13.65,11/14.95,12/16.25,13/17.55,14/18.85,15/20.15,16/21.45,17/22.75,18/24.05} {
              \draw[ink!35] (\slotx,.77)--(\slotx,.87);
              \node[font=\scriptsize,text=muted] at (\slotx,.57) {$\s$};
            }
            \node[anchor=east,font=\scriptsize,text=muted] at (0,.57) {slot $z$};
            \node[font=\scriptsize,text=ink,inner sep=0pt] at (2.6,1.66) {$c_1(0)$};
            \node[font=\scriptsize,text=ink,inner sep=0pt] at (2.6,1.34) {$\varpi_1,\,-\varpi_1$};
            \node[font=\scriptsize,text=ink,inner sep=0pt] at (7.15,1.66) {$c_1(6)$};
            \node[font=\scriptsize,text=ink,inner sep=0pt] at (7.15,1.34) {$\varpi_2-\varpi_3,\,-\varpi_1$};
            \node[font=\scriptsize,text=ink,inner sep=0pt] at (10.4,1.66) {$c_1(7)$};
            \node[font=\scriptsize,text=ink,inner sep=0pt] at (10.4,1.34) {$\varpi_1-\varpi_2,\,-\varpi_1$};
            \node[font=\scriptsize,text=ink,inner sep=0pt] at (12.35,1.66) {$c_1(9)$};
            \node[font=\scriptsize,text=ink,inner sep=0pt] at (12.35,1.34) {$-\varpi_1,\,-\varpi_1$};
            \node[font=\scriptsize,text=ink,inner sep=0pt] at (14.3,1.66) {$c_1(10)$};
            \node[font=\scriptsize,text=ink,inner sep=0pt] at (14.3,1.34) {$-\varpi_1,\,\varpi_1-\varpi_2$};
            \node[font=\scriptsize,text=ink,inner sep=0pt] at (17.55,1.66) {$c_1(14)$};
            \node[font=\scriptsize,text=ink,inner sep=0pt] at (17.55,1.34) {$-\varpi_1,\,\varpi_2-\varpi_3$};
            \node[font=\scriptsize,text=black,inner sep=0pt] at (22.1,1.66) {$\textcolor{accent}{t_1^{\mathrm r}}\mathrel{\textcolor{accent}{=}}c_1(18)$};
            \node[font=\scriptsize,text=ink,inner sep=0pt] at (22.1,1.34) {$-\varpi_1,\,\varpi_1$};
            \node[font=\scriptsize,text=ink,inner sep=0pt] at (1.3,2.66) {$c_2(0)$};
            \node[font=\scriptsize,text=ink,inner sep=0pt] at (1.3,2.34) {$\varpi_2,\,-\varpi_2$};
            \node[font=\scriptsize,text=ink,inner sep=0pt] at (4.55,2.66) {$c_2(2)$};
            \node[font=\scriptsize,text=ink,inner sep=0pt] at (4.55,2.34) {$\varpi_1+\varpi_2-\varpi_3,\,-\varpi_2$};
            \node[font=\scriptsize,text=ink,inner sep=0pt] at (8.45,2.66) {$c_2(5)$};
            \node[font=\scriptsize,text=ink,inner sep=0pt] at (8.45,2.34) {$\varpi_1-\varpi_3,\,-\varpi_2$};
            \node[font=\scriptsize,text=ink,inner sep=0pt] at (12.35,2.66) {$c_2(8)$};
            \node[font=\scriptsize,text=ink,inner sep=0pt] at (12.35,2.34) {$-\varpi_2,\,-\varpi_2$};
            \node[font=\scriptsize,text=ink,inner sep=0pt] at (16.25,2.66) {$c_2(11)$};
            \node[font=\scriptsize,text=ink,inner sep=0pt] at (16.25,2.34) {$-\varpi_2,\,\varpi_1-\varpi_3$};
            \node[font=\scriptsize,text=ink,inner sep=0pt] at (20.15,2.66) {$c_2(14)$};
            \node[font=\scriptsize,text=ink,inner sep=0pt] at (20.15,2.34) {$-\varpi_2,\,\varpi_1+\varpi_2-\varpi_3$};
            \node[font=\scriptsize,text=ink,inner sep=0pt] at (23.53,2.66) {$c_2(17)$};
            \node[font=\scriptsize,text=ink,inner sep=0pt] at (23.53,2.34) {$-\varpi_2,\,\varpi_2$};
            \node[font=\scriptsize,text=ink,inner sep=0pt] at (0.65,3.66) {$c_3(0)$};
            \node[font=\scriptsize,text=ink,inner sep=0pt] at (0.65,3.34) {$\varpi_3,\,-\varpi_3$};
            \node[font=\scriptsize,text=ink,inner sep=0pt] at (2.6,3.66) {$c_3(1)$};
            \node[font=\scriptsize,text=ink,inner sep=0pt] at (2.6,3.34) {$2\varpi_2-\varpi_3,\,-\varpi_3$};
            \node[font=\scriptsize,text=ink,inner sep=0pt] at (5.85,3.66) {$c_3(3)$};
            \node[font=\scriptsize,text=ink,inner sep=0pt] at (5.85,3.34) {$2\varpi_1-\varpi_3,\,-\varpi_3$};
            \node[font=\scriptsize,text=ink,inner sep=0pt] at (12.35,3.66) {$c_3(6)$};
            \node[font=\scriptsize,text=ink,inner sep=0pt] at (12.35,3.34) {$-\varpi_3,\,-\varpi_3$};
            \node[font=\scriptsize,text=ink,inner sep=0pt] at (18.85,3.66) {$c_3(13)$};
            \node[font=\scriptsize,text=ink,inner sep=0pt] at (18.85,3.34) {$-\varpi_3,\,2\varpi_1-\varpi_3$};
            \node[font=\scriptsize,text=ink,inner sep=0pt] at (22.1,3.66) {$c_3(16)$};
            \node[font=\scriptsize,text=ink,inner sep=0pt] at (22.1,3.34) {$-\varpi_3,\,2\varpi_2-\varpi_3$};
            \node[font=\scriptsize,text=ink,inner sep=0pt] at (24.18,3.66) {$c_3(18)$};
            \node[font=\scriptsize,text=ink,inner sep=0pt] at (24.18,3.34) {$-\varpi_3,\,\varpi_3$};

          \end{tikzpicture}%
        }
      \endgroup
      \par\medskip
      \begin{minipage}[c]{.63\linewidth}
        \setlength{\captionindent}{0pt}
        \setlength{\abovecaptionskip}{0pt}
        \caption[The BD word seed for the $Sp_6$ example]{%
          The word $\bfi=(-3,-2,-3,-1,-2,-3,-1,-2,-1,1,2,1,3,2,1,3,2,3)$ and the extended word $\bfi^\gamma=(-3,-2,-3,-1,-2,-3,-1,-2,-1,1,2,1,3,2,1,3,2,3,\textcolor{accent}{-1})$ for the group $\mcG=Sp_6$ and the BD pair $\gamma=(\{1\to2\},0)$ in \cref{ex: Sp6}.
          The vertical segment at position $j$ represents the letter~$i_j$ at level $x=|i_j|$ and, for $j\in[1,18]$, is labeled by its root; below each chamber name are its row and column weights $\lambda^\rr,\lambda^\cc$.
          On the horizontal line labeled $M_{xy}$ bullets and circles depict the values~$\varphi_{xy}(j)=1$ and $\varphi_{xy}(j)=-1$, respectively.
          The marked positions give $M_{12}$ and $M_{23}$ for the original word, with an additional mark at position~$19$ for the virtual letter.
          Blue identifies the inversion letter $i_7=-1$, at the position $p^\rr_1=7$, whose root is $\alpha_2$, the thawed chamber $t^\rr_1=c_1(18)$, and the appended virtual letter $-1$ (dashed).
          The shaded chambers are exactly those with $e^\rr_1(k)=[\lambda^\rr_k(\alpha_2^\vee)]_+>0$: the exponent is $1$ at $c_1(6)$, $c_2(0)$ and $c_2(2)$, and $2$ at $c_3(1)$.
          Their minors acquire the corresponding power of~$\Phi_\gamma(t^\rr_1;\bfi)=\Delta_{-\varpi_1,\varpi_1}$.
          The exchange matrix on the right has rows indexed by all chambers and columns by the exchangeable chambers; dots denote zero entries, and the rules separate levels.
          The blue entries in the column of $t^\rr_1=c_1(18)$ come from the virtual letter, and the four red entries are the non-zero entries of the dipole matrix $D_\gamma(\bfi)$.}
        \label{fig: Sp6 word seed}
      \end{minipage}\hfill
      \begin{minipage}[c]{.34\linewidth}
        \begingroup
          \centering
          {\small $\widetilde B_\gamma(\mathbf i)$}\par\smallskip
          \setlength{\arraycolsep}{1.5pt}
          \renewcommand{\arraystretch}{.75}
          \resizebox{\linewidth}{!}{%
            $\begin{array}{r|rrrrrr|rrrrr|rrrrr}
              & \rotatebox{270}{$c_1(6)$} & \rotatebox{270}{$c_1(7)$} & \rotatebox{270}{$c_1(9)$} & \rotatebox{270}{$c_1(10)$} & \rotatebox{270}{$c_1(14)$} & \rotatebox{270}{$\textcolor{accent}{c_1(18)}$} & \rotatebox{270}{$c_2(2)$} & \rotatebox{270}{$c_2(5)$} & \rotatebox{270}{$c_2(8)$} & \rotatebox{270}{$c_2(11)$} & \rotatebox{270}{$c_2(14)$} & \rotatebox{270}{$c_3(1)$} & \rotatebox{270}{$c_3(3)$} & \rotatebox{270}{$c_3(6)$} & \rotatebox{270}{$c_3(13)$} & \rotatebox{270}{$c_3(16)$} \\ \hline
              c_1(0) & -1 & \cdot & \cdot & \cdot & \cdot & \cdot & 1 & \cdot & \cdot & \cdot & \cdot & \cdot & \cdot & \cdot & \cdot & \cdot \\
              c_1(6) & \cdot & -1 & \cdot & \cdot & \cdot & \textcolor{red}{1} & -1 & 1 & \cdot & \cdot & \cdot & \cdot & \cdot & \cdot & \cdot & \cdot \\
              c_1(7) & 1 & \cdot & -1 & \cdot & \cdot & \textcolor{red}{-1} & \cdot & -1 & 1 & \cdot & \cdot & \cdot & \cdot & \cdot & \cdot & \cdot \\
              c_1(9) & \cdot & 1 & \cdot & 1 & \cdot & \cdot & \cdot & \cdot & -1 & \cdot & \cdot & \cdot & \cdot & \cdot & \cdot & \cdot \\
              c_1(10) & \cdot & \cdot & -1 & \cdot & 1 & \cdot & \cdot & \cdot & 1 & -1 & \cdot & \cdot & \cdot & \cdot & \cdot & \cdot \\
              c_1(14) & \cdot & \cdot & \cdot & -1 & \cdot & \textcolor{accent}{1} & \cdot & \cdot & \cdot & 1 & -1 & \cdot & \cdot & \cdot & \cdot & \cdot \\
              c_1(18) & \textcolor{red}{-1} & \textcolor{red}{1} & \cdot & \cdot & -1 & \cdot & \cdot & \cdot & \cdot & \cdot & 1 & \cdot & \cdot & \cdot & \cdot & \cdot \\
              \hline
              c_2(0) & \cdot & \cdot & \cdot & \cdot & \cdot & \cdot & -1 & \cdot & \cdot & \cdot & \cdot & 2 & \cdot & \cdot & \cdot & \cdot \\
              c_2(2) & 1 & \cdot & \cdot & \cdot & \cdot & \cdot & \cdot & -1 & \cdot & \cdot & \cdot & -2 & 2 & \cdot & \cdot & \cdot \\
              c_2(5) & -1 & 1 & \cdot & \cdot & \cdot & \cdot & 1 & \cdot & -1 & \cdot & \cdot & \cdot & -2 & 2 & \cdot & \cdot \\
              c_2(8) & \cdot & -1 & 1 & -1 & \cdot & \cdot & \cdot & 1 & \cdot & 1 & \cdot & \cdot & \cdot & -2 & \cdot & \cdot \\
              c_2(11) & \cdot & \cdot & \cdot & 1 & -1 & \cdot & \cdot & \cdot & -1 & \cdot & 1 & \cdot & \cdot & 2 & -2 & \cdot \\
              c_2(14) & \cdot & \cdot & \cdot & \cdot & 1 & \textcolor{accent}{-1} & \cdot & \cdot & \cdot & -1 & \cdot & \cdot & \cdot & \cdot & 2 & -2 \\
              c_2(17) & \cdot & \cdot & \cdot & \cdot & \cdot & \cdot & \cdot & \cdot & \cdot & \cdot & -1 & \cdot & \cdot & \cdot & \cdot & 2 \\
              \hline
              c_3(0) & \cdot & \cdot & \cdot & \cdot & \cdot & \cdot & \cdot & \cdot & \cdot & \cdot & \cdot & -1 & \cdot & \cdot & \cdot & \cdot \\
              c_3(1) & \cdot & \cdot & \cdot & \cdot & \cdot & \cdot & 1 & \cdot & \cdot & \cdot & \cdot & \cdot & -1 & \cdot & \cdot & \cdot \\
              c_3(3) & \cdot & \cdot & \cdot & \cdot & \cdot & \cdot & -1 & 1 & \cdot & \cdot & \cdot & 1 & \cdot & -1 & \cdot & \cdot \\
              c_3(6) & \cdot & \cdot & \cdot & \cdot & \cdot & \cdot & \cdot & -1 & 1 & -1 & \cdot & \cdot & 1 & \cdot & 1 & \cdot \\
              c_3(13) & \cdot & \cdot & \cdot & \cdot & \cdot & \cdot & \cdot & \cdot & \cdot & 1 & -1 & \cdot & \cdot & -1 & \cdot & 1 \\
              c_3(16) & \cdot & \cdot & \cdot & \cdot & \cdot & \cdot & \cdot & \cdot & \cdot & \cdot & 1 & \cdot & \cdot & \cdot & -1 & \cdot \\
              c_3(18) & \cdot & \cdot & \cdot & \cdot & \cdot & \cdot & \cdot & \cdot & \cdot & \cdot & \cdot & \cdot & \cdot & \cdot & \cdot & -1 \\
            \end{array}$%
          }
        \par\endgroup
      \end{minipage}
    \end{figure}
  \end{landscape}
}

  The row weights pair positively with $\alpha_2^\vee$ only at $c_1(6)$, $c_2(0)$, $c_2(2)$ and $c_3(1)$, so $\Phi_\gamma(k;\bfi)=\Delta(k;\bfi)$ except for
  \[
    \begin{aligned}
      \Phi_\gamma(c_1(6);\bfi)&=\Delta_{\varpi_2-\varpi_3,-\varpi_1}\Delta_{-\varpi_1,\varpi_1},
      &
      \Phi_\gamma(c_2(0);\bfi)&=\Delta_{\varpi_2,-\varpi_2}\Delta_{-\varpi_1,\varpi_1},
      \\
      \Phi_\gamma(c_2(2);\bfi)&=\Delta_{\varpi_1+\varpi_2-\varpi_3,-\varpi_2}\Delta_{-\varpi_1,\varpi_1},
      &
      \Phi_\gamma(c_3(1);\bfi)&=\Delta_{2\varpi_2-\varpi_3,-\varpi_3}\Delta_{-\varpi_1,\varpi_1}^2,
    \end{aligned}
  \]
  the last exponent being $(2\varpi_2-\varpi_3)(\alpha_2^\vee)=2$.
  For the aperiodic pair $(\{1\to2\},\{1\to2\})$ the column root contributes analogously with $t^\cc_1=c_1(0)$, $\Phi_\gamma(c_3(17);\bfi)=\Delta_{-\varpi_3,2\varpi_2-\varpi_3}\Delta_{\varpi_1,-\varpi_1}^2$, and $\Phi_\gamma(t^\cc_1;\bfi)=\Delta_{\varpi_1,-\varpi_1}$, while $(\{1\to2\},\{2\to1\})$ is not aperiodic.
\end{example}

\begin{lemma}
  \label{lem: local entries}
  \begin{enumerate}[(a)]
    \item For $k\in\mcI$ and $\ell$ $\bfi$-exchangeable, the entries of $\tB(\bfi^\gamma)_\mcI$ and of $\tB(\bfi)$ at $(k,\ell)$ coincide; the entries of $\tB(\bfi^\gamma)_\mcI$ between chambers at distinct levels $x,y$ with $a_{xy}=0$ vanish; and if the Cartan matrix is symmetric, the principal part of $\tB_\gamma(\bfi)$ is skew-symmetric.
    \item Every non-zero entry of $D_\gamma(\bfi)$ involves a thawed chamber; so the submatrices of $\tB_\gamma(\bfi)$ and of $\tB(\bfi)$ with rows and columns the $\bfi$-exchangeable chambers coincide.
  \end{enumerate}
\end{lemma}
\begin{proof}
  (a) The entries in the column $\ell$ come at equal levels from the two letters bounding $\ell$, and at distinct levels from pairs $p<p'$ with $p\in\ell$.
  As $\ell$ contains no end slot of $\bfi$, these letters are letters of $\bfi$, and so is $p'$, which does not exceed the letter bounding $\ell$ on the right.
  The other claims follow directly from \eqref{eq: B entries} and the antisymmetry of the dipole entries.
  Part (b) follows directly from \cref{def: BD word seed}(c).
\end{proof}

For a word $\bfi$ let $\bfi^T$ be the word obtained by reversing the order of the letters and changing their signs.
Its slots correspond to those of $\bfi$ in reverse order.
For a chamber $k$ of $\bfi$ let $k^T$ be the chamber of $\bfi^T$ at the same level formed by the slots corresponding to those of $k$.
For a double reduced word, indexed by $[1,L]$, $\bfi^T=(-i_L,\dots,-i_1)$: the letter $i_j$ becomes the letter $-i_j$ at the position $L+1-j$, and~$c_x(z)^T=c_x(L-z)$.
For a BD pair $\gamma$ let $\gamma^T:=(\gamma_{\cc},\gamma_{\rr})$; by \cref{lem: mirror} below, $\bfi\mapsto\bfi^T$ and $\gamma\mapsto\gamma^T$ give the transpose on functions.

\begin{lemma}
  \label{lem: mirror}
  \begin{enumerate}[(a)]
    \item For every word $\bfi$, $k\mapsto k^T$ is a bijection $C(\bfi)\to C(\bfi^T)$ which preserves levels and exchangeability, and $\tB(\bfi^T)$ is obtained from $\tB(\bfi)$ by relabeling rows and columns along it.
    \item For a double reduced word $\bfi$, under these identifications the inversion letters, thawed chambers and exponents of the row roots of $\gamma$ in $\bfi$ become those of the column roots of $\gamma^T$ in $\bfi^T$, and conversely, and $(\bfi^\gamma)^T=(\bfi^T)^{\gamma^T}$ up to translation.
      The matrix $\tB_{\gamma^T}(\bfi^T)$ is obtained from $\tB_\gamma(\bfi)$ by relabeling rows and columns along $k\mapsto k^T$, and $\Phi_{\gamma^T}(k^T;\bfi^T)=\Phi_\gamma(k;\bfi)\circ T$ for every chamber $k$, $T$ being the transpose.
    \item If $h_\gamma$ is birational, so is $h_{\gamma^T}$, and $x_{\gamma^T}(k^T;\bfi^T)=x_\gamma(k;\bfi)\circ T$.
  \end{enumerate}
\end{lemma}
\begin{proof}
  We first compare the exchange matrices under reversal, then follow the weights, inversion letters and recursive functions through the same correspondence.

  (a) Index $\bfi^T$ by $-J$, with the letter $-i_j$ at the position $-j$ and the slot $-1-z$ corresponding to $z$.
  The letter $i_j$, at the level $x$ and separating $c_x(j-1)$ from $c_x(j)$, becomes the letter $-i_j$ at the position~$-j$, separating $c_x(j)^T$ from $c_x(j-1)^T$; so the entries at equal levels are preserved.
  Consecutive positions~$p<p'$ at the levels $x,y$ become $-p'<-p$, with $\varphi_{xy}$ changing sign, and $c_x(p)^T=c_x(p'-1)^T$ is the chamber of~$\bfi^T$ containing the slot $-p'$.
  So the entries at distinct levels are preserved.

  (b) The negative letters of $\bfi^T$ at positions $\le L-z$ are the negatives of the positive letters of $\bfi$ at positions $>z$, in the opposite order, so $u_{\le L-z}(\bfi^T)=v_{>z}(\bfi)$, and similarly $v_{>L-z}(\bfi^T)=u_{\le z}(\bfi)$.
  So the weights of $k^T$ are those of $k$ in the opposite order, $\Delta(k^T;\bfi^T)=\Delta(k;\bfi)\circ T$ by \cref{sec: Lie}, and $e^\rr_a(k)$ for~$\gamma$ equals $e^\cc_a(k^T)$ for $\gamma^T$, and conversely.
  The inversion letters correspond, as $u_{\le j-1}(\bfi)=v_{>L+1-j}(\bfi^T)$, and so do the thawed chambers, as $k\mapsto k^T$ swaps initial and final chambers.
  The virtual letters of $\bfi^\gamma$ go to those of $(\bfi^T)^{\gamma^T}$, in the required order, which gives $(\bfi^\gamma)^T=(\bfi^T)^{\gamma^T}$ up to translation, compatibly with the identifications of chambers.
  Finally, for a position $p$ of $\bfi$ and $p'=L+1-p$, $c_x(p-1)^T=c_x(p')$ and~$c_x(p)^T=c_x(p'-1)$, so $D_{\gamma^T}(\bfi^T)$ is $D_\gamma(\bfi)$ relabeled by $k\mapsto k^T$.
  With (a) applied to $\bfi^\gamma$, this gives the claim on $\tB$, and \eqref{eq: Phi} gives the claim on $\Phi$. (c) follows from (b) and \cref{lem: h transpose}.
\end{proof}

\subsection{Moves}
\label{sec: moves}
We prove that moves of words act on the seeds $\Sigma^\Phi_\gamma(\bfi)$ by mutations, by tracking the virtual letters, dipoles and monomial factors.
In particular, these seeds are mutation equivalent for all double reduced words.

\begin{theorem}
  \label{thm: moves}
  Let $\gamma$ be an aperiodic BD pair of $\mcG$, and let $\bfi'$ be obtained from a double reduced word $\bfi$ by a move.
  Along the correspondence of chambers:
  \begin{enumerate}[(a)]
    \item if the move is a commutation, then $\Sigma^\Phi_\gamma(\bfi')=\Sigma^\Phi_\gamma(\bfi)$;
    \item if it is a sign swap or a braid move of length $3$, then $\Sigma^\Phi_\gamma(\bfi')=\mu_k(\Sigma^\Phi_\gamma(\bfi))$, $k$ being the replaced chamber;
    \item if it is a braid move of length $4$, then $\Sigma^\Phi_\gamma(\bfi')=\mu_k\mu_\ell\mu_k(\Sigma^\Phi_\gamma(\bfi))$, where $k$ and $\ell$ are the replaced chambers at the levels of the long and of the short root, and the replaced chambers of $\bfi'$ correspond to them by level.
  \end{enumerate}
  If $h_\gamma$ is birational, the same holds for $\Sigma_\gamma(\bfi)$ and $\Sigma_\gamma(\bfi')$.
\end{theorem}
The last claim holds because composing with $h_\gamma^{-1}$ commutes with mutations.

We first compare $\Sigma^\Phi_\gamma(\bfi)$ with $\Sigma(\bfi)$ for a single word (\cref{lem: exponents,cor: Fraser}).
We then prove the statements on the matrices directly and use this comparison to obtain those on the functions.

For a double reduced word $\bfi$, with $\lambda^\rr_k,\lambda^\cc_k$ the weights of the chambers $k$ and $b_{k\ell}$ the entries of $\tB(\bfi)$,
\begin{equation}
  \label{eq: homogeneity}
  \sum_k\lambda^\rr_k b_{k\ell}=0,
  \qquad
  \sum_k\lambda^\cc_k b_{k\ell}=0
\end{equation}
for every exchangeable chamber $\ell$.
This is \cite[Lemma~4.22]{GSV10} and it expresses the bihomogeneity of the exchange relations of \cite{BFZ05} for the left and right actions of $\mcH$.
We also use the second identity for arbitrary words, with $u_{\le z}$, $v_{>z}$ and the weights defined by the same formulas.
Indeed, let $\ell$ be at the level $y$, bounded by positions $p<p'$.
At a positive letter $i_j=y'$, $(v_{>j-1}-v_{>j})\varpi_{y'}=-v_{>j}\alpha_{y'}$ and~$(v_{>j-1}-v_{>j})\alpha_y=-a_{y'y}v_{>j}\alpha_{y'}$; at a negative letter both differences vanish.
Summing \eqref{eq: B entries} by parts and telescoping over $p<j<p'$ therefore gives
\[
  2\sum_k\lambda^\cc_k b_{k\ell}
  =\bigl(\varepsilon(i_p)v_{>p}-\varepsilon(i_{p'})v_{>p'-1}\bigr)
  \Bigl(\sum_x a_{xy}\varpi_x-\alpha_y\Bigr)=0.
\]

The following lemma follows immediately from the root-sequence property of reduced words recalled in \cref{sec: Lie}.
Geometrically, a reduced gallery from the dominant to the opposite Weyl chamber crosses each reflecting hyperplane exactly once.

\begin{lemma}
  \label{lem: sign}
  Let $\bfi$ be a double reduced word, $g\in[1,r]$, $i_{p_0}$ the positive letter with root $\alpha_g$, and $x$ its level.
  For $k\in C(\bfi)$, $\lambda^\cc_k(\alpha_g^\vee)\ge0$ if $k$ has a slot $\ge p_0$, and $\lambda^\cc_k(\alpha_g^\vee)\le0$ if $k$ has a slot $\le p_0-1$; moreover,
  \[
    \lambda^\cc_{c_x(p_0-1)}(\alpha_g^\vee)=-1,
    \qquad \lambda^\cc_{c_x(p_0)}(\alpha_g^\vee)=1.
  \]
  In particular, $\lambda^\cc_k(\alpha_g^\vee)=0$ for every chamber $k$ with slots on both sides of $p_0$.
\end{lemma}
\begin{proof}

  Let $\beta^\vee_z:=v_{>z}^{-1}(\alpha_g^\vee)$; then $\lambda^\cc_{c_y(z)}(\alpha_g^\vee)$ is the coefficient of $\alpha_y^\vee$ in $\beta^\vee_z$.
  A positive letter $i_j$ at the level~$y$ gives $\beta^\vee_{j-1}=s_y(\beta^\vee_j)$, which has the sign of $\beta^\vee_j$ unless $\beta^\vee_j=\pm\alpha_y^\vee$, i.e.\ the root $v_{>j}(\alpha_y)$ of $i_j$ is $\pm\alpha_g$.
  Roots of positive letters being positive and pairwise distinct, this only happens for $j=p_0$, with $\beta^\vee_{p_0}=\alpha_x^\vee$.
  Negative letters do not change $\beta^\vee_z$, and $\beta^\vee_L=\alpha_g^\vee$; so $\beta^\vee_z$ is positive for $z\ge p_0$ and negative for $z\le p_0-1$, and $\beta^\vee_{p_0-1}=-\alpha_x^\vee$.
\end{proof}
Applying \cref{lem: sign} to $\bfi^T$, whose column weights at $k^T$ are the row weights at $k$, gives the analogue for a negative letter $i_{p_0}$ with root $\alpha_g$.
Thus $\lambda^\rr_k(\alpha_g^\vee)\le0$ if $k$ has a slot $\ge p_0$, $\lambda^\rr_k(\alpha_g^\vee)\ge0$ if $k$ has a slot~$\le p_0-1$, and $\lambda^\rr_{c_x(p_0-1)}(\alpha_g^\vee)=1$, $\lambda^\rr_{c_x(p_0)}(\alpha_g^\vee)=-1$.

For a root $a$ of either kind we write $p_a$, $t_a$ and $e_a$ for $p^\rr_a$, $t^\rr_a$ and $e^\rr_a$, resp.\ $p^\cc_a$, $t^\cc_a$ and $e^\cc_a$.

\begin{lemma}
  \label{lem: exponents}
  Let $\bfi$ be a double reduced word and $a$ a row or column root, with inversion letter $i_{p_a}$ at the level $x$.
  If $b_{k\ell}$ are the entries of $\tB(\bfi)$ then, for every $\bfi$-exchangeable chamber $\ell$,
  \[
    \sum_{k\in\mcI}e_a(k)\,b_{k\ell}=
    \begin{cases}
      -\varepsilon(i_{p_a})&\text{if }\ell=c_x(p_a-1),\\
      \varepsilon(i_{p_a})&\text{if }\ell=c_x(p_a),\\
      0&\text{otherwise}.
    \end{cases}
  \]
\end{lemma}
\begin{proof}
  By \cref{lem: mirror}(b) we may assume that $a$ is a column root, so that $\varepsilon(i_{p_a})=1$; let $g:=\gamma_{\cc}(a)$ and~$p:=p_a$, so that $e_a(k)=[\lambda^\cc_k(\alpha_g^\vee)]_+$, and recall \cref{lem: sign}.
  An entry $b_{k\ell}$ vanishes unless $k$ and $\ell$ are at the same level and separated by one letter, or at distinct levels with a common slot.
  If all the slots of $\ell$ are $\ge p$, apply \eqref{eq: homogeneity} to the word $(i_p,\dots,i_L)$.
  Its slots are $[p-1,L]$, and its chambers are the intersections with $[p-1,L]$ of the chambers of $\bfi$ with a slot $\ge p-1$, with the same column weights and the same entries in the column $\ell$.
  Among them the only chamber with negative pairing is $c_x(p-1)$, whose pairing is $-1$; so $\sum_ke_a(k)b_{k\ell}=b_{c_x(p-1),\ell}$, which is $1$ if $\ell=c_x(p)$ and $0$ otherwise.
  If all the slots of $\ell$ are~$\le p-1$, a chamber $k$ with $e_a(k)\neq0$ has no slot in common with $\ell$, so $\sum_ke_a(k)b_{k\ell}$ is~$0$ unless $\ell=c_x(p-1)$, in which case it is $b_{c_x(p),c_x(p-1)}=-1$.
  Otherwise $\ell$ has the slots $p-1$ and $p$, and is at a level $y\neq x$.
  Apply \eqref{eq: homogeneity} to the word $\hat{\boldsymbol{\imath}}:=(-y,i_p,\dots,i_L)$, with $-y$ at the position $p-1$, and to its exchangeable chamber $\hat\ell$ formed by the slots of $\ell$ which are $\ge p-1$.
  The letter $-y$ changes no column weight; so, by the proof of \cref{lem: sign}, the chambers of $\hat{\boldsymbol{\imath}}$ with a slot $\le p-1$ have $\lambda^\cc_k(\alpha_g^\vee)=0$, except the initial chamber~$\hat k$ at the level $x$, with~$-1$.
  The remaining chambers are those of $\bfi$ whose slots are all $\ge p$, with the same entries in the columns $\hat\ell$ and~$\ell$.
  So $\sum_ke_a(k)b_{k\ell}$ is the entry of $\tB(\hat{\boldsymbol{\imath}})$ at $(\hat k,\hat\ell)$, given by the consecutive elements~$p-1<p$ of $M_{xy}(\hat{\boldsymbol{\imath}})$; it vanishes, as $\varphi_{xy}(p-1)=-1=\varphi_{xy}(p)$.
\end{proof}

For a double reduced word $\bfi$ let $\Sigma^\Phi_{\gamma,\mathrm{fr}}(\bfi)$ be the seed obtained from $\Sigma^\Phi_\gamma(\bfi)$ by freezing the thawed chambers.
This seed and $\Sigma(\bfi)$ have the $\bfi$-exchangeable chambers as exchangeable vertices.
For a seed~$(\bx,\tB)$ and an exchangeable index $\ell$ let $\hat y_\ell:=\prod_kx_k^{b_{k\ell}}$.

\begin{corollary}
  \label{cor: Fraser}
  For every double reduced word $\bfi$, $\Sigma(\bfi)$ and $\Sigma^\Phi_{\gamma,\mathrm{fr}}(\bfi)$ have the same principal part and the same $\hat y_\ell$ for every $\bfi$-exchangeable chamber $\ell$.
  So the identity of $\bbC(\mcG)$ is a quasi-isomorphism, in the sense of \cite[Definition~3.5]{Fra16}, between the seed patterns of $\Sigma(\bfi)$ and of $\Sigma^\Phi_{\gamma,\mathrm{fr}}(\bfi)$.
\end{corollary}
\begin{proof}

  The principal parts agree by \cref{lem: local entries}(b).
  By \cref{lem: local entries}(a), the column $\ell$ of $\tB_\gamma(\bfi)$ is the column $\ell$ of $\tB(\bfi)$ plus the dipole entries $d_{t_a\ell}$.
  As $\Phi_\gamma(k;\bfi)=\Delta(k;\bfi)\prod_a\Phi_\gamma(t_a;\bfi)^{e_a(k)}$ for every $k\in\mcI$, the quotient of the two $\hat y_\ell$ is $\prod_a\Phi_\gamma(t_a;\bfi)^{\sum_ke_a(k)b_{k\ell}+d_{t_a\ell}}$, which is $1$ by \cref{lem: exponents} and \cref{def: BD word seed}(c).
  For~$\gamma=(0,0)$ the two seeds coincide.
  Otherwise $\mcG$ has rank at least $2$, and so do the seeds, as required in \cite[Proposition~2.3]{Fra16}.
  The last claim then follows from \cite[Propositions~2.3 and~3.2]{Fra16}, applied to the identity in both directions: each $\Delta(k;\bfi)$ is $\Phi_\gamma(k;\bfi)$ times a Laurent monomial in the frozen variables of $\Sigma^\Phi_{\gamma,\mathrm{fr}}(\bfi)$, and conversely, by the unitriangularity after \cref{def: BD word seed}.
\end{proof}

\begin{proof}[Proof of \cref{thm: moves}(a)]
  Let the commutation swap the letters $i_j,i_{j+1}$, at the levels $x\neq y$.
  If both are negative, then $a_{xy}=0$, the roots $u_{\le j-1}(\alpha_x)$ and $u_{\le j}(\alpha_y)=u_{\le j-1}(\alpha_y)$ of the two letters are swapped, and the other roots do not change.
  Similarly, if both are positive, the roots $v_{>j}(\alpha_x)=v_{>j+1}(\alpha_x)$ and~$v_{>j+1}(\alpha_y)$ are swapped.
  If the letters have opposite signs, no root changes.
  So the inversion letters of $\bfi$ are those of $\bfi'$, possibly moved from $j$ to $j+1$ or conversely, and the thawed chambers correspond.
  The order of two inversion letters changes only if they are $i_j,i_{j+1}$, negative say.
  Then their roots are orthogonal, so they are $\alpha_{\gamma_{\rr}(a)}$, $\alpha_{\gamma_{\rr}(a')}$ with $a,a'$ not adjacent, $\gamma_{\rr}$ being an isometry, and the virtual letters of $a$ and $a'$ commute.
  So $(\bfi')^\gamma$ is obtained from $\bfi^\gamma$ by commutations, and $\tB((\bfi')^\gamma)_\mcI=\tB(\bfi^\gamma)_\mcI$ by \cref{lem: commutation}.
  The two chambers separated by an inversion letter correspond to the two chambers separated by the same letter in $\bfi'$, so $D_\gamma(\bfi')=D_\gamma(\bfi)$.
  The elements $u_{\le z}$, $v_{>z}$ do not change for slots $z\neq j$, and every chamber contains a slot $z\neq j$.
  So the weights, the minors and the exponents of corresponding chambers coincide.
\end{proof}

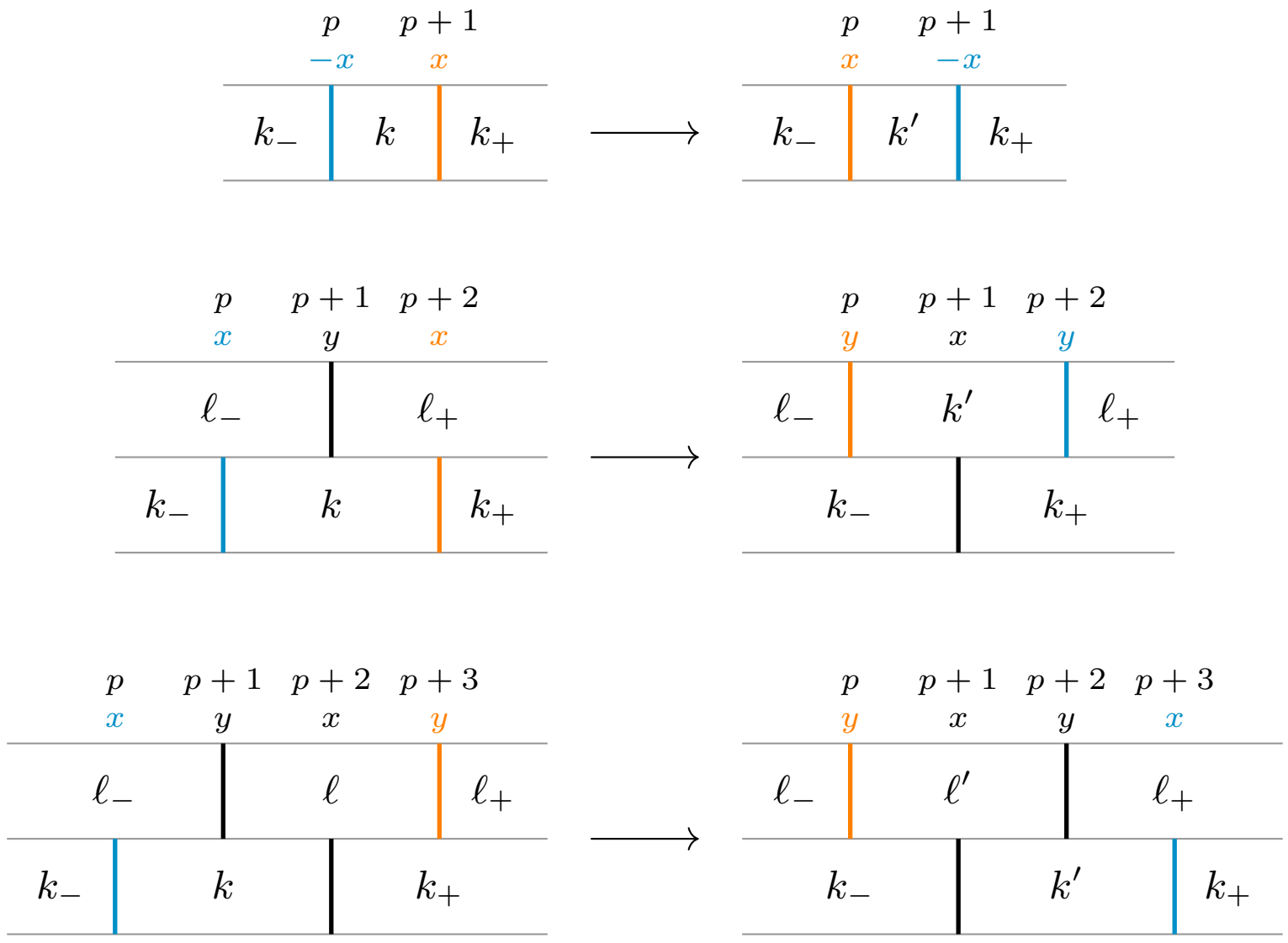
\begin{figure}[tbp]
  \centering
  \begin{tikzpicture}[x=.85cm,y=.75cm,every node/.append style={text height=1.5ex,text depth=.25ex}]
    \draw[black!40] (-3.9,7.9) -- (-0.9,7.9);
    \draw[black!40] (-3.9,8.9) -- (-0.9,8.9);
    \draw[accent,line width=1pt] (-2.9,7.9) -- (-2.9,8.9);
    \node[font=\scriptsize,text=accent] at (-2.9,9.2) {$-x$};
    \draw[orange,line width=1pt] (-1.9,7.9) -- (-1.9,8.9);
    \node[font=\scriptsize,text=orange] at (-1.9,9.2) {$x$};
    \node[font=\scriptsize] at (-2.9,9.6) {$p$};
    \node[font=\scriptsize] at (-1.9,9.6) {$p+1$};
    \node[font=\small] at (-3.4,8.4) {$k_-$};
    \node[font=\small] at (-2.4,8.4) {$k$};
    \node[font=\small] at (-1.4,8.4) {$k_+$};
    \draw[->] (-0.5,8.4) -- (0.5,8.4);
    \draw[black!40] (0.9,7.9) -- (3.9,7.9);
    \draw[black!40] (0.9,8.9) -- (3.9,8.9);
    \draw[orange,line width=1pt] (1.9,7.9) -- (1.9,8.9);
    \node[font=\scriptsize,text=orange] at (1.9,9.2) {$x$};
    \draw[accent,line width=1pt] (2.9,7.9) -- (2.9,8.9);
    \node[font=\scriptsize,text=accent] at (2.9,9.2) {$-x$};
    \node[font=\scriptsize] at (1.9,9.6) {$p$};
    \node[font=\scriptsize] at (2.9,9.6) {$p+1$};
    \node[font=\small] at (1.4,8.4) {$k_-$};
    \node[font=\small] at (2.4,8.4) {$k'$};
    \node[font=\small] at (3.4,8.4) {$k_+$};
    \draw[black!40] (-4.9,4) -- (-0.9,4);
    \draw[black!40] (-4.9,5) -- (-0.9,5);
    \draw[black!40] (-4.9,6) -- (-0.9,6);
    \draw[accent,line width=1pt] (-3.9,4) -- (-3.9,5);
    \node[font=\scriptsize,text=accent] at (-3.9,6.3) {$x$};
    \draw[black,line width=1pt] (-2.9,5) -- (-2.9,6);
    \node[font=\scriptsize,text=black] at (-2.9,6.3) {$y$};
    \draw[orange,line width=1pt] (-1.9,4) -- (-1.9,5);
    \node[font=\scriptsize,text=orange] at (-1.9,6.3) {$x$};
    \node[font=\scriptsize] at (-3.9,6.7) {$p$};
    \node[font=\scriptsize] at (-2.9,6.7) {$p+1$};
    \node[font=\scriptsize] at (-1.9,6.7) {$p+2$};
    \node[font=\small] at (-4.4,4.5) {$k_-$};
    \node[font=\small] at (-2.9,4.5) {$k$};
    \node[font=\small] at (-1.4,4.5) {$k_+$};
    \node[font=\small] at (-3.9,5.5) {$\ell_-$};
    \node[font=\small] at (-1.9,5.5) {$\ell_+$};
    \draw[->] (-0.5,5) -- (0.5,5);
    \draw[black!40] (0.9,4) -- (4.9,4);
    \draw[black!40] (0.9,5) -- (4.9,5);
    \draw[black!40] (0.9,6) -- (4.9,6);
    \draw[orange,line width=1pt] (1.9,5) -- (1.9,6);
    \node[font=\scriptsize,text=orange] at (1.9,6.3) {$y$};
    \draw[black,line width=1pt] (2.9,4) -- (2.9,5);
    \node[font=\scriptsize,text=black] at (2.9,6.3) {$x$};
    \draw[accent,line width=1pt] (3.9,5) -- (3.9,6);
    \node[font=\scriptsize,text=accent] at (3.9,6.3) {$y$};
    \node[font=\scriptsize] at (1.9,6.7) {$p$};
    \node[font=\scriptsize] at (2.9,6.7) {$p+1$};
    \node[font=\scriptsize] at (3.9,6.7) {$p+2$};
    \node[font=\small] at (1.9,4.5) {$k_-$};
    \node[font=\small] at (3.9,4.5) {$k_+$};
    \node[font=\small] at (1.4,5.5) {$\ell_-$};
    \node[font=\small] at (2.9,5.5) {$k'$};
    \node[font=\small] at (4.4,5.5) {$\ell_+$};
    \draw[black!40] (-5.9,0) -- (-0.9,0);
    \draw[black!40] (-5.9,1) -- (-0.9,1);
    \draw[black!40] (-5.9,2) -- (-0.9,2);
    \draw[accent,line width=1pt] (-4.9,0) -- (-4.9,1);
    \node[font=\scriptsize,text=accent] at (-4.9,2.3) {$x$};
    \draw[black,line width=1pt] (-3.9,1) -- (-3.9,2);
    \node[font=\scriptsize,text=black] at (-3.9,2.3) {$y$};
    \draw[black,line width=1pt] (-2.9,0) -- (-2.9,1);
    \node[font=\scriptsize,text=black] at (-2.9,2.3) {$x$};
    \draw[orange,line width=1pt] (-1.9,1) -- (-1.9,2);
    \node[font=\scriptsize,text=orange] at (-1.9,2.3) {$y$};
    \node[font=\scriptsize] at (-4.9,2.7) {$p$};
    \node[font=\scriptsize] at (-3.9,2.7) {$p+1$};
    \node[font=\scriptsize] at (-2.9,2.7) {$p+2$};
    \node[font=\scriptsize] at (-1.9,2.7) {$p+3$};
    \node[font=\small] at (-5.4,0.5) {$k_-$};
    \node[font=\small] at (-3.9,0.5) {$k$};
    \node[font=\small] at (-1.9,0.5) {$k_+$};
    \node[font=\small] at (-4.9,1.5) {$\ell_-$};
    \node[font=\small] at (-2.9,1.5) {$\ell$};
    \node[font=\small] at (-1.4,1.5) {$\ell_+$};
    \draw[->] (-0.5,1) -- (0.5,1);
    \draw[black!40] (0.9,0) -- (5.9,0);
    \draw[black!40] (0.9,1) -- (5.9,1);
    \draw[black!40] (0.9,2) -- (5.9,2);
    \draw[orange,line width=1pt] (1.9,1) -- (1.9,2);
    \node[font=\scriptsize,text=orange] at (1.9,2.3) {$y$};
    \draw[black,line width=1pt] (2.9,0) -- (2.9,1);
    \node[font=\scriptsize,text=black] at (2.9,2.3) {$x$};
    \draw[black,line width=1pt] (3.9,1) -- (3.9,2);
    \node[font=\scriptsize,text=black] at (3.9,2.3) {$y$};
    \draw[accent,line width=1pt] (4.9,0) -- (4.9,1);
    \node[font=\scriptsize,text=accent] at (4.9,2.3) {$x$};
    \node[font=\scriptsize] at (1.9,2.7) {$p$};
    \node[font=\scriptsize] at (2.9,2.7) {$p+1$};
    \node[font=\scriptsize] at (3.9,2.7) {$p+2$};
    \node[font=\scriptsize] at (4.9,2.7) {$p+3$};
    \node[font=\small] at (1.9,0.5) {$k_-$};
    \node[font=\small] at (3.9,0.5) {$k'$};
    \node[font=\small] at (5.4,0.5) {$k_+$};
    \node[font=\small] at (1.4,1.5) {$\ell_-$};
    \node[font=\small] at (2.9,1.5) {$\ell'$};
    \node[font=\small] at (4.9,1.5) {$\ell_+$};
  \end{tikzpicture}
  \caption{%
    The moves in the proof of \cref{thm: moves}, starting at the position $p$, with the chambers named as in the proof: the sign swap $(-x,x)\to(x,-x)$, and the braid moves~$(x,y,x)\to(y,x,y)$ and $(x,y,x,y)\to(y,x,y,x)$.
    Above each letter are its position and its value; the colored letters are those whose root can be simple, and letters with the same root have the same color.}
  \label{fig: moves}
\end{figure}

\begin{proof}[Proof of the statements on the matrices in \cref{thm: moves}(b)]
  We separate the ordinary word matrix from the dipole corrections.
  After checking that they do not cancel at the mutated chamber, we compare the corrections and the order of virtual letters.

  Under the identification of \cref{lem: mirror}, sign swaps become sign swaps, braid moves of negative letters become braid moves of positive letters, and row roots become column roots.
  As mutation is involutive, the reverse of a move satisfies the same identity.
  Thus it suffices to treat the sign swap $(-x,x)\to(x,-x)$ and the positive braid move $(x,y,x)\to(y,x,y)$, $a_{xy}=a_{yx}=-1$, both starting at $p$.
  Let $p'$ be the last position of the move and $v:=v_{>p'}$.

  \emph{Local data.} \Cref{fig: moves} shows the chambers named below.
  For the sign swap, let $k_-=c_x(p-1)$, $k=c_x(p)$ and $k_+=c_x(p+1)$; in $\bfi'$ the replaced chamber is $k'=c_x(p)$.
  By \eqref{eq: B entries}, the non-zero entries of the column~$k$ of $\tB(\bfi)$ are $b_{k_-k}=b_{k_+k}=-1$ and $b_{c_{y'}(p),k}=-a_{y'x}$ for the levels $y'\neq x$.
  For the braid move, let $k_-=c_x(p-1)$, $k=c_x(p)$, $k_+=c_x(p+2)$, $\ell_-=c_y(p)$ and $\ell_+=c_y(p+1)$; in $\bfi'$ the replaced chamber is $k'=c_y(p)$.
  The non-zero entries of the column $k$ of $\tB(\bfi)$ are $b_{k_-k}=b_{\ell_+k}=1$ and $b_{k_+k}=b_{\ell_-k}=-1$, as for a level $y'\notin\{x,y\}$ the letters at $p$ and $p+2$ have the same value of $\varphi_{y'x}$.

  \emph{Dipoles at $k$.} By \cref{lem: local entries}, the column $k$ of $\tB_\gamma(\bfi)$ is the column $k$ of $\tB(\bfi)$ plus the dipole entries $d_{t_ak}$ of the roots $a$ with $p_a\in\{p,p'\}$.
  In the sign swap these are a row root with $p_a=p$ and a column root with $p_b=p+1$.
  In the braid move they are column roots $b$ with $p_b\in\{p,p+2\}$, the roots of the letters at $p+2,p+1,p$ being $v(\alpha_x)$, $v(\alpha_x+\alpha_y)$, $v(\alpha_y)$.
  For a column root $b$, $d_{t_bk}=1$ if $p_b=p'$ and $d_{t_bk}=-1$ if~$p_b=p$, by \cref{def: BD word seed}(c).
  No dipole entry of a column root $b$ cancels an entry of $\tB(\bfi)$.
  Indeed, an entry of the opposite sign at $(t_b,k)$ would force $t_b=k_-$ in a sign swap.
  In a braid move it would force $t_b=\ell_-$ if $p_b=p+2$, resp.\ $t_b=k_-$ if $p_b=p$.
  These are the only possibilities because the other chambers with a non-zero entry in the column $k$ are not initial.
  With $g:=\gamma_{\cc}(b)$, $\lambda^\cc_{t_b}(\alpha_g^\vee)=-\varpi_b(\alpha_g^\vee)=0$ by \cref{lem: tails}, as $\gamma_{\cc}$ is nilpotent.
  On the other hand, $\lambda^\cc_{k_-}(\alpha_g^\vee)=-1$ in the first and third cases by \cref{lem: sign}.
  Indeed, $k_-$ is $c_x(p_b-1)$ in the third case, and has the column weight of $c_x(p_b-1)=k$ in the first, the letter $i_p$ being negative.
  In the second case $\lambda^\cc_{\ell_-}(\alpha_g^\vee)=-1$, as $v_{>p}^{-1}(\alpha_g^\vee)=s_y(-\alpha_x^\vee)=-\alpha_x^\vee-\alpha_y^\vee$.
  So $t_b$ is none of these chambers.
  Row roots follow by \cref{lem: mirror}.

  \emph{The matrix $\tB(\bfi^\gamma)_\mcI$.} The ordinary word mutation gives the extended matrix with the old virtual order.
  It remains to account for the dipoles and any change in that order.
  By \cref{thm: SSVZ}, $\mu_k(\tB(\bfi^\gamma))$ is the matrix of the word obtained from $\bfi^\gamma$ by the move, i.e.\ of $(\bfi')^\gamma$ computed with the virtual letters of $\bfi^\gamma$, in the same order.
  This is $\tB((\bfi')^\gamma)$ unless the move changes the order of the virtual letters.
  As $k\in \mcI_\ex$, restricting to the rows $\mcI$ and the columns $\mcI_\ex$ commutes with $\mu_k$.

  \emph{Dipoles.} If no inversion letter is a letter of the move, no dipole entry involves $k$, the inversion letters, the chambers they separate and the order of the virtual letters do not change; therefore we have $\mu_k(\tB_\gamma(\bfi))=\mu_k(\tB(\bfi^\gamma)_\mcI)+D_\gamma(\bfi)=\tB_\gamma(\bfi')$.
  Otherwise, as no dipole entry cancels an entry of $\tB(\bfi)$ in the column $k$, the corrections in the mutation rule split into those of $\tB(\bfi^\gamma)_\mcI$ and those involving a dipole entry at $k$.
  The latter change the columns $t_a$ of the roots $a$ with $p_a$ in the move as follows, which gives~$D_\gamma(\bfi')$.
  In \cref{fig: moves} each inversion letter of the move has the color of the inversion letter it becomes in $\bfi'$.
  For a column root $b$ in the sign swap, the column $t_b$ of the correction has entries $-1$ at $k$ and $1$ at~$k_+$.
  It becomes the column with $1$ at $k'$ and $-1$ at $k_-$, the entry at $k_+$ canceling.
  Indeed in $\bfi'$, $p_b=p$, $c_x(p-1)=k_-$ and $c_x(p)=k'$.
  In the braid move, if $p_b=p+2$, the column $t_b$, with $-1$ at $k$ and $1$ at~$k_+$, becomes the column with $1$ at $k'$ and $-1$ at $\ell_-$, and in $\bfi'$, $p_b=p$, $c_y(p-1)=\ell_-$ and $c_y(p)=k'$.
  If~$p_b=p$, the column $t_b$, with $-1$ at $k_-$ and $1$ at $k$, becomes the column with $-1$ at $k'$ and $1$ at $\ell_+$, and in $\bfi'$, $p_b=p+2$, $c_y(p+1)=k'$ and $c_y(p+2)=\ell_+$.
  If in a sign swap $D_a$ and $D_b$ have entries at $k$ for a row root $a$ and a column root $b$, then $d_{t_ak}=d_{t_bk}=1$.
  Hence the mutation at $k$ does not change the entry at $(t_a,t_b)$.
  In a braid move two matrices $D_b$, $D_{b'}$ have entries at $k$ only if $p_b=p+2$ and $p_{b'}=p$.
  Then $d_{t_bk}=1$ and $d_{t_{b'}k}=-1$, the mutation at $k$ adds $1$ at $(t_b,t_{b'})$ and $-1$ at $(t_{b'},t_b)$, and this is accounted for by the change of the order of the virtual letters, treated next.

  \emph{Virtual letters.} The remaining correction occurs when two inversion letters change order; reversing their virtual letters gives exactly this correction.

  The order of two inversion letters of the same kind changes only in this last case.
  Then $\alpha_{\gamma_{\cc}(b)}=v(\alpha_x)$ and $\alpha_{\gamma_{\cc}(b')}=v(\alpha_y)$, so $b,b'$ are joined by a single bond, as $x,y$ are, and the virtual letters of $b$ and $b'$ are at the levels $x':=b^*$ and $y':=b'^*$, with $a_{x'y'}=-1$.
  In $\bfi'$ the inversion letters of $b$ and $b'$ are at $p$ and $p+2$, so their order is reversed, and so is the order of the two virtual letters.
  In $M_{x'y'}(\bfi^\gamma)$ the two virtual letters are positive.
  The only entries depending on their order are those at $(t_b,t_{b'})$ and $(t_{b'},t_b)$; by skew-symmetry it suffices to treat the first.
  It is given by the pair formed by the second virtual letter and the next element $j$ of $M_{x'y'}$.
  Indeed, the entries for the slots before it involve chambers of $\bfi^\gamma$ which are not chambers of $\bfi$.
  In $\bfi$, $p_{b'}<p_b$, so the virtual letter of $b$ comes first.
  By \eqref{eq: B entries} the entry is $\frac12(\varphi_{x'y'}(j)-1)$, and in $\bfi'$, where the virtual letter of $b$ comes second, it is $\frac12(\varphi_{x'y'}(j)+1)$, larger by $1$.
\end{proof}

\begin{proof}[Proof of the statements on the matrices in \cref{thm: moves}(c)]
  We again separate the word matrix from the dipoles, now through three mutations.
  At most one dipole is relevant, and we follow it explicitly.

  By \cref{lem: mirror}, and since $\mu_k\mu_\ell\mu_k$ is an involution, so that the reverse of a move satisfies the same identity, it suffices to treat the positive braid move $(x,y,x,y)\to(y,x,y,x)$ starting at $p$, with $\alpha_x$ long, so that $a_{xy}=-1$ and $a_{yx}=-2$.
  We keep the notation of the proof of (b), with $p'=p+3$ and $v:=v_{>p'}$.

  \emph{Local data.} As in \cref{fig: moves}, let $k_-=c_x(p-1)$, $k=c_x(p)$, $k_+=c_x(p+2)$, $\ell_-=c_y(p)$, $\ell=c_y(p+1)$ and~$\ell_+=c_y(p+3)$; in $\bfi'$ the replaced chambers are $k'=c_x(p+1)$ and $\ell'=c_y(p)$, corresponding to $k$ and $\ell$ by level.
  By \eqref{eq: B entries}, the non-zero entries of the columns $k$ and $\ell$ of $\tB(\bfi)$ are $b_{k_-k}=1$, $b_{k_+k}=-1$, $b_{\ell_-k}=-2$, $b_{\ell k}=2$, and $b_{\ell_-\ell}=1$, $b_{\ell_+\ell}=-1$, $b_{k\ell}=-1$, $b_{k_+\ell}=1$.
  So along $\mu_k,\mu_\ell,\mu_k$ those of the columns of the mutated chambers in the matrices obtained from $\tB(\bfi)$ are these, then $b_{k\ell}=1$, $b_{\ell_-\ell}=b_{\ell_+\ell}=-1$, then $b_{k_-k}=-1$, $b_{k_+k}=1$, $b_{\ell_+k}=-2$, $b_{\ell k}=2$.

  \emph{Roots in the move.} We identify the possible inversion letters in the move and exclude cancellation with the ordinary matrix at either replaced chamber.

  The roots of the letters at $p,p+1,p+2,p+3$ are $v(\alpha_x)$, $v(\alpha_x+\alpha_y)$, $v(\alpha_x+2\alpha_y)$ and $v(\alpha_y)$, and only the first and the last can be simple, the others being sums of two of these positive roots.
  So the only roots whose inversion letter is a letter of the move are column roots $b$ with $p_b=p$ or $p_b=p+3$, colored in \cref{fig: moves}.
  There is at most one such root: otherwise $\gamma_{\cc}$ would map two simple roots joined by a double bond to the simple roots $v(\alpha_x)$, $v(\alpha_y)$, also joined by a double bond.
  As $\gamma_{\cc}$ is an isometry and $\mcG$ has only one such pair, $\gamma_{\cc}$ would fix both of its roots, contradicting nilpotency.
  In particular the move does not change the order of the virtual letters, and the dipole entries at $k$ and $\ell$ are those of $D_b$.
  If there is no such root, there are none, and the mutations of $\tB_\gamma(\bfi)$ are those of $\tB(\bfi^\gamma)_\mcI$.
  We may therefore assume that $b$ exists.
  As $t_b$ is initial, its entries in the columns $k$ and $\ell$ of $\tB(\bfi)$ vanish unless $t_b=k_-$ or~$t_b=\ell_-$.
  Put $g:=\gamma_{\cc}(b)$; then $\lambda^\cc_{t_b}(\alpha_g^\vee)=0$ by \cref{lem: tails} and nilpotency, while the level $b^*$ of $t_b$ has the root length of $\alpha_g$, as $\gamma_{\cc}$ is an isometry.
  If $p_b=p$, root lengths exclude $t_b=\ell_-$, and $\lambda^\cc_{k_-}(\alpha_g^\vee)=-1$ by \cref{lem: sign} excludes $t_b=k_-$.
  If $p_b=p+3$, root lengths exclude $t_b=k_-$, and $v_{>p}^{-1}(\alpha_g^\vee)=-2\alpha_x^\vee-\alpha_y^\vee$ gives $\lambda^\cc_{\ell_-}(\alpha_g^\vee)=-1$, excluding $t_b=\ell_-$.
  Thus both entries vanish.

  \emph{Mutation calculation.} As the ordinary row vanishes, the dipole correction can be followed separately.

  By \cref{lem: local entries}, the matrices obtained from $\tB_\gamma(\bfi)$ and from $\tB(\bfi)$, in the columns $k$ and $\ell$, initially differ only in the row $t_b$.
  Write $\delta=(\delta_k,\delta_\ell)$ for this correction, and $E$ for the correction to the matrix obtained from $\tB(\bfi^\gamma)_\mcI$ due to $D_b$; initially $E=D_b$, and $(E_{t_bk},E_{t_b\ell})=\delta$.
  The tables below record these corrections along the three mutations.
  Entries of $E$ not listed are zero, except those in the column $t_b$, which follow from the same mutations.
  As for $D_b$, entries in columns outside $\mcI_\ex$ are omitted.
  The row~$t_b$ of the matrices obtained from $\tB(\bfi)$ remains zero in the columns $k,\ell$ throughout the three mutations.
  Indeed, each mutation is at a column where that row is zero, so the matrix mutation formula leaves the row unchanged.
  Hence no cancellation occurs at the mutated chamber.
  Thus the mutations of $\tB_\gamma(\bfi)$ split into those of $\tB(\bfi^\gamma)_\mcI$ and the corrections $E$.
  The corrections in the two cases are:
  \[
    \begin{array}{c|c|c}
      \multicolumn{3}{c}{p_b=p}\\[2pt]
      \text{stage}&\delta&\text{other entries of }E\\ \hline
      \text{initial}&(-1,0)&E_{t_bk_-}=1\\
      \mu_k&(1,-1)&\text{none}\\
      \mu_\ell&(-1,1)&\text{none}\\
      \mu_k&(1,0)&E_{t_bk_+}=-1
    \end{array}
    \qquad
    \begin{array}{c|c|c}
      \multicolumn{3}{c}{p_b=p+3}\\[2pt]
      \text{stage}&\delta&\text{other entries of }E\\ \hline
      \text{initial}&(0,1)&E_{t_b\ell_+}=-1\\
      \mu_k&(0,1)&E_{t_b\ell_+}=-1\\
      \mu_\ell&(0,-1)&E_{t_b\ell_-}=1\\
      \mu_k&(0,-1)&E_{t_b\ell_-}=1
    \end{array}
  \]

  \emph{Matrices.} By \cref{thm: SSVZ} for $\bfi^\gamma$, whose virtual letters keep their order, $\mu_k\mu_\ell\mu_k(\tB(\bfi^\gamma)_\mcI)=\tB((\bfi')^\gamma)_\mcI$, and the other dipoles do not meet $k$ or $\ell$ and remain unchanged.
  If $p_b=p$, the final correction in the column $t_b$ is $-1$ at $k'$ and $1$ at $k_+$, and in $\bfi'$ the inversion letter is at $p+3$, separating $k'=c_x(p+2)$ from $k_+=c_x(p+3)$.
  If $p_b=p+3$, it is $1$ at $\ell'$ and $-1$ at $\ell_-$, and the inversion letter is at $p$, separating~$\ell_-=c_y(p-1)$ from $\ell'=c_y(p)$.
  Thus in both cases the correction is $D_b(\bfi')$, as required.
\end{proof}

To recover the functions from \cref{cor: Fraser}, we use the following lemma.

\begin{lemma}
  \label{lem: BFZ indecomposable}
  For every double reduced word $\bfi$ the principal part of $\tB(\bfi)$ is block-indecomposable.
  In particular, if $\mcG$ has rank at least $2$, every $\bfi$-exchangeable chamber $k$ has an $\bfi$-exchangeable chamber $\ell$ with~$b_{\ell k}\neq0$.
\end{lemma}
\begin{proof}
  We first rule out a separation of two levels in a reduced word.
  We then use this observation to find an exchangeable connection between every pair of adjacent levels.

  A reduced word for $w_0$ cannot have all its letters at a level $x$ before all its letters at another level $y$.
  Otherwise $w_0=uv$ with $u$ not involving $s_y$ and $v$ not involving $s_x$.
  Then $v$ fixes $\varpi_x$, and so~$\varpi_x+\varpi_{x^*}=\varpi_x-w_0(\varpi_x)=\varpi_x-u(\varpi_x)$.
  As $s_{y'}(\lambda)=\lambda-\lambda(\alpha_{y'}^\vee)\alpha_{y'}$, this is a combination of the $\alpha_{y'}$ with $s_{y'}$ in $u$, so its coefficient at $\alpha_y$ is $0$.
  But the coefficients of the fundamental weights in the simple roots are strictly positive, the root system being irreducible \cite[Theorem~4.3 (Fin)]{Kac90}.

  At every level there are letters of both signs, so at least one exchangeable chamber, and consecutive exchangeable chambers $k,k'$ at the same level have $b_{kk'}=\pm1$.
  Let $x,y$ be levels with $a_{xy}\neq0$.
  We will find exchangeable chambers $k,\ell$ at the levels $x$ and $y$ with $b_{k\ell}\neq0$.
  Let $F_x,L_x$ and $F_y,L_y$ be the first and last positions of letters at these levels, and assume $F_x<F_y$.
  Then $F_y<L_x$, as otherwise every letter at the level $x$ comes before every letter at the level $y$.
  For consecutive elements $p<p'$ of $M_{xy}$ with~$p\in[F_y,\min(L_x,L_y)-1]$, the chambers $c_x(p)$ and $c_y(p)$ are exchangeable, and $b_{c_x(p)c_y(p)}$ vanishes only if $\varphi_{xy}(p)=\varphi_{xy}(p')$.
  If all these entries vanish, $\varphi_{xy}$ is constant, equal to $c$ say, on $M_{xy}\cap[F_y,\min(L_x,L_y)]$.
  If $L_y<L_x$, then all the letters at the level $y$ have the sign $c$, which is impossible.
  If $L_x<L_y$, then the letters at the level $x$ after $F_y$ have the sign $-c$ and those at the level $y$ before $L_x$ have the sign $c$.
  So in the reduced word for $w_0$ formed by the letters of sign $c$ all the letters at the level $x$ come before those at the level $y$, which is impossible.
\end{proof}

\begin{proof}[Proof of the statements on the functions in \cref{thm: moves}(b) and (c)]
  The matrix calculation and quasi-isomorphism leave only possible frozen monomial factors at the replaced chambers.
  We eliminate these factors by comparing exchange monomials and using block-indecomposability.

  If $\gamma=(0,0)$, these statements are \cref{thm: SSVZ}; so we may assume that $\mcG$ has rank at least $2$.
  Let $\mu$ be the mutation, or the sequence of mutations, of the statement, and $S:=\mu(\Sigma^\Phi_\gamma(\bfi))$.
  By the statements on the matrices, $S$ has the matrix $\tB_\gamma(\bfi')$.
  Its functions at the chambers which are not replaced are those of~$\Sigma^\Phi_\gamma(\bfi')$.
  Indeed, the weights of such a chamber can be computed at a slot outside the move, where $u_{\le z}$ and $v_{>z}$ do not change.
  Also, the $\Phi_\gamma(t_a;\bfi)$ do not depend on the word by \cref{lem: tails}.
  Freezing the thawed chambers commutes with $\mu$; let $S_\mathrm{fr}$ be $S$ with them frozen.
  The quasi-isomorphism of \cref{cor: Fraser} compares $S_\mathrm{fr}$ with $\mu(\Sigma(\bfi))=\Sigma(\bfi')$ (\cref{thm: SSVZ}): their $\hat y_\ell$ agree and corresponding functions differ by Laurent monomials in the frozen variables.
  Applying \cref{cor: Fraser} at $\bfi'$ gives the same comparison with $\Sigma^\Phi_{\gamma,\mathrm{fr}}(\bfi')$.
  The frozen variables of all three seeds generate the same group of Laurent monomials by unitriangularity.
  So $S_\mathrm{fr}$ and $\Sigma^\Phi_{\gamma,\mathrm{fr}}(\bfi')$ have the same matrix and the same $\hat y_\ell$.
  Their functions differ at most at the replaced chambers $k$ of $\bfi'$, by Laurent monomials $\prod_tx_t^{n_t(k)}$ in the frozen variables $x_t$.
  Comparing the $\hat y_\ell$, and as the variables of a seed are algebraically independent, $\sum_k n_t(k)b_{k\ell}=0$ for every frozen $t$ and every $\bfi'$-exchangeable $\ell$, the $b_{k\ell}$ being the entries of $\tB(\bfi')$.
  If $k'$ is the only replaced chamber, this gives $n_t(k')=0$ for some $\ell$ with $b_{k'\ell}\neq0$, which exists by \cref{lem: BFZ indecomposable} and skew-symmetrizability.
  For a braid move of length $4$, the choices $\ell=k'$ and $\ell=\ell'$ give $n_t(\ell')b_{\ell'k'}=0$ and $n_t(k')b_{k'\ell'}=0$, with~$b_{\ell'k'}\neq0\neq b_{k'\ell'}$.
  So $S=\Sigma^\Phi_\gamma(\bfi')$.
\end{proof}

\begin{remark}
  The statements on the functions in \cref{thm: moves} can also be proved without quasi-isomorphisms.
  Substituting $\Phi_\gamma(k;\bfi)=\Delta(k;\bfi)\prod_a\Phi_\gamma(t_a;\bfi)^{e_a(k)}$ in the exchange relations of $\Sigma^\Phi_\gamma(\bfi)$ at the replaced chambers reduces them to those of $\Sigma(\bfi)$ (\cref{thm: SSVZ}) and to identities among the exponents $e_a$.
  These identities can be checked by a direct, but more tedious, inspection of the position of the inversion letters relative to the move.
\end{remark}

\begin{corollary}
  \label{cor: mutation equivalent}
  For every $\mcG$ and every aperiodic BD pair $\gamma$, the seeds $\Sigma^\Phi_\gamma(\bfi)$ of all double reduced words~$\bfi$ are related by sequences of mutations, up to relabeling; so are the seeds $\Sigma_\gamma(\bfi)$ if $h_\gamma$ is birational.
\end{corollary}
\begin{proof}

  Any two double reduced words are related by moves (\cref{sec: double words}).
  Moves of length at most $4$ are treated in \cref{thm: moves}.
  Braid moves of length $6$ only occur in type $G_2$, where the only BD pair is $(0,0)$, as BD maps are isometries and the two simple roots have different lengths.
  Then $\Sigma^\Phi_\gamma(\bfi)$ is the seed of \cite[Theorem~2.10]{BFZ05}, and \cref{thm: SSVZ} applies.
\end{proof}

Suppose that $h_\gamma$ is birational.
Then $f\mapsto f\circ h_\gamma^{-1}$ is a quasi-isomorphism from the seed pattern of~$\Sigma(\bfi)$ to that of $\Sigma_\gamma(\bfi)$ with the thawed chambers frozen.
Consider the cluster variables obtained without mutating at the thawed chambers.
Up to Laurent monomials in the frozen variables and in the functions of the thawed chambers, they are those of the cluster structure of \cite{BFZ05}, composed with $h_\gamma^{-1}$.

\begin{lemma}
  \label{lem: indecomposable}
  For every double reduced word $\bfi$, the principal part of $\tB_\gamma(\bfi)$ is block-indecomposable.
\end{lemma}
\begin{proof}
  By \cref{lem: local entries}(b), the submatrix of $\tB_\gamma(\bfi)$ on the $\bfi$-exchangeable chambers is the principal part of $\tB(\bfi)$, which is block-indecomposable by \cref{lem: BFZ indecomposable}.
  It remains to join each thawed chamber to an~$\bfi$-exchangeable chamber.

  Let $a$ be a row root, $t=t^\rr_a$, $i_j$ the last letter at the level $a^*$, and $k=c_{a^*}(j-1)$, which is $\bfi$-exchangeable as every level carries letters of both signs.
  By \eqref{eq: B entries}, $\tB(\bfi^\gamma)_\mcI$ has the entry $\varepsilon(i_j)$ at $(k,t)$.
  As $k$ is not thawed, a dipole entry at $(k,t)$ can only come from $D_a$, and it cancels $\varepsilon(i_j)$ only if the inversion letter $i_{p_a}$ of $a$ is the last negative letter at the level $a^*$.
  Then $w_0=u\,s_{a^*}v$ with $u=u_{\le p_a-1}$ and $v$ not involving~$s_{a^*}$, so $\alpha_{\gamma_{\rr}(a)}=u(\alpha_{a^*})=-w_0v^{-1}(\alpha_{a^*})$.
  Here $v^{-1}(\alpha_{a^*})$ is a positive root with coefficient $1$ at $\alpha_{a^*}$, and $-w_0$ permutes the simple roots, hence preserves heights; so $v^{-1}(\alpha_{a^*})=\alpha_{a^*}$ and $\gamma_{\rr}(a)=a$, contradicting nilpotency.
  Column roots $b$ are treated symmetrically, with $t^\cc_b=c_{b^*}(0)$, the first letter at the level $b^*$ and the first positive letter there.
\end{proof}

\section{The BD word seeds of \texorpdfstring{$SL_N$}{SLN}}
\label{sec: type A}
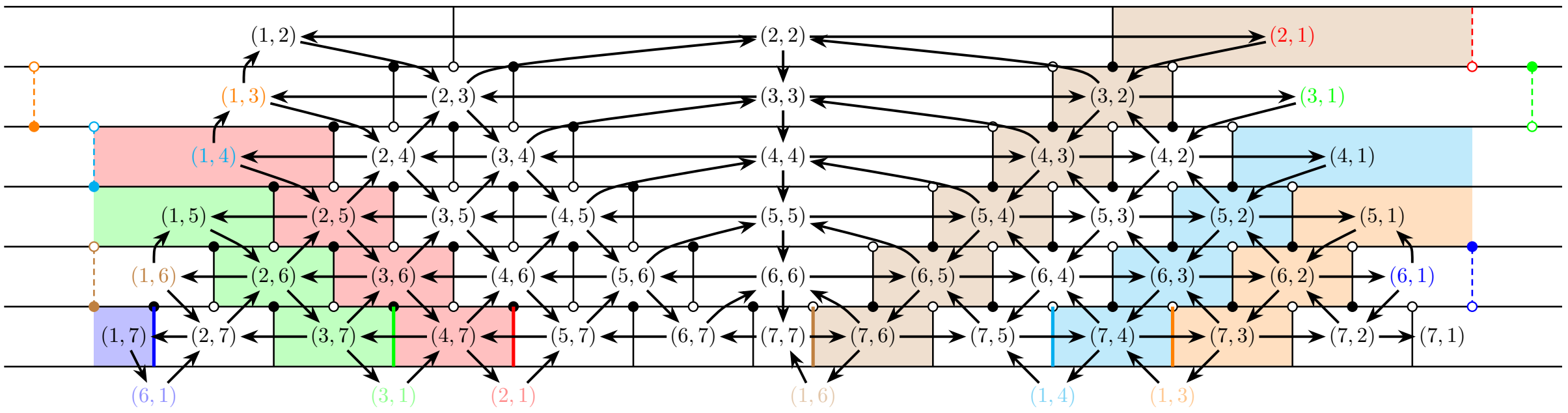
\begin{figure}[!h]
  \centering
  \begingroup
    \resizebox{\linewidth}{!}{%
      \begin{tikzpicture}
        \fill[blue, opacity=0.25] (0,1) rectangle (1,2);
        \fill[green, opacity=0.25] (3,1) rectangle (5,2);
        \fill[red, opacity=0.25] (5,1) rectangle (7,2);
        \fill[brown, opacity=0.25] (12,1) rectangle (14,2);
        \fill[cyan, opacity=0.25] (16,1) rectangle (18,2);
        \fill[orange, opacity=0.25] (18,1) rectangle (20,2);
        \fill[green, opacity=0.25] (2,2) rectangle (4,3);
        \fill[red, opacity=0.25] (4,2) rectangle (6,3);
        \fill[brown, opacity=0.25] (13,2) rectangle (15,3);
        \fill[cyan, opacity=0.25] (17,2) rectangle (19,3);
        \fill[orange, opacity=0.25] (19,2) rectangle (21,3);
        \fill[green, opacity=0.25] (0,3) rectangle (3,4);
        \fill[red, opacity=0.25] (3,3) rectangle (5,4);
        \fill[brown, opacity=0.25] (14,3) rectangle (16,4);
        \fill[cyan, opacity=0.25] (18,3) rectangle (20,4);
        \fill[orange, opacity=0.25] (20,3) rectangle (23,4);
        \fill[red, opacity=0.25] (0,4) rectangle (4,5);
        \fill[brown, opacity=0.25] (15,4) rectangle (17,5);
        \fill[cyan, opacity=0.25] (19,4) rectangle (23,5);
        \fill[brown, opacity=0.25] (16,5) rectangle (18,6);
        \fill[brown, opacity=0.25] (17,6) rectangle (23,7);
        \draw[thick] (24.5,1) -- (-1.5,1);
        \draw[thick] (24.5,2) -- (-1.5,2);
        \draw[thick] (24.5,3) -- (-1.5,3);
        \draw[thick] (24.5,4) -- (-1.5,4);
        \draw[thick] (24.5,5) -- (-1.5,5);
        \draw[thick] (24.5,6) -- (-1.5,6);
        \draw[thick] (24.5,7) -- (-1.5,7);
        \draw[thick] (22,1) -- (22,2);
        \draw[draw=black,fill=white,line width=0.70pt] (22,2) circle[radius=2.20pt];
        \draw[thick] (21,2) -- (21,3);
        \draw[draw=black,fill=black,line width=0.70pt] (21,2) circle[radius=2.20pt];
        \draw[draw=black,fill=white,line width=0.70pt] (21,3) circle[radius=2.20pt];
        \draw[thick] (20,1) -- (20,2);
        \draw[draw=black,fill=white,line width=0.70pt] (20,2) circle[radius=2.20pt];
        \draw[thick] (20,3) -- (20,4);
        \draw[draw=black,fill=black,line width=0.70pt] (20,3) circle[radius=2.20pt];
        \draw[draw=black,fill=white,line width=0.70pt] (20,4) circle[radius=2.20pt];
        \draw[thick] (19,2) -- (19,3);
        \draw[draw=black,fill=black,line width=0.70pt] (19,2) circle[radius=2.20pt];
        \draw[draw=black,fill=white,line width=0.70pt] (19,3) circle[radius=2.20pt];
        \draw[thick] (18,1) -- (18,2);
        \draw[draw=black,fill=white,line width=0.70pt] (18,2) circle[radius=2.20pt];
        \draw[thick] (19,4) -- (19,5);
        \draw[draw=black,fill=black,line width=0.70pt] (19,4) circle[radius=2.20pt];
        \draw[draw=black,fill=white,line width=0.70pt] (19,5) circle[radius=2.20pt];
        \draw[thick] (18,3) -- (18,4);
        \draw[draw=black,fill=black,line width=0.70pt] (18,3) circle[radius=2.20pt];
        \draw[draw=black,fill=white,line width=0.70pt] (18,4) circle[radius=2.20pt];
        \draw[thick] (17,2) -- (17,3);
        \draw[draw=black,fill=black,line width=0.70pt] (17,2) circle[radius=2.20pt];
        \draw[draw=black,fill=white,line width=0.70pt] (17,3) circle[radius=2.20pt];
        \draw[thick] (16,1) -- (16,2);
        \draw[draw=black,fill=white,line width=0.70pt] (16,2) circle[radius=2.20pt];
        \draw[thick] (18,5) -- (18,6);
        \draw[draw=black,fill=black,line width=0.70pt] (18,5) circle[radius=2.20pt];
        \draw[draw=black,fill=white,line width=0.70pt] (18,6) circle[radius=2.20pt];
        \draw[thick] (17,4) -- (17,5);
        \draw[draw=black,fill=black,line width=0.70pt] (17,4) circle[radius=2.20pt];
        \draw[draw=black,fill=white,line width=0.70pt] (17,5) circle[radius=2.20pt];
        \draw[thick] (16,3) -- (16,4);
        \draw[draw=black,fill=black,line width=0.70pt] (16,3) circle[radius=2.20pt];
        \draw[draw=black,fill=white,line width=0.70pt] (16,4) circle[radius=2.20pt];
        \draw[thick] (15,2) -- (15,3);
        \draw[draw=black,fill=black,line width=0.70pt] (15,2) circle[radius=2.20pt];
        \draw[draw=black,fill=white,line width=0.70pt] (15,3) circle[radius=2.20pt];
        \draw[thick] (14,1) -- (14,2);
        \draw[draw=black,fill=white,line width=0.70pt] (14,2) circle[radius=2.20pt];
        \draw[thick] (17,6) -- (17,7);
        \draw[draw=black,fill=black,line width=0.70pt] (17,6) circle[radius=2.20pt];
        \draw[thick] (16,5) -- (16,6);
        \draw[draw=black,fill=black,line width=0.70pt] (16,5) circle[radius=2.20pt];
        \draw[draw=black,fill=white,line width=0.70pt] (16,6) circle[radius=2.20pt];
        \draw[thick] (15,4) -- (15,5);
        \draw[draw=black,fill=black,line width=0.70pt] (15,4) circle[radius=2.20pt];
        \draw[draw=black,fill=white,line width=0.70pt] (15,5) circle[radius=2.20pt];
        \draw[thick] (14,3) -- (14,4);
        \draw[draw=black,fill=black,line width=0.70pt] (14,3) circle[radius=2.20pt];
        \draw[draw=black,fill=white,line width=0.70pt] (14,4) circle[radius=2.20pt];
        \draw[thick] (13,2) -- (13,3);
        \draw[draw=black,fill=black,line width=0.70pt] (13,2) circle[radius=2.20pt];
        \draw[draw=black,fill=white,line width=0.70pt] (13,3) circle[radius=2.20pt];
        \draw[thick] (12,1) -- (12,2);
        \draw[draw=black,fill=white,line width=0.70pt] (12,2) circle[radius=2.20pt];
        \draw[thick] (11,1) -- (11,2);
        \draw[draw=black,fill=black,line width=0.70pt] (11,2) circle[radius=2.20pt];
        \draw[thick] (10,2) -- (10,3);
        \draw[draw=black,fill=white,line width=0.70pt] (10,2) circle[radius=2.20pt];
        \draw[draw=black,fill=black,line width=0.70pt] (10,3) circle[radius=2.20pt];
        \draw[thick] (9,3) -- (9,4);
        \draw[draw=black,fill=white,line width=0.70pt] (9,3) circle[radius=2.20pt];
        \draw[draw=black,fill=black,line width=0.70pt] (9,4) circle[radius=2.20pt];
        \draw[thick] (8,4) -- (8,5);
        \draw[draw=black,fill=white,line width=0.70pt] (8,4) circle[radius=2.20pt];
        \draw[draw=black,fill=black,line width=0.70pt] (8,5) circle[radius=2.20pt];
        \draw[thick] (7,5) -- (7,6);
        \draw[draw=black,fill=white,line width=0.70pt] (7,5) circle[radius=2.20pt];
        \draw[draw=black,fill=black,line width=0.70pt] (7,6) circle[radius=2.20pt];
        \draw[thick] (6,6) -- (6,7);
        \draw[draw=black,fill=white,line width=0.70pt] (6,6) circle[radius=2.20pt];
        \draw[thick] (9,1) -- (9,2);
        \draw[draw=black,fill=black,line width=0.70pt] (9,2) circle[radius=2.20pt];
        \draw[thick] (8,2) -- (8,3);
        \draw[draw=black,fill=white,line width=0.70pt] (8,2) circle[radius=2.20pt];
        \draw[draw=black,fill=black,line width=0.70pt] (8,3) circle[radius=2.20pt];
        \draw[thick] (7,3) -- (7,4);
        \draw[draw=black,fill=white,line width=0.70pt] (7,3) circle[radius=2.20pt];
        \draw[draw=black,fill=black,line width=0.70pt] (7,4) circle[radius=2.20pt];
        \draw[thick] (6,4) -- (6,5);
        \draw[draw=black,fill=white,line width=0.70pt] (6,4) circle[radius=2.20pt];
        \draw[draw=black,fill=black,line width=0.70pt] (6,5) circle[radius=2.20pt];
        \draw[thick] (5,5) -- (5,6);
        \draw[draw=black,fill=white,line width=0.70pt] (5,5) circle[radius=2.20pt];
        \draw[draw=black,fill=black,line width=0.70pt] (5,6) circle[radius=2.20pt];
        \draw[thick] (7,1) -- (7,2);
        \draw[draw=black,fill=black,line width=0.70pt] (7,2) circle[radius=2.20pt];
        \draw[thick] (6,2) -- (6,3);
        \draw[draw=black,fill=white,line width=0.70pt] (6,2) circle[radius=2.20pt];
        \draw[draw=black,fill=black,line width=0.70pt] (6,3) circle[radius=2.20pt];
        \draw[thick] (5,3) -- (5,4);
        \draw[draw=black,fill=white,line width=0.70pt] (5,3) circle[radius=2.20pt];
        \draw[draw=black,fill=black,line width=0.70pt] (5,4) circle[radius=2.20pt];
        \draw[thick] (4,4) -- (4,5);
        \draw[draw=black,fill=white,line width=0.70pt] (4,4) circle[radius=2.20pt];
        \draw[draw=black,fill=black,line width=0.70pt] (4,5) circle[radius=2.20pt];
        \draw[thick] (5,1) -- (5,2);
        \draw[draw=black,fill=black,line width=0.70pt] (5,2) circle[radius=2.20pt];
        \draw[thick] (4,2) -- (4,3);
        \draw[draw=black,fill=white,line width=0.70pt] (4,2) circle[radius=2.20pt];
        \draw[draw=black,fill=black,line width=0.70pt] (4,3) circle[radius=2.20pt];
        \draw[thick] (3,3) -- (3,4);
        \draw[draw=black,fill=white,line width=0.70pt] (3,3) circle[radius=2.20pt];
        \draw[draw=black,fill=black,line width=0.70pt] (3,4) circle[radius=2.20pt];
        \draw[thick] (3,1) -- (3,2);
        \draw[draw=black,fill=black,line width=0.70pt] (3,2) circle[radius=2.20pt];
        \draw[thick] (2,2) -- (2,3);
        \draw[draw=black,fill=white,line width=0.70pt] (2,2) circle[radius=2.20pt];
        \draw[draw=black,fill=black,line width=0.70pt] (2,3) circle[radius=2.20pt];
        \draw[thick] (1,1) -- (1,2);
        \draw[draw=black,fill=black,line width=0.70pt] (1,2) circle[radius=2.20pt];
        \draw[red, thick, densely dashed] (23,6) -- (23,7);
        \draw[draw=red,fill=white,line width=0.70pt] (23,6) circle[radius=2.20pt];
        \draw[green, thick, densely dashed] (24,5) -- (24,6);
        \draw[draw=green,fill=white,line width=0.70pt] (24,5) circle[radius=2.20pt];
        \draw[draw=green,fill=green,line width=0.70pt] (24,6) circle[radius=2.20pt];
        \draw[blue, thick, densely dashed] (23,2) -- (23,3);
        \draw[draw=blue,fill=white,line width=0.70pt] (23,2) circle[radius=2.20pt];
        \draw[draw=blue,fill=blue,line width=0.70pt] (23,3) circle[radius=2.20pt];
        \draw[brown, thick, densely dashed] (0,2) -- (0,3);
        \draw[draw=brown,fill=brown,line width=0.70pt] (0,2) circle[radius=2.20pt];
        \draw[draw=brown,fill=white,line width=0.70pt] (0,3) circle[radius=2.20pt];
        \draw[cyan, thick, densely dashed] (0,4) -- (0,5);
        \draw[draw=cyan,fill=cyan,line width=0.70pt] (0,4) circle[radius=2.20pt];
        \draw[draw=cyan,fill=white,line width=0.70pt] (0,5) circle[radius=2.20pt];
        \draw[orange, thick, densely dashed] (-1,5) -- (-1,6);
        \draw[draw=orange,fill=orange,line width=0.70pt] (-1,5) circle[radius=2.20pt];
        \draw[draw=orange,fill=white,line width=0.70pt] (-1,6) circle[radius=2.20pt];
        \draw[red, ultra thick] (7,1.0) -- (7,2.0);
        \draw[green, ultra thick] (5,1.0) -- (5,2.0);
        \draw[blue, ultra thick] (1,1.0) -- (1,2.0);
        \draw[orange, ultra thick] (18,1.0) -- (18,2.0);
        \draw[cyan, ultra thick] (16,1.0) -- (16,2.0);
        \draw[brown, ultra thick] (12,1.0) -- (12,2.0);
        \tikzset{
          vlabel/.style={inner sep=1pt, outer sep=.5pt, font=\normalsize},
        }
        \node[vlabel, text=black] (dwd0) at (22.5,1.5) {$(7,1)$};
        \node[vlabel, text=black] (dwd1) at (21,1.5) {$(7,2)$};
        \node[vlabel, text=black] (dwd2) at (19,1.5) {$(7,3)$};
        \node[vlabel, text=black] (dwd3) at (17,1.5) {$(7,4)$};
        \node[vlabel, text=black] (dwd4) at (15,1.5) {$(7,5)$};
        \node[vlabel, text=black] (dwd5) at (13,1.5) {$(7,6)$};
        \node[vlabel, text=black] (dwd6) at (11.5,1.5) {$(7,7)$};
        \node[vlabel, text=black] (dwd7) at (10,1.5) {$(6,7)$};
        \node[vlabel, text=black] (dwd8) at (8,1.5) {$(5,7)$};
        \node[vlabel, text=black] (dwd9) at (6,1.5) {$(4,7)$};
        \node[vlabel, text=black] (dwd10) at (4,1.5) {$(3,7)$};
        \node[vlabel, text=black] (dwd11) at (2,1.5) {$(2,7)$};
        \node[vlabel, text=black] (dwd12) at (0.5,1.5) {$(1,7)$};
        \node[vlabel, text=blue] (dwd13) at (22,2.5) {$(6,1)$};
        \node[vlabel, text=black] (dwd14) at (20,2.5) {$(6,2)$};
        \node[vlabel, text=black] (dwd15) at (18,2.5) {$(6,3)$};
        \node[vlabel, text=black] (dwd16) at (16,2.5) {$(6,4)$};
        \node[vlabel, text=black] (dwd17) at (14,2.5) {$(6,5)$};
        \node[vlabel, text=black] (dwd18) at (11.5,2.5) {$(6,6)$};
        \node[vlabel, text=black] (dwd19) at (9,2.5) {$(5,6)$};
        \node[vlabel, text=black] (dwd20) at (7,2.5) {$(4,6)$};
        \node[vlabel, text=black] (dwd21) at (5,2.5) {$(3,6)$};
        \node[vlabel, text=black] (dwd22) at (3,2.5) {$(2,6)$};
        \node[vlabel, text=brown] (dwd23) at (1,2.5) {$(1,6)$};
        \node[vlabel, text=black] (dwd24) at (21.5,3.5) {$(5,1)$};
        \node[vlabel, text=black] (dwd25) at (19,3.5) {$(5,2)$};
        \node[vlabel, text=black] (dwd26) at (17,3.5) {$(5,3)$};
        \node[vlabel, text=black] (dwd27) at (15,3.5) {$(5,4)$};
        \node[vlabel, text=black] (dwd28) at (11.5,3.5) {$(5,5)$};
        \node[vlabel, text=black] (dwd29) at (8,3.5) {$(4,5)$};
        \node[vlabel, text=black] (dwd30) at (6,3.5) {$(3,5)$};
        \node[vlabel, text=black] (dwd31) at (4,3.5) {$(2,5)$};
        \node[vlabel, text=black] (dwd32) at (1.5,3.5) {$(1,5)$};
        \node[vlabel, text=black] (dwd33) at (21,4.5) {$(4,1)$};
        \node[vlabel, text=black] (dwd34) at (18,4.5) {$(4,2)$};
        \node[vlabel, text=black] (dwd35) at (16,4.5) {$(4,3)$};
        \node[vlabel, text=black] (dwd36) at (11.5,4.5) {$(4,4)$};
        \node[vlabel, text=black] (dwd37) at (7,4.5) {$(3,4)$};
        \node[vlabel, text=black] (dwd38) at (5,4.5) {$(2,4)$};
        \node[vlabel, text=cyan] (dwd39) at (2,4.5) {$(1,4)$};
        \node[vlabel, text=green] (dwd40) at (20.5,5.5) {$(3,1)$};
        \node[vlabel, text=black] (dwd41) at (17,5.5) {$(3,2)$};
        \node[vlabel, text=black] (dwd42) at (11.5,5.5) {$(3,3)$};
        \node[vlabel, text=black] (dwd43) at (6,5.5) {$(2,3)$};
        \node[vlabel, text=orange] (dwd44) at (2.5,5.5) {$(1,3)$};
        \node[vlabel, text=red] (dwd45) at (20,6.5) {$(2,1)$};
        \node[vlabel, text=black] (dwd46) at (11.5,6.5) {$(2,2)$};
        \node[vlabel, text=black] (dwd47) at (3,6.5) {$(1,2)$};
        \draw[-Stealth, very thick] (dwd1) -- (dwd0);
        \draw[-Stealth, very thick] (dwd14) -- (dwd13);
        \draw[-Stealth, very thick] (dwd2) -- (dwd1);
        \draw[-Stealth, very thick] (dwd25) -- (dwd24);
        \draw[-Stealth, very thick] (dwd15) -- (dwd14);
        \draw[-Stealth, very thick] (dwd3) -- (dwd2);
        \draw[-Stealth, very thick] (dwd34) -- (dwd33);
        \draw[-Stealth, very thick] (dwd26) -- (dwd25);
        \draw[-Stealth, very thick] (dwd16) -- (dwd15);
        \draw[-Stealth, very thick] (dwd4) -- (dwd3);
        \draw[-Stealth, very thick] (dwd41) -- (dwd40);
        \draw[-Stealth, very thick] (dwd35) -- (dwd34);
        \draw[-Stealth, very thick] (dwd27) -- (dwd26);
        \draw[-Stealth, very thick] (dwd17) -- (dwd16);
        \draw[-Stealth, very thick] (dwd5) -- (dwd4);
        \draw[-Stealth, very thick] (dwd46) -- (dwd45);
        \draw[-Stealth, very thick] (dwd42) -- (dwd41);
        \draw[-Stealth, very thick] (dwd36) -- (dwd35);
        \draw[-Stealth, very thick] (dwd28) -- (dwd27);
        \draw[-Stealth, very thick] (dwd18) -- (dwd17);
        \draw[-Stealth, very thick] (dwd6) -- (dwd5);
        \draw[-Stealth, very thick] (dwd6) -- (dwd7);
        \draw[-Stealth, very thick] (dwd18) -- (dwd19);
        \draw[-Stealth, very thick] (dwd28) -- (dwd29);
        \draw[-Stealth, very thick] (dwd36) -- (dwd37);
        \draw[-Stealth, very thick] (dwd42) -- (dwd43);
        \draw[-Stealth, very thick] (dwd46) -- (dwd47);
        \draw[-Stealth, very thick] (dwd7) -- (dwd8);
        \draw[-Stealth, very thick] (dwd19) -- (dwd20);
        \draw[-Stealth, very thick] (dwd29) -- (dwd30);
        \draw[-Stealth, very thick] (dwd37) -- (dwd38);
        \draw[-Stealth, very thick] (dwd43) -- (dwd44);
        \draw[-Stealth, very thick] (dwd8) -- (dwd9);
        \draw[-Stealth, very thick] (dwd20) -- (dwd21);
        \draw[-Stealth, very thick] (dwd30) -- (dwd31);
        \draw[-Stealth, very thick] (dwd38) -- (dwd39);
        \draw[-Stealth, very thick] (dwd9) -- (dwd10);
        \draw[-Stealth, very thick] (dwd21) -- (dwd22);
        \draw[-Stealth, very thick] (dwd31) -- (dwd32);
        \draw[-Stealth, very thick] (dwd10) -- (dwd11);
        \draw[-Stealth, very thick] (dwd22) -- (dwd23);
        \draw[-Stealth, very thick] (dwd11) -- (dwd12);
        \draw[-Stealth, very thick] (dwd13) .. controls (21.5,2.0) .. (dwd1);
        \draw[-Stealth, very thick] (dwd1) .. controls (20.5,2.0) .. (dwd14);
        \draw[-Stealth, very thick] (dwd14) .. controls (19.5,2.0) .. (dwd2);
        \draw[-Stealth, very thick] (dwd2) .. controls (18.5,2.0) .. (dwd15);
        \draw[-Stealth, very thick] (dwd15) .. controls (17.5,2.0) .. (dwd3);
        \draw[-Stealth, very thick] (dwd3) .. controls (16.5,2.0) .. (dwd16);
        \draw[-Stealth, very thick] (dwd16) .. controls (15.5,2.0) .. (dwd4);
        \draw[-Stealth, very thick] (dwd4) .. controls (14.5,2.0) .. (dwd17);
        \draw[-Stealth, very thick] (dwd17) .. controls (13.5,2.0) .. (dwd5);
        \draw[-Stealth, very thick] (dwd5) .. controls (12.5,2.0) .. (dwd18);
        \draw[-Stealth, very thick] (dwd18) .. controls (11.5,2.0) .. (dwd6);
        \draw[-Stealth, very thick] (dwd7) .. controls (10.5,2.0) .. (dwd18);
        \draw[-Stealth, very thick] (dwd19) .. controls (9.5,2.0) .. (dwd7);
        \draw[-Stealth, very thick] (dwd8) .. controls (8.5,2.0) .. (dwd19);
        \draw[-Stealth, very thick] (dwd20) .. controls (7.5,2.0) .. (dwd8);
        \draw[-Stealth, very thick] (dwd9) .. controls (6.5,2.0) .. (dwd20);
        \draw[-Stealth, very thick] (dwd21) .. controls (5.5,2.0) .. (dwd9);
        \draw[-Stealth, very thick] (dwd10) .. controls (4.5,2.0) .. (dwd21);
        \draw[-Stealth, very thick] (dwd22) .. controls (3.5,2.0) .. (dwd10);
        \draw[-Stealth, very thick] (dwd11) .. controls (2.5,2.0) .. (dwd22);
        \draw[-Stealth, very thick] (dwd23) .. controls (1.5,2.0) .. (dwd11);
        \draw[-Stealth, very thick] (dwd13) .. controls (22,3.0) .. (dwd24);
        \draw[-Stealth, very thick] (dwd24) .. controls (20.5,3.0) .. (dwd14);
        \draw[-Stealth, very thick] (dwd14) .. controls (19.5,3.0) .. (dwd25);
        \draw[-Stealth, very thick] (dwd25) .. controls (18.5,3.0) .. (dwd15);
        \draw[-Stealth, very thick] (dwd15) .. controls (17.5,3.0) .. (dwd26);
        \draw[-Stealth, very thick] (dwd26) .. controls (16.5,3.0) .. (dwd16);
        \draw[-Stealth, very thick] (dwd16) .. controls (15.5,3.0) .. (dwd27);
        \draw[-Stealth, very thick] (dwd27) .. controls (14.5,3.0) .. (dwd17);
        \draw[-Stealth, very thick] (dwd17) .. controls (13.5,3.0) .. (dwd28);
        \draw[-Stealth, very thick] (dwd28) .. controls (11.5,3.0) .. (dwd18);
        \draw[-Stealth, very thick] (dwd19) .. controls (9.5,3.0) .. (dwd28);
        \draw[-Stealth, very thick] (dwd29) .. controls (8.5,3.0) .. (dwd19);
        \draw[-Stealth, very thick] (dwd20) .. controls (7.5,3.0) .. (dwd29);
        \draw[-Stealth, very thick] (dwd30) .. controls (6.5,3.0) .. (dwd20);
        \draw[-Stealth, very thick] (dwd21) .. controls (5.5,3.0) .. (dwd30);
        \draw[-Stealth, very thick] (dwd31) .. controls (4.5,3.0) .. (dwd21);
        \draw[-Stealth, very thick] (dwd22) .. controls (3.5,3.0) .. (dwd31);
        \draw[-Stealth, very thick] (dwd32) .. controls (2.5,3.0) .. (dwd22);
        \draw[-Stealth, very thick] (dwd23) .. controls (1,3.0) .. (dwd32);
        \draw[-Stealth, very thick] (dwd33) .. controls (19.5,4.0) .. (dwd25);
        \draw[-Stealth, very thick] (dwd25) .. controls (18.5,4.0) .. (dwd34);
        \draw[-Stealth, very thick] (dwd34) .. controls (17.5,4.0) .. (dwd26);
        \draw[-Stealth, very thick] (dwd26) .. controls (16.5,4.0) .. (dwd35);
        \draw[-Stealth, very thick] (dwd35) .. controls (15.5,4.0) .. (dwd27);
        \draw[-Stealth, very thick] (dwd27) .. controls (14.5,4.0) .. (dwd36);
        \draw[-Stealth, very thick] (dwd36) .. controls (11.5,4.0) .. (dwd28);
        \draw[-Stealth, very thick] (dwd29) .. controls (8.5,4.0) .. (dwd36);
        \draw[-Stealth, very thick] (dwd37) .. controls (7.5,4.0) .. (dwd29);
        \draw[-Stealth, very thick] (dwd30) .. controls (6.5,4.0) .. (dwd37);
        \draw[-Stealth, very thick] (dwd38) .. controls (5.5,4.0) .. (dwd30);
        \draw[-Stealth, very thick] (dwd31) .. controls (4.5,4.0) .. (dwd38);
        \draw[-Stealth, very thick] (dwd39) .. controls (3.5,4.0) .. (dwd31);
        \draw[-Stealth, very thick] (dwd40) .. controls (18.5,5.0) .. (dwd34);
        \draw[-Stealth, very thick] (dwd34) .. controls (17.5,5.0) .. (dwd41);
        \draw[-Stealth, very thick] (dwd41) .. controls (16.5,5.0) .. (dwd35);
        \draw[-Stealth, very thick] (dwd35) .. controls (15.5,5.0) .. (dwd42);
        \draw[-Stealth, very thick] (dwd42) .. controls (11.5,5.0) .. (dwd36);
        \draw[-Stealth, very thick] (dwd37) .. controls (7.5,5.0) .. (dwd42);
        \draw[-Stealth, very thick] (dwd43) .. controls (6.5,5.0) .. (dwd37);
        \draw[-Stealth, very thick] (dwd38) .. controls (5.5,5.0) .. (dwd43);
        \draw[-Stealth, very thick] (dwd44) .. controls (4.5,5.0) .. (dwd38);
        \draw[-Stealth, very thick] (dwd39) .. controls (2,5.0) .. (dwd44);
        \draw[-Stealth, very thick] (dwd45) .. controls (17.5,6.0) .. (dwd41);
        \draw[-Stealth, very thick] (dwd41) .. controls (16.5,6.0) .. (dwd46);
        \draw[-Stealth, very thick] (dwd46) .. controls (11.5,6.0) .. (dwd42);
        \draw[-Stealth, very thick] (dwd43) .. controls (6.5,6.0) .. (dwd46);
        \draw[-Stealth, very thick] (dwd47) .. controls (5.5,6.0) .. (dwd43);
        \draw[-Stealth, very thick] (dwd44) .. controls (2.5,6.0) .. (dwd47);
        \node[vlabel, text=blue!45] (sdwd13) at (1,0.5) {$(6,1)$};
        \node[vlabel, text=brown!45] (sdwd23) at (12,0.5) {$(1,6)$};
        \node[vlabel, text=cyan!45] (sdwd39) at (16,0.5) {$(1,4)$};
        \node[vlabel, text=green!45] (sdwd40) at (5,0.5) {$(3,1)$};
        \node[vlabel, text=orange!45] (sdwd44) at (18,0.5) {$(1,3)$};
        \node[vlabel, text=red!45] (sdwd45) at (7,0.5) {$(2,1)$};
        \draw[-Stealth, very thick] (sdwd45) -- (dwd8);
        \draw[-Stealth, very thick] (dwd9) -- (sdwd45);
        \draw[-Stealth, very thick] (sdwd40) -- (dwd9);
        \draw[-Stealth, very thick] (dwd10) -- (sdwd40);
        \draw[-Stealth, very thick] (sdwd13) -- (dwd11);
        \draw[-Stealth, very thick] (dwd12) -- (sdwd13);
        \draw[-Stealth, very thick] (sdwd44) -- (dwd3);
        \draw[-Stealth, very thick] (dwd2) -- (sdwd44);
        \draw[-Stealth, very thick] (sdwd39) -- (dwd4);
        \draw[-Stealth, very thick] (dwd3) -- (sdwd39);
        \draw[-Stealth, very thick] (sdwd23) -- (dwd6);
        \draw[-Stealth, very thick] (dwd5) -- (sdwd23);
      \end{tikzpicture}
    }
  \endgroup
  \par\medskip
  \caption[The extended word $\bfi_\mathrm{GSV}^\gamma$ for the running example of GSV23]{%
    The extended word $\bfi_\mathrm{GSV}^\gamma$ for the running example of \cite{GSV23}, that is $\mcG=SL_7$ and $\gamma=(\{1\to4,2\to3,5\to1\},\{2\to3,3\to4,5\to6\})$, with chambers labeled by $\psi$ and $\tB_\gamma(\bfi_\mathrm{GSV})$ drawn as a quiver.
    The conventions are those of \cref{fig: Sp6 word seed}, with one color for each root of $\gamma$; compare with \cite[Figure~6]{GSV23}.}
  \label{fig: GSV23 word seed}
\end{figure}

We now show that the BD word seeds in type A belong to the cluster structure of \cite{GSV23}.
Let~$\mcG=SL_N$ with its standard pinning, let $\gamma$ be an aperiodic BD pair of type $A_{N-1}$, and let $\mcC_\gamma$ be the cluster structure of \cref{thm: GSV}.
Here $i^*=N-i$, $h_\gamma$ is birational by \cref{thm: GSV h}, $W=S_N$ acts on $[1,N]$, and $\Delta_{u\varpi_x,v\varpi_x}$ is the minor with rows $u([1,x])$ and columns $v([1,x])$; moreover
\begin{equation}
  \label{eq: pairing A}
  (w\varpi_x)(\alpha^\vee_g)=[g\in w([1,x])]-[g+1\in w([1,x])].
\end{equation}
We first identify one BD word seed with the initial seed of \cite{GSV23}.
The result for all words then follows from \cref{thm: moves}.
Following \cite{GSV23}, we describe skew-symmetric matrices by quivers: there are $b_{k\ell}$ arrows $k\to\ell$ if $b_{k\ell}>0$, and arrows between frozen vertices are not recorded.
The reader can follow along by comparing \cref{fig: GSV23 word seed} with \cite[Figure~6]{GSV23}.

\subsection{Identification of the initial seed}
\label{sec: base case}
Put $\mathbf e_c=(-c,\dots,-2,-1)$ and
\[
  \bfi_\mathrm{GSV}=(\mathbf e_1,\dots,\mathbf e_{N-1},\mathbf e_{N-1}^T,\dots,\mathbf e_1^T),
\]
so that $\bfi_\mathrm{GSV}^T=\bfi_\mathrm{GSV}$.
For $(i,j)\in[1,N]^2\setminus\{(1,1)\}$ let $F_{ij}=\Delta_{[i,i+s-1],[j,j+s-1]}$, where $s=N+1-\max(i,j)$.
We number the chambers at level $x$ by $j=0,\dots,2N-2x$ in increasing order of their slots.
Define
\begin{equation}
  \label{eq: psi}
  \psi_x(j)=
  \begin{cases}
    (j+1,N-x+1)&0\le j\le N-x,\\
    (N-x+1,2N-2x+1-j)&N-x\le j\le2N-2x.
  \end{cases}
\end{equation}
Thus each level traverses the hook $\max(i,j)=N-x+1$, and $\psi$ is a bijection onto $[1,N]^2\setminus\{(1,1)\}$, with initial and final chambers on the first row and column, respectively.
As $k\mapsto k^T$ maps the chamber numbered $j$ at level $x$ to the one numbered $2N-2x-j$, $\psi(k^T)=\psi(k)^T$.
\Cref{fig: GSV23 word seed} shows $\bfi_\mathrm{GSV}^\gamma$, $\tB_\gamma(\bfi_\mathrm{GSV})$ and $\psi$ for the running example of \cite{GSV23}.
The BD triples of \cite{GSV23} are, in our notation,
\[
  \mathbf\Gamma^\rr=(\Gamma_{\rr},\gamma_{\rr}(\Gamma_{\rr}),\gamma_{\rr}),\qquad
  \mathbf\Gamma^\cc=(\gamma_{\cc}(\Gamma_{\cc}),\Gamma_{\cc},\gamma_{\cc}^{-1}).
\]

\begin{proposition}
  \label{prop: base case}
  Via $\psi$, $\Sigma_\gamma(\bfi_\mathrm{GSV})$ is the initial seed of \cite[Theorem~4.11]{GSV23}.
\end{proposition}
\begin{proof}
  As $\bfi_\mathrm{GSV}^T=\bfi_\mathrm{GSV}$ and $\psi(k^T)=\psi(k)^T$, \cref{lem: mirror} gives $\Delta(k^T;\bfi_\mathrm{GSV})=\Delta(k;\bfi_\mathrm{GSV})\circ T$ and in~$\tB(\bfi_\mathrm{GSV})$ we have $b_{k^T\ell^T}=b_{k\ell}$.
  For the same reason, the inversion letters, thawed chambers, dipoles and virtual letters of the column roots of $\gamma$ are the transposes of those of the row roots of $\gamma^T$.

  \emph{Minors and tails.} We compute the chamber minors and their exponents, which reproduce the recursion of \cite{GSV23}.
  The block $\mathbf e_c$ moves the entry in position $c+1$ of the row permutation to position $1$.
  Consequently, for the chamber $k$ numbered $j\le N-x$ at level $x$, the row set is $[j+1,j+x]$ and the column set is $[N-x+1,N]$, so $\Delta(k;\bfi_\mathrm{GSV})=F_{\psi(k)}$, and this holds for all $k$ as $F_{ij}\circ T=F_{ji}$.
  Applying \eqref{eq: pairing A} to the rows and columns of $F_{ij}$ gives, for $\psi(k)=(i,j)$,
  \[
    e^\rr_a(k)=[i-j=\gamma_{\rr}(a)-N],\qquad
    e^\cc_b(k)=[i-j=N-\gamma_{\cc}(b)],
    \qquad
    \psi(t^\rr_a)=(a+1,1),\quad \psi(t^\cc_b)=(1,b+1).
  \]
  Thus \eqref{eq: Phi} reads $\Phi_\gamma(k;\bfi_\mathrm{GSV})=F_{ij}\pi_{i-j}$, where $\pi_0=1$ and
  \begin{equation}
    \label{eq: tails recursion}
    \pi_{g-N}=\begin{cases}F_{a+1,1}\pi_a&g=\gamma_{\rr}(a),\\1&g\notin\gamma_{\rr}(\Gamma_{\rr}),\end{cases}
    \qquad
    \pi_{N-g}=\begin{cases}F_{1,b+1}\pi_{-b}&g=\gamma_{\cc}(b),\\1&g\notin\gamma_{\cc}(\Gamma_{\cc}).
    \end{cases}
  \end{equation}
  This recursion terminates by \cref{lem: tails}, and by \cite[Theorem~4.4, (4.1), and (4.6)]{GSV23}, the cluster variables $f_{ij}$ of the initial seed of \cite[Theorem~4.11]{GSV23} satisfy $f_{ij}\circ h_\gamma=F_{ij}\pi_{i-j}$.
  So $x_\gamma(k;\bfi_\mathrm{GSV})=f_{\psi(k)}$.

  \emph{The exchange matrix.} By \eqref{eq: B entries}, a negative letter $i_p$ at level $x$ gives the arrow $c_x(p)\to c_x(p-1)$.
  Under~$\psi$ these are the arrows $(i,N-x+1)\to(i-1,N-x+1)$ along the hook, and those given by positive letters are their transposes.

  For arrows between adjacent levels, fix $1\le x\le N-2$.
  By \eqref{eq: B entries}, a consecutive pair $p<p'$ in $M_{x,x+1}$ with $\varphi_{x,x+1}$-values $(1,-1)$ gives $c_{x+1}(p)\to c_x(p)$; values $(-1,1)$ give the reverse arrow, and equal values give none.
  By transposition, we need only consider pairs in the negative half and the pair crossing the middle of the word.

  At these two levels, the negative half consists of $-x$ in $\mathbf e_x$, followed by $-(x+1),-x$ in each block $\mathbf e_c$ for $x+1\le c\le N-1$.
  Between $\mathbf e_{c-1}$ and $\mathbf e_c$, the consecutive letters $-x,-(x+1)$ have values $(1,-1)$; applying \eqref{eq: psi} gives
  \[
    (c-x,N-x)\longrightarrow(c-x+1,N-x+1).
  \]
  Within $\mathbf e_c$, the consecutive letters $-(x+1),-x$ have values $(-1,1)$ and give
  \[
    (c-x+1,N-x+1)\longrightarrow(c-x+1,N-x).
  \]
  Finally, the last $-x$ in $\mathbf e_{N-1}$ and the first $x$ in $\mathbf e_{N-1}^T$ are consecutive in $M_{x,x+1}$, with values $(1,-1)$, and give
  \[
    (N-x,N-x)\longrightarrow(N-x+1,N-x+1).
  \]
  These arrows and their transposes are exactly the arrows between the hooks $N-x+1$ and $N-x$ in $Q_{\varnothing,\varnothing}$ \cite[Section~4.3 and Figure~5]{GSV23}.

  The row inversion for $g=\gamma_{\rr}(a)$ is the last letter of $\mathbf e_g$, at level $1$: the letter $-x$ of $\mathbf e_c$ swaps the values~$c-x+1$ and $c+1$, which are consecutive precisely when $x=1$.
  Its adjacent chambers map to~$(g,N)$ and $(g+1,N)$.
  So the thawing and dipoles of the row roots, and by transposition those of the column roots, give exactly the GSV modifications
  \begin{equation}
    \label{eq: GSV arrows}
    (\gamma_{\rr}(a),N)\to(a+1,1)\to(\gamma_{\rr}(a)+1,N),\qquad
    (N,\gamma_{\cc}(b))\to(1,b+1)\to(N,\gamma_{\cc}(b)+1).
  \end{equation}

  By \cref{lem: local entries}, it remains to check the entries between end chambers.
  On a connected component $[d,d'-1]\subset\Gamma_{\rr}$, adjacent virtual letters occur in decreasing order of their $\gamma_{\rr}$-images.
  For $a,a+1$ in this component, the last real letter and the two virtual letters on levels $a^*-1,a^*$ have $\varphi$-values $-1,+1,-1$ when $\gamma_{\rr}$ preserves orientation, and $-1,-1,+1$ when it reverses orientation.
  By \eqref{eq: B entries}, these give, respectively, the arrow $(a+2,1)\to(a+1,1)$ and no arrow between the old end chambers.
  At the lower end of the component the values $-1,+1$ give $(d+1,1)\to(d,1)$; at the upper end the values $-1,-1$ give no arrow.
  Thus the added first-column arrows are
  \[
    \begin{cases}
      (a+2,1)\to(a+1,1),\quad d-1\le a\le d'-2,&\gamma_{\rr}\text{ preserves orientation},\\
      (d+1,1)\to(d,1),&\gamma_{\rr}\text{ reverses orientation}.
    \end{cases}
  \]
  Arrows incident with $(1,1)$ are omitted.
  The column virtual letters give the transposes, and no initial and final chambers share a slot, so the extended word adds no entries between them.
  These are precisely the remaining modifications turning $Q_{\varnothing,\varnothing}$ into the quiver $Q_{\mathbf\Gamma^\rr,\mathbf\Gamma^\cc}$ of \cite[Section~4.3]{GSV23}.
  This completes the seed identification.
\end{proof}

\begin{corollary}
  \label{cor: type A seeds}
  For every double reduced word $\bfi$, $\Sigma_\gamma(\bfi)$ is a seed of $\mcC_\gamma$.
  Commutations preserve the seed, and sign swaps and braid moves mutate it at the replaced chamber, with the chamber correspondences of \cref{sec: word seeds}.
\end{corollary}
\begin{proof}
  Every double reduced word is connected to $\bfi_\mathrm{GSV}$ by these moves.
  The claim follows by applying \cref{thm: moves}, using birationality of $h_\gamma$ and the fact that type $A$ is simply-laced.
\end{proof}

\subsection{The global toric action}
\label{sec: toric A}
We now prove \cref{thm: GSV}(iii) for arbitrary aperiodic pairs, using only parts (i), (ii) and the seed just identified.
Equivariance gives homogeneity of the initial cluster.
Regularity and completeness force homogeneity of exchange relations.
Let $y_{\rr}\in\mcH_{\gamma_{\rr}}$ and $y_{\cc}\in\mcH_{\gamma_{\cc}}$.
As $\mcH$ normalizes $\mcN_\pm^{\Gamma_\bullet}$ and lies in $\mcB_\pm$, the products $y_{\rr}Uy_{\cc}=\big(y_{\rr}V_{\rr}(U)y_{\rr}^{-1}\big)\big(y_{\rr}V_{\rr}(U)^{-1}Uy_{\cc}\big)=\big(y_{\rr}UV_{\cc}(U)^{-1}y_{\cc}\big)\big(y_{\cc}^{-1}V_{\cc}(U)y_{\cc}\big)$ are again factorizations as in \cref{sec: h}, so
\[
  V_{\rr}(y_{\rr}Uy_{\cc})=y_{\rr}V_{\rr}(U)y_{\rr}^{-1},\qquad V_{\cc}(y_{\rr}Uy_{\cc})=y_{\cc}^{-1}V_{\cc}(U)y_{\cc},
\]
and the same formulas for $\overline{V}_{\rr},\overline{V}_{\cc}$.
The defining character equations of these tori imply that $\boldsymbol\gamma_{\rr}$ and $\boldsymbol\gamma_{\cc}$ commute with the respective conjugations.
The formulas for $H_{\rr},H_{\cc}$ therefore give $h_\gamma(y_{\rr}Uy_{\cc})=y_{\rr}h_\gamma(U)y_{\cc}$.
For any double reduced word $\bfi$, $\Phi_\gamma(k;\bfi)$ is a monomial in generalized minors, so birationality and equivariance of $h_\gamma$ imply that $x_\gamma(k;\bfi)$ is $\mcH_\gamma$-homogeneous.
In a basis of the character lattice of $\mcH_\gamma$, the weights of the initial cluster variables are the rows of a matrix $P$, whose local toric action (\cref{sec: cluster conventions}) is the action of $\mcH_\gamma$.

To prove $\tB^TP=0$, suppose the two monomials in an exchange relation $x_kx'_k=M_++M_-$ have different weights.
Regularity of $x'_k$ and the weight decomposition of $\bbC[SL_N]$ then imply that $M_+/x_k$ is regular.
By completeness, it belongs to the Laurent ring of the adjacent cluster, where
\[
  \frac{M_+}{x_k}=\frac{M_+x'_k}{M_++M_-}.
\]
The exchange column is nonzero by full rank, so $M_++M_-$ is a nonunit binomial in that ring.
Its monomials have disjoint supports, so in that factorial ring it is coprime to every monomial.
As $x'_k$ and the exchangeable variables are units there, $M_+x'_k$ is a unit times a monomial in the frozen variables, and so are its divisors.
So the displayed fraction does not belong to the ring, a contradiction.
Thus $\tB^TP=0$, and the action is global.

Finally, as the initial cluster generates $\bbC(SL_N)$, the rank of $P$ is the dimension of the image of $\mcH_\gamma$ in the automorphisms of $\bbC(SL_N)$.
If $(y_{\rr},y_{\cc})$ acts trivially, then $y_{\rr}Uy_{\cc}=U$ for all $U$, so $y_{\cc}=y_{\rr}^{-1}$ and $y_{\rr}$ is central; hence this image has dimension $\dim\mcH_\gamma=2(N-1)-|\Gamma_{\rr}|-|\Gamma_{\cc}|$.
This is the number of frozen variables, the largest possible rank (\cref{sec: cluster conventions}), so the global toric action is maximal.

\section{The symplectic group and its BD pairs}
\label{sec: Sp}
\subsection{The symplectic group as a fixed point set}
\label{sec: fixed points}
For $n=1$, $Sp_2=SL_2$ and the only BD pair is $(0,0)$.
From now on $n\ge2$, $N=2n$, $\tG=\SL$ with the standard pinning, and $\widetilde{\mfg}=\mathfrak{sl}_{2n}$.
Let $\tau(k):=2n-k$ on~$[1,2n-1]$; we also write $\tau$ for the induced involutions of the simple roots, of the weights, and of~$S_{2n}$ (conjugation by $w_0$).
Thus $\tau(k)=k^*$ in type $A_{2n-1}$; we retain the notation $\tau$ for this involution to distinguish it from $*$ in type $C_n$, where $k^*=k$.
Let
\[
  J:=\sum_{i=1}^{N}(-1)^{i+1}E_{i,N+1-i},
  \qquad
  \vt(U):=J\,(U^T)^{-1}J^{-1}.
\]
The alternating signs are chosen so that $\vt$ preserves the standard pinning.
Then $J^T=-J=J^{-1}$, $\vt$ is an involutive automorphism of $\tG$, and its fixed point set is $\mcG:=\{U\mid UJU^T=J\}=\Sp$.
Its differential is $\Vt(\xi)=-J\xi^TJ^{-1}$ for $\xi\in\widetilde{\mfg}$.

\begin{lemma}
  \label{lem: pinning}
  $\Vt(E_k)=E_{\tau(k)}$, $\Vt(F_k)=F_{\tau(k)}$ and $\Vt(\alpha_k^\vee)=\alpha_{\tau(k)}^\vee$.
  Hence $\vt$ preserves $\widetilde{\mcH}$, $\widetilde{\mcB}_\pm$, $\widetilde{\mcN}_\pm$, commutes with the transpose and satisfies $\vt(\overline{w})=\overline{\tau(w)}$ for $w\in S_{2n}$, and $\Vt$ preserves the trace form.
\end{lemma}
\begin{proof}
  We have $JE_{ij}J^{-1}=(-1)^{i+j}E_{w_0(i)w_0(j)}$, so $\Vt(E_k)=-JE_{k+1,k}J^{-1}=E_{N-k,N-k+1}=E_{\tau(k)}$, and similarly for $F_k$.
  Since $\overline{s_k}=\exp(-E_k)\exp(F_k)\exp(-E_k)$, we have $\vt(\overline{s_k})=\overline{s_{\tau(k)}}$.
  Moreover we have~$\vt(U)^T=J^{-T}U^{-1}J^T=JU^{-1}J^{-1}=\vt(U^T)$.
  Also $\operatorname{tr}(\Vt(\xi)\Vt(\eta))=\operatorname{tr}(J\xi^T\eta^TJ^{-1})=\operatorname{tr}(\xi\eta)$, so $\Vt$ preserves the trace form.
\end{proof}

As $\mcG$ is the fixed point set of $\vt$, $\mfg=\mathrm{Lie}(\mcG)=\widetilde{\mfg}^\Vt$; let $\mfp$ be the $(-1)$-eigenspace of $\Vt$.
As $\Vt$ preserves the trace form, the decomposition $\widetilde{\mfg}=\mfg\oplus\mfp$ is orthogonal for it.
We equip $\mfg$ with the restricted trace form and write $\pr:=\frac12(1+\Vt)$ for the orthogonal projection onto $\mfg$.
The group $\mcG$ is simply connected of type $C_n$, with maximal torus $\mcH=\widetilde{\mcH}^\vt$, Cartan subalgebra $\mfh=\widetilde{\mfh}^\Vt$, $\mcB_\pm=\widetilde{\mcB}_\pm\cap\mcG$, $\mcN_\pm=\widetilde{\mcN}_\pm\cap\mcG$, and Chevalley generators
\[
  E^C_k:=E_k+E_{\tau(k)},\quad F^C_k:=F_k+F_{\tau(k)}\quad(k<n),
  \qquad
  E^C_n:=E_n,\quad F^C_n:=F_n.
\]
The simple roots are $\alpha^C_k:=\alpha_k|_\mfh=\alpha_{\tau(k)}|_\mfh$, $k\in[1,n]$, the long one being $\alpha^C_n$, and the fundamental weights are $\varpi^C_k=\varpi_k|_\mfh=\varpi_{\tau(k)}|_\mfh$.
The Weyl group of $\mcG$ is $W^\tau=\{w\in S_{2n}\mid ww_0=w_0w\}$, with simple reflections $s^C_k=s_ks_{\tau(k)}$ ($k<n$) and $s^C_n=s_n$, and longest element $w_0$.
Since $E_k,F_k$ commute with $E_{\tau(k)},F_{\tau(k)}$ for $k<n$, the representative of $s^C_k$ in $\mcG$ is $\overline{s_k}\,\overline{s_{\tau(k)}}$, so the representative of $w\in W^\tau$ in $\mcG$ is its representative $\overline{w}$ in $\tG$.
The transpose of $\tG$ preserves $\mcG$ and restricts to the transpose of $\mcG$.

\begin{lemma}
  \label{lem: minors restrict}
  Let $x\in[1,2n-1]$ and $x'=\min(x,\tau(x))$.
  For $u,v\in W^\tau$, $\Delta_{u\varpi_x,v\varpi_x}|_\mcG=\Delta_{u\varpi^C_{x'},v\varpi^C_{x'}}$, and for $u,v\in S_{2n}$, $\Delta_{u\varpi_x,v\varpi_x}\circ\vt=\Delta_{\tau(u)\varpi_{\tau(x)},\tau(v)\varpi_{\tau(x)}}$.
\end{lemma}
\begin{proof}
  The Gauss factors of an element of $\mcG\cap\widetilde{\mcN}_-\widetilde{\mcH}\widetilde{\mcN}_+$ lie in $\mcG$ by uniqueness and \cref{lem: pinning}, and for $y=\mathrm{diag}(t_1,\dots,t_n,t_n^{-1},\dots,t_1^{-1})\in\mcH$ we have $y^{\varpi_x}=t_1\cdots t_{x'}=y^{\varpi^C_{x'}}$.
  This gives the first claim for~$u=v=e$, and the general case follows as $u,v$ have the same representatives in $\mcG$ and $\tG$.
  To prove the second claim, write $U=U_-U_0U_+$ with $U_0=\mathrm{diag}(t_1,\dots,t_N)$.
  Then $\vt(U)=\vt(U_-)\vt(U_0)\vt(U_+)$ is the Gauss decomposition of $\vt(U)$, and $\vt(U_0)=\mathrm{diag}(t_N^{-1},\dots,t_1^{-1})$.
  Since $t_1\cdots t_N=1$, we obtain
  \[
    \vt(U_0)^{\varpi_x}=(t_{N-x+1}\cdots t_N)^{-1}=t_1\cdots t_{N-x}=U_0^{\varpi_{\tau(x)}}.
  \]
  This proves the case $u=v=e$.
  The general case follows from $\vt(\overline{u})=\overline{\tau(u)}$.
\end{proof}

\subsection{BD maps of type \texorpdfstring{$C_n$}{Cn} and their lifts}
We lift a symplectic Belavin--Drinfeld map by reflecting it across the middle of the type A diagram.
Lifting both maps gives a pair that commutes with the involution and has the same aperiodicity condition.

A BD map of type $C_n$ does not involve $\alpha^C_n$: as an isometry, it could only map the unique long simple root to itself, which nilpotency forbids, and similarly for $\gamma_\bullet^{-1}$.
So it is a BD map of type $A_{n-1}$ on~$\{\alpha^C_1,\dots,\alpha^C_{n-1}\}$, identified as before with $[1,n-1]$.

\begin{definition}
  \label{def: lift}
  The \newword{lift} of a BD map $\gamma_\bullet:\Gamma_\bullet\to[1,n-1]$ of type $C_n$ is the map $\tg_\bullet$ on $\Gamma_\bullet\sqcup\tau(\Gamma_\bullet)\subseteq[1,2n-1]$ given by $\tg_\bullet(k)=\gamma_\bullet(k)$ and $\tg_\bullet(\tau(k))=\tau(\gamma_\bullet(k))$.
  The lift of a BD pair $\gamma=(\gamma_{\rr},\gamma_{\cc})$ is $\tg:=(\tg_{\rr},\tg_{\cc})$.
\end{definition}

For example, in $Sp_6$ the choice $\gamma_{\rr}=\{1\to2\}$ gives $\tg_{\rr}=\{1\to2,5\to4\}$.

\begin{proposition}
  \label{prop: lift}
  Let $\gamma_\bullet$ be a BD map of type $C_n$.
  \begin{enumerate}[(a)]
    \item $\tg_\bullet$ is a BD map of type $A_{2n-1}$ commuting with $\tau$.
    \item As endomorphisms of $\widetilde{\mfg}$, $\tg_\bullet\,\Vt=\Vt\,\tg_\bullet$, and $\tg_\bullet|_\mfg=\gamma_\bullet$; in particular $C_{\tg_\bullet}$ commutes with $\pr$ on $\widetilde{\mfn}_+\oplus\widetilde{\mfn}_-$ and restricts to $C_{\gamma_\bullet}$ on $\mfn_+\oplus\mfn_-$.
    \item $\widetilde{\mfh}_{\tg_\bullet}\cap\mfh=\mfh_{\gamma_\bullet}$ and $\widetilde{\mcH}_{\tg_\bullet}\cap\mcG=\mcH_{\gamma_\bullet}$.
  \end{enumerate}
  A BD pair $\gamma$ of type $C_n$ is aperiodic if and only if its lift $\tg$ is.
\end{proposition}
\begin{proof}
  (a) Both $\Gamma_\bullet$ and $\gamma_\bullet(\Gamma_\bullet)$ lie in $[1,n-1]$, while their $\tau$-images lie in $[n+1,2n-1]$.
  Thus the two halves of the domain are not adjacent, nor are their images, so $\tg_\bullet$ is an isometry.
  It is nilpotent since~$\tg_\bullet^m(\tau(k))=\tau(\gamma_\bullet^m(k))$.

  (b) The subalgebras $\widetilde{\mfg}_{\Gamma_\bullet}$ and $\widetilde{\mfg}_{\tau(\Gamma_\bullet)}=\Vt(\widetilde{\mfg}_{\Gamma_\bullet})$ commute, $\Vt\tg_\bullet\Vt$ and $\tg_\bullet$ agree on their Chevalley generators by \cref{lem: pinning}, and the orthogonal projection onto their sum commutes with the isometry $\Vt$.
  So $\tg_\bullet\Vt=\Vt\tg_\bullet$.
  The fixed points of $\Vt$ in $\widetilde{\mfg}_{\Gamma_\bullet}\oplus\widetilde{\mfg}_{\tau(\Gamma_\bullet)}$ form the subalgebra of $\mfg$ generated by the $E^C_k,F^C_k$, for $k\in\Gamma_\bullet$, on which $\tg_\bullet$ acts by $E^C_k\mapsto E^C_{\gamma_\bullet(k)}$, $F^C_k\mapsto F^C_{\gamma_\bullet(k)}$.
  The orthogonal projection of $\mfg$ onto it is the restriction of the projection of $\widetilde{\mfg}$ onto $\widetilde{\mfg}_{\Gamma_\bullet}\oplus\widetilde{\mfg}_{\tau(\Gamma_\bullet)}$; hence $\tg_\bullet|_\mfg=\gamma_\bullet$.

  (c) On $\mfh$ the conditions $\alpha(\zeta)=\tg_\bullet(\alpha)(\zeta)$ for $\alpha=\alpha_k$ and $\alpha=\alpha_{\tau(k)}$ both read $\alpha^C_k(\zeta)=\gamma_\bullet(\alpha^C_k)(\zeta)$.
  For the tori, write $y=\mathrm{diag}(t_1,\dots,t_n,t_n^{-1},\dots,t_1^{-1})\in\mcH$ and use the coordinates $q_i=t_i/t_{i+1}$ for $i<n$, and~$q_n=t_n$, with inverse $t_i=q_i\cdots q_n$.
  The character equations for $\gamma_\bullet$, and those for $\tg_\bullet$ restricted to $\mcH$, both become $q_i=q_{\gamma_\bullet(i)}$ for $i\in\Gamma_\bullet$.
  Their full solution subgroup in $\mcH$ is connected, since these equations merely identify coordinates, and hence equals $\mcH_{\gamma_\bullet}$.
  By connectedness it lies in the identity component~$\widetilde{\mcH}_{\tg_\bullet}$ of the ambient solution subgroup.
  Conversely, every element of $\widetilde{\mcH}_{\tg_\bullet}\cap\mcG$ satisfies these equations.
  This proves the torus equality.

  Finally, $-w_0$ acts as $\tau$ on the simple roots of $\SL$ and trivially on those of $\mcG$.
  So the map~${\alpha_k,\alpha_{\tau(k)}\mapsto\alpha^C_k}$ intertwines the alternating maps of $\tg$ and of $\gamma$, and the domain of the former is the preimage of the domain of the latter.
  A periodic point of the former maps to a periodic point of the latter.
  Conversely, if $x$ is periodic for the latter, the former maps the finite preimage of the orbit of $x$ injectively, hence bijectively, into itself, so it has periodic points.
\end{proof}

\begin{remark}
  \label{rmk: aperiodic not PL}
  A nontrivial aperiodic BD pair of type $C_n$ does not give a Poisson--Lie bracket on $\Sp$.
  Indeed, for a Poisson--Lie bracket every $\{f,g\}$ vanishes at the identity, whereas at the identity $\{\cdot,\cdot\}_{\gamma,s}$ is given by $r_\cc-r_\rr$, whose non-Cartan part vanishes only if $\gamma_{\cc}=\gamma_{\rr}^{-1}$.
  Since $-w_0$ acts trivially on the simple roots of type $C_n$, every $\alpha\in\Gamma_{\rr}$ then gives the alternating two-cycle $\alpha^\rr\mapsto\gamma_{\rr}(\alpha)^\cc\mapsto\alpha^\rr$.
  Thus an aperiodic pair can give a Poisson--Lie bracket only if $\gamma_{\rr}=\gamma_{\cc}=0$.
\end{remark}

\subsection{Cartan data and restriction of brackets}
Let $\mfk$ be the $(-1)$-eigenspace of $\Vt$ on $\widetilde{\mfh}$, so $\widetilde{\mfh}=\mfh\oplus\mfk$ is an orthogonal decomposition.
For a $\Vt$-invariant $\widetilde{s}\in\widetilde{\mfh}\wedge\widetilde{\mfh}$, put $s:=(\pr\otimes\pr)\widetilde{s}\in\mfh\wedge\mfh$; then $\widetilde{s}-s\in\mfk\wedge\mfk$.

\begin{proposition}
  \label{prop: cartan}
  For a BD map $\gamma_\bullet$ of type $C_n$, the map $\widetilde{s}\mapsto s$ sends the $\Vt$-invariant elements of $\mcS(\tg_\bullet)$ onto $\mcS(\gamma_\bullet)$.
\end{proposition}
\begin{proof}

  For a $\Vt$-invariant $\widetilde{s}$ and $\lambda\in\widetilde{\mfh}^*$, the $\mfh$-component of $2(\lambda\otimes1)\widetilde{s}$ is $2(\lambda|_\mfh\otimes1)s$.
  For $\alpha=\alpha_k$ or $\alpha_{\tau(k)}$ with $k\in\Gamma_\bullet$, we have $(\alpha-\tg_\bullet(\alpha))|_\mfh=\alpha^C_k-\gamma_\bullet(\alpha^C_k)$ and $\pr\big(((\alpha+\tg_\bullet(\alpha))\otimes1)\Cas_0\big)=\big((\alpha^C_k+\gamma_\bullet(\alpha^C_k))\otimes1\big)\Cas^\mfg_0$, as $\widetilde{\mfh}=\mfh\oplus\mfk$ is orthogonal, $\Cas^\mfg_0$ being the Cartan part of the Casimir element of $\mfg$ for the restricted form.
  So the $\mfh$-components of the equations \eqref{eq: system s} for $\tg_\bullet$ are the equations for $\gamma_\bullet$, and the $\mfk$-components only involve $\widetilde{s}-s$; hence $s\in\mcS(\gamma_\bullet)$.
  Conversely, $\mcS(\tg_\bullet)$ is a non-empty affine space stable under $\Vt\otimes\Vt$, so it contains the $\Vt$-invariant element $\widetilde{s}^0:=\frac12(\widetilde{s}+(\Vt\otimes\Vt)\widetilde{s})$ for any $\widetilde{s}\in\mcS(\tg_\bullet)$.
  Replacing its $\mfh\wedge\mfh$ component by any $s\in\mcS(\gamma_\bullet)$ gives the invariant solution $\widetilde{s}^0+s-(\pr\otimes\pr)\widetilde{s}^0$, which projects to $s$.
\end{proof}

For a BD pair $\gamma$ of type $C_n$ we always take $\widetilde{s}=(\widetilde{s}_{\rr},\widetilde{s}_{\cc})\in\mcS(\tg)$ with $\Vt$-invariant components, and~$s=(s_{\rr},s_{\cc})\in\mcS(\gamma)$ with $s_\bullet:=(\pr\otimes\pr)\widetilde{s}_\bullet$.
Every $s\in\mcS(\gamma)$ arises in this way.

\begin{proposition}
  \label{prop: restriction}
  \begin{enumerate}[(a)]
    \item $\vt$ is a Poisson automorphism of $(\tG,\{\cdot,\cdot\}_{\tg,\widetilde{s}})$.
    \item For $\vt$-invariant $f,g\in\bbC[\tG]$, $\{f,g\}_{\tg,\widetilde{s}}|_\mcG=\{f|_\mcG,g|_\mcG\}_{\gamma,s}$.
  \end{enumerate}
\end{proposition}
\begin{proof}
  (a) Since $\partial^\rr_\xi(f\circ\vt)=(\partial^\rr_{\Vt\xi}f)\circ\vt$, and similarly for $\partial^\cc$, it suffices that $\Vt\otimes\Vt$ fixes $r_\rr$ and $r_\cc$.
  This holds because $\Vt$ preserves the trace form, hence $\Cas$, preserves $\widetilde{\mfn}_\pm$, hence $\Cas_{-+}$ and $\Cas_{+-}$, commutes with~$C_{\tg_{\rr}}$,~$C_{\tg_{\cc}}$ and with the transpose, and fixes $\widetilde{s}_{\rr},\widetilde{s}_{\cc}$.

  (b) Let $U\in\mcG$.
  For any $\xi\in\widetilde{\mfg}$ we have $\partial^\rr_\xi f(U)=\partial^\rr_{\Vt\xi}f(U)$ because $f=f\circ\vt$ and $\vt(U)=U$, hence~$\partial^\rr_\xi f(U)=\partial^\rr_{\pr\xi}f(U)$.
  For $\xi\in\mfg$ the curve $\exp(t\xi)U$ lies in $\mcG$, hence $\partial^\rr_\xi f(U)=\partial^\rr_\xi(f|_\mcG)(U)$.
  The same holds for $\partial^\cc$.
  For any tensor $\Psi=\sum_u\xi_u\otimes\eta_u\in\widetilde{\mfg}\otimes\widetilde{\mfg}$, we therefore have
  \[
    \Psi(f,g)^\rr(U)=\sum_u\partial^\rr_{\pr\xi_u}(f|_\mcG)(U)\,\partial^\rr_{\pr\eta_u}(g|_\mcG)(U)=((\pr\otimes\pr)\Psi)(f|_\mcG,g|_\mcG)^\rr(U),
  \]
  and similarly for column derivatives.
  So it suffices to show that ${(\pr\otimes\pr)r_{\tg_\bullet,\widetilde{s}_\bullet}=r^\mfg_{\gamma_\bullet,s_\bullet}}$, the $R$-matrix of~$\mfg$, for each of the two BD maps (the transpose commutes with $\pr$).
  By orthogonality the Casimir element of $\widetilde{\mfg}$ is $\Cas^\mfg+\Cas^\mfp$, so $(\pr\otimes\pr)\Cas=\Cas^\mfg$, and $\pr$ maps the components $\Cas_0,\Cas_{-+},\Cas_{+-}$ to those of~$\Cas^\mfg$, as it preserves $\widetilde{\mfh}$ and $\widetilde{\mfn}_\pm$.
  By \cref{prop: lift}(b), $(\pr\otimes\pr)(1\otimes C_{\tg_\bullet})\Cas_{-+}=(1\otimes C_{\gamma_\bullet})\Cas^\mfg_{-+}$, and similarly for~$\Cas_{+-}$; finally $(\pr\otimes\pr)\widetilde{s}_\bullet=s_\bullet$.
\end{proof}

\begin{remark}
  \label{rmk: not Poisson subvariety}
  $\mcG$ is not a Poisson subvariety of $(\tG,\{\cdot,\cdot\}_{\tg,\widetilde{s}})$, already for $\gamma=(0,0)$.
  For $n=2$ and $\widetilde{s}=0$, write $U=(x_{ij})$ and
  \[
    q:=(UJU^T-J)_{13}=x_{11}x_{34}-x_{12}x_{33}+x_{13}x_{32}-x_{14}x_{31}.
  \]
  Then $q|_{Sp_4}=0$, whereas $\{x_{11},q\}_{(0,0),0}(I+E_{12}+E_{34})=-\tfrac12$, with $I+E_{12}+E_{34}\in Sp_4$, so the defining ideal of $Sp_4$ is not a Poisson ideal.
  This failure is independent of the chosen compatible pinnings.
  Indeed, changing the pinning transports $\vt$, its fixed subgroup and the associated Poisson bracket by the same automorphism of $\tG$.
  \Cref{prop: restriction} says instead that the restriction $\bbC[\tG]^\vt\to\bbC[\mcG]$, which is surjective, is a Poisson morphism.
\end{remark}

\subsection{Lifted words}
We lift symplectic words by replacing each short-root letter with its commuting pair in type A. To apply the ambient word-seed construction, we must check that the lifted word is reduced.

A double reduced word of type $C_n$ is a shuffle $\bfi$ of a reduced word for $w_0\in W^\tau$ in the alphabet $[-n,-1]$ and one in the alphabet $[1,n]$.
Two such words are related by commutations, sign swaps, braid moves of length $3$, $(k,k+1,k)\leftrightarrow(k+1,k,k+1)$ with $k+1<n$, and of length $4$, $(n-1,n,n-1,n)\leftrightarrow(n,n-1,n,n-1)$.

\begin{definition}
  \label{def: lifted word}
  The \newword{lift} $\widetilde{\boldsymbol{\jmath}}$ of an arbitrary word $\bfj$ in $\pm[1,n]$ is obtained by replacing each letter $\pm k$, with~$k<n$, by the two letters $\pm k,\pm\tau(k)$, and keeping each $\pm n$.
\end{definition}

\begin{lemma}
  \label{lem: lifted word}
  If $\bfi$ is a double reduced word of type $C_n$, its lift $\tbfi$ is a double reduced word for $(w_0,w_0)$ in~$\SL$, and replacing every letter $\pm k$ of $\tbfi$ by $\pm\tau(k)$ gives a word related to $\tbfi$ by commutations.
\end{lemma}
\begin{proof}
  Braid moves preserve the number of occurrences of the long-root reflection $s_n^C$.
  Comparing with the reduced expression $(s_1^C\cdots s_n^C)^n$, every reduced expression for $w_0$ in $W^\tau$ therefore has $n$ long-root factors and $n^2-n$ short-root factors.
  Lifting keeps each long-root factor and replaces each short-root factor by two simple reflections.
  Since $s^C_k=s_ks_{\tau(k)}$ in $S_{2n}$, the lifted expression still has product $w_0$, and its length is $n+2(n^2-n)=\binom{2n}{2}$, so it is reduced.
  Applying this to the positive and negative subwords proves the first claim.
  For the second, $\tau$ fixes the letters $\pm n$ and swaps the two letters in each lifted pair; these commute because $k$ and $\tau(k)$ are not adjacent for $k<n$.
\end{proof}

\subsection{The map \texorpdfstring{$h_\gamma$}{h}}
Let $\gamma$ be an aperiodic BD pair of type $C_n$ and $\tg$ its lift.

\begin{proposition}
  \label{prop: h theta}
  $h_{\tg}\circ\vt=\vt\circ h_{\tg}$, and $h_{\tg}$ is defined on a dense open subset of $\mcG$, which it maps into $\mcG$; on $\mcG$ it agrees with the map $h_\gamma$ of \eqref{eq: h} for $\mcG$ and $\gamma$.
\end{proposition}
\begin{proof}
  We use the notation from the definition of $h_{\tg}$ in \cref{sec: h}; in particular we write $\sigma_{\rr}=w_0^{\widetilde{\Gamma}_{\rr}}w_0$ and $\sigma_{\cc}=w_0^{\widetilde{\Gamma}_{\cc}}w_0$.
  The sets $\widetilde{\Gamma}_{\rr}=\Gamma_{\rr}\sqcup\tau(\Gamma_{\rr})$ and $\widetilde{\Gamma}_{\cc}=\Gamma_{\cc}\sqcup\tau(\Gamma_{\cc})$ are $\tau$-stable, so $w_0^{\widetilde{\Gamma}_{\rr}}$, $w_0^{\widetilde{\Gamma}_{\cc}}$, $\sigma_{\rr}$, $\sigma_{\cc}$ lie in~$W^\tau$ and their representatives are fixed by $\vt$.
  Moreover, $\vt$ preserves $\widetilde{\mcN}_\pm^{\widetilde{\Gamma}_{\rr}}$, $\widetilde{\mcN}_\pm^{\widetilde{\Gamma}_{\cc}}$ and the cells $\widetilde{\mcB}_\pm\overline{w}\widetilde{\mcB}_\pm$ for~$w\in W^\tau$, and it commutes with $\boldsymbol\gamma_{\rr},\boldsymbol\gamma_{\cc}$ by \cref{prop: lift}(b).
  By uniqueness of the factorizations in \eqref{eq: h}, $V_{\rr}(\vt(U))=\vt(V_{\rr}(U))$, $\overline{V}_{\rr}(\vt(U))=\vt(\overline{V}_{\rr}(U))$, and similarly for columns; so $h_{\tg}\circ\vt=\vt\circ h_{\tg}$.
  In particular $h_{\tg}$ maps the fixed points of $\vt$ in its domain into $\mcG$.
  It remains to show that the domains of $H_{\rr}$ and $H_{\cc}$ each meet $\mcG$.
  Their intersections with $\mcG$ are then nonempty open subsets, and irreducibility of $\mcG$ ensures that their intersection is nonempty and dense, with $h_{\tg}$ defined there by \eqref{eq: h}.
  By \cref{lem: h transpose}, as the transpose preserves $\mcG$, it suffices to treat $H_{\rr}$; write $K=\widetilde{\Gamma}_{\rr}$.

  Put $y_k(t):=\exp(tF_k)$, $y^C_k:=y_k y_{\tau(k)}$ for $k<n$, and $y^C_n:=y_n$.
  Choose a reduced word $(j_1,\dots,j_{\len(w_0)})$ for $w_0$ in $W^\tau$ whose first $\len(w_0^K)$ letters form a reduced word for $w_0^K\in W^\tau$, $\len$ being the length function of $W^\tau$.
  As in \cref{lem: lifted word}, its lift is a reduced word $\widetilde{\boldsymbol{\jmath}}$ for $w_0$ in $S_{2n}$ beginning with one for $w_0^K$.
  By \cite[Proposition~2.18]{FZ99} (transposed), the map $(t_m)\mapsto\prod_m y_{\widetilde{\jmath}_m}(t_m)$ is an isomorphism from a torus onto an open subset $\mathcal{U}_-$ of $\widetilde{\mcN}_-\cap\widetilde{\mcB}_+\overline{w_0}\widetilde{\mcB}_+$, hence of $\widetilde{\mcN}_-$.
  For $U=U_-U_0U_+$ in the big cell $\widetilde{\mcN}_-\widetilde{\mcH}\widetilde{\mcN}_+$, let $U_-$ denote its negative Gauss factor.
  On the open set $\mathcal{U}$ of those $U$ with $U_-\in\mathcal{U}_-$, write $U_-=\prod_m y_{\widetilde{\jmath}_m}(t_m)$ and let~$V(U)$ be the product of the factors coming from the letters of $\widetilde{\boldsymbol{\jmath}}$ lifted from $j_1,\dots,j_{\len(w_0^K)}$.
  The $t_m$ are regular functions on $\mathcal{U}$, as $U\mapsto U_-$ is regular on the big cell, so this defines a morphism $V\colon\mathcal{U}\to\widetilde{\mcN}_-^K$; moreover the product of the remaining factors lies in $\widetilde{\mcB}_+\overline{w_0^Kw_0}\widetilde{\mcB}_+$, since $y_k(t)\in\widetilde{\mcB}_+\overline{s_k}\widetilde{\mcB}_+$ for $t\neq0$ and Bruhat cells multiply along reduced words.
  Similarly $V(U)\in\widetilde{\mcB}_+^K\overline{w_0^K}\widetilde{\mcB}_+^K=\widetilde{\mcN}_+^K\overline{w_0^K}^{\,-1}\widetilde{\mcB}_+^K$, so $V(U)\overline{w_0^K}$ lies in the big cell $\widetilde{\mcN}_+^K\widetilde{\mcB}_-^K$ of $\tG^K$, and $\overline{V}(U):=\lfloor V(U)\overline{w_0^K}\rfloor_+$ defines a morphism $\overline{V}$ on $\mathcal{U}$, as Gauss factors are regular on the big cell.
  Finally $\mathcal{U}$ contains the elements $y^C_{j_1}(t_1)\cdots y^C_{j_{\len(w_0)}}(t_{\len(w_0)})$ of $\mcG$ with $t_m\neq0$.
  By the defining factorizations, $V=V_{\rr}$ and $\overline{V}=\overline{V}_{\rr}$, so $H_{\rr}$ is a morphism on $\mathcal{U}$.

  On the common domain constructed above, let $U\in\mcG$.
  The factors $V_{\rr}(U)$, $M_{\rr}=V_{\rr}(U)^{-1}U$ and $\overline{V}_{\rr}(U)$ lie in $\mcG$, being fixed by $\vt$.
  Moreover,
  \[
    \widetilde{\mcN}_\pm^{\widetilde{\Gamma}_{\rr}}\cap\mcG=\mcN_\pm^{\Gamma_{\rr}}:
  \]
  both sides are connected unipotent subgroups.
  The left one is the image under the $\vt$-equivariant isomorphism $\exp$ of the $\Vt$-fixed part of $\widetilde{\mfn}_\pm^{\widetilde{\Gamma}_{\rr}}$.
  Both have the same Lie algebra $\mfn_\pm^{\Gamma_{\rr}}$, by the fixed-point calculation in the proof of \cref{prop: lift}(b).
  Thus $V_{\rr}(U)\in\mcN_-^{\Gamma_{\rr}}$ and $\overline{V}_{\rr}(U)\in\mcN_+^{\Gamma_{\rr}}$.
  Also $M_{\rr}\in\mcB_+\overline{\sigma_{\rr}}\mcB_+$, since
  \[
    \widetilde{\mcB}_+\overline{\sigma_{\rr}}\widetilde{\mcB}_+\cap\mcG
    =\mcB_+\overline{\sigma_{\rr}}\mcB_+,
  \]
  by the Bruhat decompositions of $\tG$ and $\mcG$ and the inclusions $\mcB_+\overline{w}\mcB_+\subseteq\widetilde{\mcB}_+\overline{w}\widetilde{\mcB}_+$ for $w\in W^\tau$.
  The representative of $w_0^{\widetilde{\Gamma}_{\rr}}$ is that of the longest element $w_0^{\Gamma_{\rr}}$ of $W_{\Gamma_{\rr}}$.
  Uniqueness of the factorizations therefore identifies $V_{\rr}(U)$ and $\overline{V}_{\rr}(U)$ with the factors defining $h_\gamma$ in \eqref{eq: h}.
  The same argument applies to columns.
  Finally, $\boldsymbol\gamma_{\rr},\boldsymbol\gamma_{\cc}$ restrict to the corresponding homomorphisms for $\mcG$ by \cref{prop: lift}(b), so~$h_{\tg}(U)=h_\gamma(U)$.
\end{proof}

\section{Folding seeds}
\label{sec: folding}
In this section $Y$ is an irreducible affine variety with an involutive automorphism $\vt$, and $X\subseteq Y$ is an irreducible closed subvariety fixed pointwise by $\vt$.
We write $\overline{f}:=f|_X$.
In our application $Y=\SL$ and~$X=\Sp$.
We adapt Dupont's folding construction \cite{Dup08} to seeds of functions.
For the definitions concerning matrices, we use the restatement in \cite[Section~5.1]{NS14}, whose notation is close to ours; see also \cite[Section~2.3]{Vie18} for a further account.

\begin{remark}
  \label{rmk: folding generality}
  One can show that the results below extend to skew-symmetrizable principal parts and to automorphisms of finite order.
  For the latter, mutations within an orbit need not commute.
  One must require that the entries between indices of each mutated orbit vanish.
  One must also impose admissibility by requiring that, in each row, the entries in the columns of that orbit agree weakly in sign.
  The Poisson and toric arguments use the same orbit-sum formulas.
  We retain skew-symmetric principal parts and involutions, which are simpler to state and suffice for our application.
\end{remark}

\begin{definition}
  \label{def: symmetric seed}
  Let $\Sigma=(\bx,\tB)$ be a seed in $\bbC(Y)$, with index set $\mcI=\mcI_\ex\sqcup\mcI_\fr$ and skew-symmetric principal part.
  The seed $\Sigma$ is \newword{$\vt$-symmetric} if there is an involution $\tau$ of $\mcI$ such that
  \begin{enumerate}[(S1)]
    \item $\tau$ preserves $\mcI_\ex$ and $\mcI_\fr$, and $b_{\tau(k)\tau(\ell)}=b_{k\ell}$ for every entry $b_{k\ell}$ of $\tB$;
    \item $x_{\tau(k)}=x_k\circ\vt$ for all $k\in\mcI$.
  \end{enumerate}
\end{definition}

The entries of $\bx$ are algebraically independent, hence pairwise distinct, so (S2) determines $\tau$: it is the permutation of $\mcI$ through which $f\mapsto f\circ\vt$ acts on $\bx$.
We call it the \newword{involution} of $\Sigma$, and its orbits the \newword{orbits} of $\Sigma$.
For $k\in\mcI_\ex$, (S1) and skew-symmetry give $b_{k\tau(k)}=b_{\tau(k)k}=-b_{k\tau(k)}$, so $b_{k\tau(k)}=0$.
Thus the mutations within any exchangeable orbit $O$ commute, and the \newword{orbit mutation} $\mu_O:=\prod_{k\in O}\mu_k$ is well defined.

\begin{definition}
  \label{def: admissible}
  Let $\Sigma$ be a $\vt$-symmetric seed with involution $\tau$.
  We say that $\tau$ is \newword{admissible at} the exchangeable orbit $O$ if, for every orbit $R$, the entries of the $R\times O$ block of $\tB$ agree weakly in sign.
\end{definition}

This is the sign condition of \cite[(5.2)]{NS14} and \cite[Definition~2.3.2(ii)]{Vie18}, imposed only on the columns of $O$ but including the frozen rows.
By (S1), for $\ell\in R$ and $k\in O$ the $R\times O$ block consists of the entries~$b_{\ell k}$ and $b_{\ell\tau(k)}$, so $\tau$ is admissible at $O$ if and only if $b_{\ell k}b_{\ell\tau(k)}\ge0$ for all $\ell\in\mcI$; in particular it is admissible at every orbit of size one.
Admissibility of $\tau$ at an exchangeable orbit $R\neq O$ need not be preserved by $\mu_O$ \cite[Remark~2.18]{Dup08}, that is, admissible automorphisms need not be stable \cite[Definition~2.3.4]{Vie18}.
We do not need stability: we apply \cref{prop: fold mutation} only to seeds whose involution is checked directly to be admissible at the mutated orbit.

\begin{definition}
  \label{def: fold}
  The \newword{fold} of a $\vt$-symmetric seed $\Sigma$ whose extended cluster consists of regular functions on $Y$ is $\overline{\Sigma}=(\overline{\bx},\overline{B})$, with $\overline{\bx}=(\overline{x}_O)_{O\in\mcI/\tau}$, $\overline{x}_O:=\overline{x}_k$, and $\overline{B}=(\overline{b}_{RO})_{R\in\mcI/\tau,O\in\mcI_\ex/\tau}$, $\overline{b}_{RO}:=\sum_{\ell\in R}b_{\ell k}$, for any $k\in O$.
\end{definition}

These definitions do not depend on $k\in O$ by (S1) and (S2), as $x_k\circ\vt$ and $x_k$ agree on $X$.
The fold is a seed in $\bbC(X)$ precisely when the functions $\overline{x}_O$ are algebraically independent and generate $\bbC(X)$.
These properties do not follow from symmetry alone.
Since $B$ is skew-symmetric, the principal part of $\overline{B}$ is skew-symmetrizable, $|O|\overline{b}_{RO}=-|R|\overline{b}_{OR}$; it is the folded matrix of \cite[(5.5)]{NS14}.

\begin{proposition}
  \label{prop: fold mutation}
  Let $\Sigma$ be a $\vt$-symmetric seed, belonging to a regular cluster structure on $Y$, whose fold is a seed in $\bbC(X)$, and let $O$ be an exchangeable orbit.
  Then $\mu_O(\Sigma)$ is $\vt$-symmetric.
  If $\tau$ is admissible at~$O$, then it is also admissible at $O$ for $\mu_O(\Sigma)$, and the extended cluster and the matrix of the fold of~$\mu_O(\Sigma)$ are obtained from those of $\overline{\Sigma}$ by the mutation in direction $O$.
\end{proposition}
\begin{proof}
  Choose $k\in O$.
  Since $b_{k\tau(k)}=0$, the mutations at the distinct indices in $O$ do not change each other's rows and columns, so $\mu_O(\tB)$ satisfies (S1).
  The new variables at $k$ and $\tau(k)$ are given by exchange relations swapped by $\tau$ and $\vt$, so (S2) holds.
  The new variables are regular since the cluster structure is, so the fold of $\mu_O(\Sigma)$ is defined.
  Admissibility at $O$ is preserved since the entries $b_{\ell k}$, $k\in O$, only change sign.
  If $\tau$ is admissible at $O$ and $k\in O$, then $\sum_{\ell\in R}[\pm b_{\ell k}]_+=[\pm\overline{b}_{RO}]_+$ for every orbit $R$, so the restriction to $X$ of the exchange relation at $k$ is the exchange relation of $\overline{\Sigma}$ in direction $O$.
  For the matrix statement, let $R,S\neq O$, with $S$ exchangeable, and choose $m\in S$.
  Admissibility at $O$ and (S1) give
  \[
    \sum_{\ell\in R,\,k\in O}[b_{\ell k}]_+[b_{km}]_+
    =[\overline{b}_{RO}]_+[\overline{b}_{OS}]_+,
  \]
  and the same identity with both entries negated.
  Summing the ordinary matrix-mutation formula therefore gives the mutation formula for $\overline{B}$.
  Entries with row or column orbit $O$ simply change sign.
  This uses only admissibility at $O$, rather than at every orbit; compare \cite[proof of Theorem~2.24]{Dup08}.
\end{proof}

\begin{proposition}
  \label{prop: fold compatible}
  Let $\{\cdot,\cdot\}$ be a Poisson bracket on $\bbC(Y)$ for which $\vt$ is a Poisson automorphism, and~$\{\cdot,\cdot\}_X$ a Poisson bracket on $\bbC(X)$ with $\{f,g\}|_X=\{\overline{f},\overline{g}\}_X$ for all $\vt$-invariant $f,g\in\bbC[Y]$.
  Let~$\Sigma$ be a $\vt$-symmetric seed whose extended cluster consists of regular functions on $Y$, such that $\overline{\bx}$ is algebraically independent and generates $\bbC(X)$, and suppose that $\{\cdot,\cdot\}$ is log-canonical with respect to~$\bx$, with coefficient matrix $\Omega$ satisfying \eqref{eq: compatibility} with all $\kappa_\ell\neq0$.
  Then $\{\cdot,\cdot\}_X$ is log-canonical with respect to $\overline{\bx}$, with coefficients~$\overline{\omega}_{OR}=\frac1{|R|}\sum_{\ell\in R}\omega_{k\ell}$, $k\in O$, and these satisfy \eqref{eq: compatibility} with $\overline{\kappa}_O=\kappa_k/|O|$, $k\in O$ exchangeable.
  In particular $\overline{B}$ has full rank, and $\{\cdot,\cdot\}_X$ is compatible with the cluster structure defined by~$\overline{\Sigma}$.
\end{proposition}
\begin{proof}

  As $\vt$ is Poisson, $\omega_{\tau(k)\tau(\ell)}=\omega_{k\ell}$ and $\kappa_{\tau(k)}=\kappa_k$, so the folded coefficients do not depend on the representative $k\in O$.
  The functions $f_O:=\prod_{k\in O}x_k$ are regular (\cref{def: fold}) and $\vt$-invariant, with~$\overline{f}_O=\overline{x}_O^{|O|}$ and $\{f_O,f_R\}=\big(\sum_{k\in O,\ell\in R}\omega_{k\ell}\big)f_Of_R$; restricting to $X$ gives $\{\overline{x}_O,\overline{x}_R\}_X=\overline{\omega}_{OR}\overline{x}_O\overline{x}_R$.
  For $k\in O$ exchangeable and any orbit $S$,
  \[
    \sum_R\overline{b}_{RO}\overline{\omega}_{RS}=\frac1{|S|}\sum_{\ell\in S}\sum_{m\in\mcI}b_{mk}\omega_{m\ell}=\frac1{|S|}\sum_{\ell\in S}\kappa_k\delta_{k\ell},
  \]
  which is $\kappa_k/|O|$ if $S=O$ and $0$ otherwise.
  As all $\overline{\kappa}_O$ are nonzero, full rank and compatibility follow from \cite[Proposition~2.5]{GSV18dd}, see \cref{sec: cluster conventions}.
\end{proof}

\begin{proposition}
  \label{prop: fold toric}
  Let $\Sigma$ be a $\vt$-symmetric seed whose fold is a seed, and $\mathsf T$ a torus acting on $Y$ such that every $x_k$ is $\mathsf T$-homogeneous of weight $w_k$, with $\tB^TP=0$ for the matrix $P=(w_k)_k$ whose rows are the coordinates of the $w_k$ in a basis of the character lattice, i.e.\ $\sum_{\ell\in\mcI}b_{\ell k}w_\ell=0$ for $k\in\mcI_\ex$.
  Let $\mathsf T_X\subseteq\mathsf T$ be a subtorus preserving $X$ and commuting with $\vt$.
  Then $\overline{x}_O$ is $\mathsf T_X$-homogeneous of weight $\overline{w}_O:=w_k|_{\mathsf T_X}$ for $k\in O$, and $\overline{B}^T\overline{P}=0$ for $\overline{P}=(\overline{w}_O)_O$; if the $w_k$ span the character space of $\mathsf T$ over $\bbC$, the $w_k|_{\mathsf T_X}$ span that of $\mathsf T_X$.
\end{proposition}

\begin{proof}
  Let $t\in\mathsf T_X$.
  As $t$ preserves $X$, for $k\in O$ and $p\in X$ we get $\overline{x}_O(t\cdot p)=x_k(t\cdot p)=t^{w_k}\,\overline{x}_O(p)$, so~$\overline{x}_O$ is $\mathsf T_X$-homogeneous of weight $w_k|_{\mathsf T_X}$.
  This weight does not depend on $k\in O$: as $t$ commutes with $\vt$,
  \[
    x_{\tau(k)}(t\cdot y)=x_k(t\cdot\vt(y))=t^{w_k}x_{\tau(k)}(y)\qquad(y\in Y),
  \]
  so the nonzero function $x_{\tau(k)}$ has $\mathsf T_X$-weight both $w_k|_{\mathsf T_X}$ and $w_{\tau(k)}|_{\mathsf T_X}$, and these agree.
  Hence $\overline{w}_R=w_\ell|_{\mathsf T_X}$ for every $\ell\in R$, and, for an exchangeable orbit $O$ and $k\in O$,
  \[
    \sum_R\overline{b}_{RO}\overline{w}_R=\sum_R\sum_{\ell\in R}b_{\ell k}\,w_\ell|_{\mathsf T_X}=\Big(\sum_{\ell\in\mcI}b_{\ell k}w_\ell\Big)\Big|_{\mathsf T_X}=0.
  \]
  Finally, $\mathsf T_X$ is a direct factor of $\mathsf T$, so restriction of characters from $\mathsf T$ to $\mathsf T_X$ is surjective, and the images of a spanning set of the character space of $\mathsf T$ span that of $\mathsf T_X$; these images are the $\overline{w}_O$.
\end{proof}

\section{Symmetric BD word seeds and their folds}
\label{sec: symmetric}
From now on $\gamma$ is an aperiodic BD pair of type $C_n$, and $\tg$ is its lift.
Let $\bfi$ be a double reduced word of type $C_n$.
The lift is order preserving: the letters lifted from $i_j$ precede those lifted from $i_{j'}$ if $j<j'$.
The letters at the level $x$ of $\tbfi$ are those lifted from the letters $\pm\min(x,\tau(x))$ of $\bfi$, with the same signs.
For a slot $z$ of $\bfi$ let $\widehat z$ be the slot of $\tbfi$ just after the letters lifted from the letters at the positions $\le z$.
Every chamber of $\tbfi$ contains such a slot, and the chamber $c_x(\widehat z)$ of $\tbfi$ only depends on the chamber $c_{x'}(z)$ of $\bfi$ for~$x'=\min(x,\tau(x))$.
We say that it \newword{lies over} $c_{x'}(z)$.
We also write $\tau$ for the induced involution of~$C(\tbfi)$, given by $\tau(c_x(\widehat z))=c_{\tau(x)}(\widehat z)$.
Its fixed points are the chambers at the level $n$.
So the $\tau$-orbits of~$C(\tbfi)$ are the sets of chambers lying over the chambers of $\bfi$, and exchangeable chambers and end chambers lie over exchangeable chambers and end chambers.
The same holds for the lift of an arbitrary word in~$\pm[1,n]$.
For $k\in C(\bfi)$ we write $O_k$ for the orbit of the chambers of $\tbfi$ lying over $k$, so that $k\mapsto O_k$ is a bijection~$C(\bfi)\to C(\tbfi)/\tau$.
We write $\widetilde{\mcI}:=C(\tbfi)$, with $\widetilde{\mcI}_\ex$ and $\widetilde{\mcI}_\fr$ as in \cref{def: BD word seed} for $\tbfi$ and $\tg$.

\subsection{Symmetry}
For a BD pair $\gamma'$ of type $A_{2n-1}$ let $\tau(\gamma'):=(\tau\gamma'_{\rr}\tau,\tau\gamma'_{\cc}\tau)$.
For a word $\bfj$ let $\tau(\bfj)$ be obtained by replacing every letter $\pm x$ by $\pm\tau(x)$.
Identify $C(\bfj)$ with $C(\tau(\bfj))$ by mapping $c_x(z)$ to~$c_{\tau(x)}(z)$.

\begin{lemma}
  \label{lem: naturality}
  Let $\gamma'$ be an aperiodic BD pair of type $A_{2n-1}$ and $\bfj$ a double reduced word for $\SL$.
  With this identification, $\tB_{\tau(\gamma')}(\tau(\bfj))=\tB_{\gamma'}(\bfj)$ and $\Phi_{\tau(\gamma')}(k;\tau(\bfj))=\Phi_{\gamma'}(k;\bfj)\circ\vt$ for every $k$.
\end{lemma}
\begin{proof}
  The elements $u_{\le z}$, $v_{>z}$ for $\tau(\bfj)$ are the images under $\tau$ of those for $\bfj$.
  Hence the weights of a chamber for $\tau(\bfj)$ are the images under $\tau=-w_0$ of its weights for $\bfj$, and $\Delta(k;\tau(\bfj))=\Delta(k;\bfj)\circ\vt$ by \cref{lem: minors restrict}.
  As $\tau(u)(\alpha_{\tau(x)})=\tau(u\alpha_x)$, the inversion letter of the root $\tau(a)$ of $\tau(\gamma')$ in $\tau(\bfj)$ is at the position of that of $a$ in $\bfj$.
  The thawed chambers, which are at the levels $\tau(a)^*=\tau(a^*)$, the extended words and the dipole matrices correspond.
  As $\tau$ preserves the pairing, the exponents correspond.
  Finally, \eqref{eq: B entries} only depends on the levels through the Cartan entries, and $a_{\tau(x)\tau(y)}=a_{xy}$.
\end{proof}

\begin{theorem}
  \label{thm: symmetric}
  The BD word seed $\Sigma_{\tg}(\tbfi)$ is a $\vt$-symmetric seed of $\mcC_{\tg}$ with involution $\tau$, and its extended cluster consists of regular functions on $\SL$.
\end{theorem}
\begin{proof}
  By \cref{cor: type A seeds,thm: GSV}, $\Sigma_{\tg}(\tbfi)$ is a seed of the regular cluster structure $\mcC_{\tg}$.
  By \cref{lem: lifted word}, $\tau(\tbfi)$ is obtained from $\tbfi$ by commutations swapping the two letters of each lifted pair, and~$\tau(\tg)=\tg$ by \cref{prop: lift}(a).
  By \cref{thm: moves}(a), these commutations preserve the exchange matrix and the functions $\Phi_{\tg}$ attached to corresponding chambers.
  The correspondence maps the chamber of $\tau(\tbfi)$ at the level $x$ containing a slot $\widehat z$ to the chamber of $\tbfi$ at the level $x$ containing $\widehat z$, as the slots~$\widehat z$ lie outside the commuted pairs.
  Composing with \cref{lem: naturality}, which maps $c_x(\widehat z)$ of $\tbfi$ to $c_{\tau(x)}(\widehat z)$ of~$\tau(\tbfi)$, we get $b_{\tau(\widetilde{k})\tau(\widetilde{\ell})}=b_{\widetilde{k}\widetilde{\ell}}$ and $\Phi_{\tg}(\tau(\widetilde{k});\tbfi)=\Phi_{\tg}(\widetilde{k};\tbfi)\circ\vt$.
  The involution $\tau$ preserves $\widetilde{\mcI}_\ex$ and $\widetilde{\mcI}_\fr$, as it maps thawed chambers to thawed chambers.
  So (S1) holds, and (S2) follows from \cref{prop: h theta} as~$x_{\tg}(\tau(\widetilde{k});\tbfi)=\Phi_{\tg}(\widetilde{k};\tbfi)\circ\vt\circ h_{\tg}^{-1}=x_{\tg}(\widetilde{k};\tbfi)\circ\vt$.
\end{proof}

\subsection{The fold}
We compute the fold of $\Sigma_{\tg}(\tbfi)$ in two steps: \cref{lem: index minors,lem: fold BFZ} treat the chamber minors and the matrices of ordinary words, and \cref{thm: fold} adds the dipoles and the exponents.

\begin{lemma}
  \label{lem: index minors}
  If $\widetilde{k}\in C(\tbfi)$ lies over $k\in C(\bfi)$, then $\Delta(\widetilde{k};\tbfi)|_\mcG=\Delta(k;\bfi)$.
  In particular the functions~$\Delta(\widetilde{k};\tbfi)|_\mcG$, one for each orbit, are algebraically independent and generate $\bbC(\mcG)$.
\end{lemma}
\begin{proof}
  Let $\widetilde{k}=c_x(\widehat z)$ and $k=c_{x'}(z)$, so that $x'=\min(x,\tau(x))$.
  The elements $u_{\le\widehat z}$, $v_{>\widehat z}$ for $\tbfi$ are the elements $u_{\le z}$, $v_{>z}$ of $W^\tau$ attached to $\bfi$, so $\Delta(\widetilde{k};\tbfi)=\Delta_{u_{\le z}\varpi_x,v_{>z}\varpi_x}$.
  The first claim now follows from \cref{lem: minors restrict}, and the second from \cite[Theorem~2.10]{BFZ05}.
\end{proof}

\begin{lemma}
  \label{lem: fold BFZ}
  Let $\bfj$ be a word in $\pm[1,n]$ and $\widetilde{\boldsymbol{\jmath}}$ its lift.
  For every $\bfj$-exchangeable $\ell$, every $k\in C(\bfj)$ and every chamber $\widetilde{\ell}$ over $\ell$, the entry of $\tB(\bfj)$ at $(k,\ell)$ is $\sum_{\widetilde{k}}b_{\widetilde{k}\widetilde{\ell}}$, the sum running over the chambers $\widetilde{k}$ over $k$ and $b_{\widetilde{k}\widetilde{\ell}}$ being the entries of $\tB(\widetilde{\boldsymbol{\jmath}})$.
\end{lemma}
\begin{proof}
  By \eqref{eq: B entries}, an entry depends only on the signed sequences of letters at its row and column levels and on the Cartan entry between them.
  We show that the lift preserves these sequences, and that counting the lifted levels that contribute recovers the Cartan entry.

  Let $k$ and $\ell$ be at the levels $x$ and $y$, and let $\widetilde{y}\in\{y,\tau(y)\}$ be the level of the chosen chamber $\widetilde{\ell}$.

  If $x=y$, only the chamber over $k$ at the level $\widetilde{y}$ can contribute: any other chamber over $k$ is at a non-adjacent level.
  The lift preserves the order and signs of the letters at each level, so the first two cases of \eqref{eq: B entries} give the required equality.
  If $x\neq y$ and $a_{xy}=0$, no level over $x$ is adjacent to $\widetilde{y}$, and both sides vanish.

  Suppose now that $a_{xy}\neq0$ and $x\neq y$.
  Only the chambers $\widetilde{k}$ at levels $\widetilde{x}\in\{x,\tau(x)\}$ adjacent to $\widetilde{y}$ can contribute.
  Restricting $\bfj$ to the levels $x,y$ and its lift to $\widetilde{x},\widetilde{y}$ gives the same signed sequence after relabeling, with $k,\ell$ corresponding to $\widetilde{k},\widetilde{\ell}$.
  The third case of \eqref{eq: B entries} therefore gives the same entry, except that the factor $-a_{xy}$ is replaced by $1$.
  There are exactly $-a_{xy}$ such levels $\widetilde{x}$: two, $n-1$ and $n+1$, if~$(x,y)=(n-1,n)$, and one otherwise.
  The claim follows by summing these contributions.
\end{proof}
An analogous folding formula for exchange matrices and compatible skew-symmetric matrices of $\mathbf{i}$-box seeds of positive braids is proved in \cite[Theorem~3.15]{JO26} for admissible foldings of Cartan data.

\begin{theorem}
  \label{thm: fold}
  Index the fold of $\Sigma_{\tg}(\tbfi)$ by $C(\bfi)$ via $O_k\mapsto k$.
  \begin{enumerate}[(a)]
    \item The exchange matrix of the fold of $\Sigma_{\tg}(\tbfi)$ is $\tB_\gamma(\bfi)$.
    \item For $\widetilde{k}$ over $k$, $\Phi_{\tg}(\widetilde{k};\tbfi)|_\mcG=\Phi_\gamma(k;\bfi)$.
    \item The map $h_\gamma$ is birational, and the extended cluster of the fold of $\Sigma_{\tg}(\tbfi)$ is $(x_\gamma(k;\bfi))_k$.
  \end{enumerate}
  In particular the fold of $\Sigma_{\tg}(\tbfi)$ is the BD word seed $\Sigma_\gamma(\bfi)$ of \cref{def: BD word seed} for $\mcG$ and $\gamma$.
\end{theorem}
\begin{proof}
  \emph{Inversion letters, thawed chambers and virtual letters.} We begin by showing that inversion letters and thawed chambers lift in pairs, so the extended words agree up to commutations.
  Let $a$ be a row root of $\gamma$, $g=\gamma_{\rr}(a)$, and $i_p$ its inversion letter, at the position $p=p^\rr_a$ and level $x$, and put $u=u_{\le p-1}\in W^\tau$, so that $u(\alpha^C_x)=\alpha^C_g$.
  As $g<n$ and $u$ preserves root lengths, $x<n$, so $i_p$ lifts to the two letters $-x,-\tau(x)$.
  Since $s^C_k=s_ks_{\tau(k)}$ for $k<n$ and $s^C_n=s_n$, the negative letters of $\tbfi$ preceding these two letters give the product $u\in S_{2n}$.
  Hence the roots of the two letters are $u\alpha_x$ and $us_x\alpha_{\tau(x)}=u\alpha_{\tau(x)}$, as $x$ and $\tau(x)$ are not adjacent.

  As $u$ preserves $\mfh$, the positive root $u\alpha_x$ restricts to $u\alpha^C_x=\alpha^C_g$.
  It is a sum of simple roots, each restricting to a simple root of type $C_n$, so it is $\alpha_g$ or $\alpha_{\tau(g)}$, and $u\alpha_{\tau(x)}=\tau(u\alpha_x)$ is the other one, as~$u\in W^\tau$.
  Since $\tg_{\rr}(a)=g$ and $\tg_{\rr}(\tau(a))=\tau(g)$, the two letters are the inversion letters of $a$ and $\tau(a)$ for $\tg$, and they separate chambers over $c_x(p-1)$ and $c_x(p)$.
  The thawed chambers of $a$ and $\tau(a)$ are the final chambers of $\tbfi$ at the levels $\tau(a)$ and $a$, i.e.\ the orbit over the thawed chamber $t^\rr_a$ of $a$.
  The same holds for column roots.

  As the roots of $\tg$ are the roots $a$ of $\gamma$ and their images $\tau(a)$, the inversion letters of $\tg$ in $\tbfi$ come in consecutive pairs, in the order of the inversion letters of $\gamma$ in $\bfi$.
  By \cref{def: BD word seed}(b), the lift of $\bfi^\gamma$ and $\tbfi^{\tg}$ thus differ only in the order of the two virtual letters of each pair, at the non-adjacent levels $a,\tau(a)$, and so are related by commutations.

  \emph{(a)} By \cref{lem: fold BFZ} for $\bfj=\bfi^\gamma$ and \cref{lem: commutation}, the orbit sums of $\tB(\tbfi^{\tg})_{\widetilde{\mcI}}$ are the entries of $\tB(\bfi^\gamma)_\mcI$.
  The nonzero entries of $D_a$ and $D_{\tau(a)}$ for $\tbfi$ are swapped by $\tau$, and their rows and columns lie over those of the nonzero entries of $D_a$ for $\bfi$, by the description of inversion letters and thawed chambers.
  None of their rows and columns is at the level $n$.
  Indeed, their levels are $a,\tau(a)$ and $x,\tau(x)$, where $a<n$ and, as shown above, $x<n$.
  Consequently, for a nonzero entry of $D_\gamma(\bfi)$ in the row $k$ and the column $\ell$, and a chamber~$\widetilde{\ell}$ over $\ell$, exactly one entry of $D_{\tg}(\tbfi)$ in a row over $k$ and the column $\widetilde{\ell}$ is nonzero.
  It has the same sign, and all nonzero entries of $D_{\tg}(\tbfi)$ arise in this way.
  So the orbit sums of $D_{\tg}(\tbfi)$ are the entries of~$D_\gamma(\bfi)$.

  \emph{(b)} Let $\widetilde{k}=c_x(\widehat z)$ lie over $k=c_{x'}(z)$, and let $u=u_{\le z}\in W^\tau$, so that the row weight of $\widetilde{k}$ is~$u\varpi_x$ by the proof of \cref{lem: index minors}, and let $a$ be a row root of $\gamma$, $g=\gamma_{\rr}(a)$.
  Then $(u\varpi_x)(\alpha^\vee_{\tau(g)})=(u\varpi_{\tau(x)})(\alpha^\vee_g)$, as~$u$ commutes with $w_0$, and $e^\rr_a(\widetilde{k})+e^\rr_{\tau(a)}(\widetilde{k})=[(u\varpi_x)(\alpha_g^\vee)]_++[(u\varpi_x)(\alpha_{\tau(g)}^\vee)]_+$.
  The sets $S=u([1,x])$ and~$S'=u([1,\tau(x)])$ are nested, so the two pairings, which are $[g\in S]-[g+1\in S]$ and $[g\in S']-[g+1\in S']$ by \eqref{eq: pairing A}, do not have opposite signs, and the sum is
  \[
    \big[(u\varpi_x)(\alpha_g^\vee+\alpha_{\tau(g)}^\vee)\big]_+=\big[(u\varpi^C_{x'})(\alpha^{C\vee}_g)\big]_+=e^\rr_a(k),
  \]
  where $x'=\min(x,\tau(x))$ and $\alpha^{C\vee}_g=[E^C_g,F^C_g]=\alpha^\vee_g+\alpha^\vee_{\tau(g)}$.
  The same calculation applies to columns.
  By \cref{thm: symmetric}, $\Phi_{\tg}(\widetilde{t};\tbfi)|_\mcG$ is the same function for the two thawed chambers $\widetilde{t}$ of an orbit.
  By \cref{lem: index minors} and induction along the recursion \eqref{eq: Phi}, which terminates by \cref{lem: tails},
  \[
    \Phi_{\tg}(\widetilde{k};\tbfi)|_\mcG=\Delta(k;\bfi)\prod_a\Phi_\gamma(t^\rr_a;\bfi)^{e^\rr_a(k)}\prod_b\Phi_\gamma(t^\cc_b;\bfi)^{e^\cc_b(k)}=\Phi_\gamma(k;\bfi).
  \]

  \emph{(c)} For $\widetilde{k}\in O_k$, let $\overline{x}_k:=x_{\tg}(\widetilde{k};\tbfi)|_\mcG$; it is regular, and $\overline{x}_k\circ h_\gamma=\Phi_{\tg}(\widetilde{k};\tbfi)|_\mcG=\Phi_\gamma(k;\bfi)$ by (b) and \cref{prop: h theta}.
  The $\Phi_\gamma(k;\bfi)$ are $\dim\mcG$ algebraically independent generators of $\bbC(\mcG)$, as noted after \cref{def: BD word seed}.
  So $h_\gamma$ is dominant, and $h_\gamma^*$ is an injective endomorphism of $\bbC(\mcG)$ whose image contains the $\Phi_\gamma(k;\bfi)$.
  Hence $h_\gamma^*$ is an automorphism, $h_\gamma$ is birational, and $\overline{x}_k=\Phi_\gamma(k;\bfi)\circ h_\gamma^{-1}=x_\gamma(k;\bfi)$.
  Together with (a) and \cref{thm: symmetric}, this shows that the fold is the seed $\Sigma_\gamma(\bfi)$.
\end{proof}

\begin{example}[\cref{ex: Sp6} continued]
  The lift of $\gamma=(\{1\to2\},0)$ is $\tg=(\{1\to2,5\to4\},0)$.
  By \cref{lem: minors restrict}, $\Delta_{-\varpi_1,\varpi_1}$ and $\Delta_{2\varpi_2-\varpi_3,-\varpi_3}$ are the restrictions of the matrix minors $\Delta_{6,1}$ and $\Delta_{124,456}$.
  The thawed chambers of the roots $5$ and $1$ of $\tg$ in $\tbfi$ are the final chambers at the levels $1$ and $5$, with minors $\Delta_{6,1}$ and $\Delta_{23456,12345}$, which agree on $Sp_6$ by \cref{lem: minors restrict}.
  The chamber of $\tbfi$ over $c_3(1)$ has the minor $\Delta_{124,456}$, and its row weight pairs to $1$ with both $\alpha_2^\vee$ and $\alpha_4^\vee$.
  Hence its function is $\Delta_{124,456}\Delta_{6,1}\Delta_{23456,12345}$, which restricts to $\Phi_\gamma(c_3(1);\bfi)$.
  The exponent $2$ is the sum of the exponents of the two thawed chambers of an orbit, as in the proof of (b).
  The six chambers over $c_1(6)$, $c_2(0)$, $c_2(2)$ have a single one of the two factors.
\end{example}

\subsection{The initial word}
We choose a word for which every dipole stays in one half of the lifted diagram.
The absence of entries between the two halves will make the involution of the lifted seed admissible at every exchangeable orbit.

Let $K:=[1,n-1]$, so that $W_K\subseteq W^\tau$ is of type $A_{n-1}$; in $S_{2n}$, $w_0^K$ is the longest element of the parabolic subgroup generated by the $s_k$, $k\neq n$.
Since $w_0$ is central in $W^\tau$, $w_0w_0^K$ is an involution; in $S_{2n}$ it is the permutation $\sigma:c\mapsto c\pm n$ swapping the two halves of $[1,2n]$ and preserving the order inside each of them.

\begin{definition}
  \label{def: initial word}
  Consider the reduced words $\mathbf l:=(1,2,\dots,n-1,1,2,\dots,n-2,\cdots,1,2,1)$ for $w_0^K$ and $\mathbf w:=(n,n-1,\dots,1,n,n-1,\dots,2,\cdots,n,n-1,n)$ for $w_0w_0^K$.
  Let $\bfi_0:=(\mathbf w^T,\mathbf l^T,\mathbf l,\mathbf w)$, so that $\bfi_0^T=\bfi_0$.
\end{definition}

The negative part of $\bfi_0$ is a reduced word for $w_0w_0^K\,w_0^K=w_0$, and so is the positive part by symmetry.
So $\bfi_0$ is a double reduced word of type $C_n$.

\begin{lemma}
  \label{lem: halves}
  In $\tbfi_0$, the row and column inversion letters of a root $a<n$ of $\tg$ lie at levels $>n$, and those of a root $a>n$ at levels $<n$.
  Consequently, the entries of $\tB_{\tg}(\tbfi_0)$ between chambers at levels $<n$ and chambers at levels $>n$ vanish, and the involution of $\Sigma_{\tg}(\tbfi_0)$ is admissible at every exchangeable orbit.
\end{lemma}
\begin{proof}
  In decreasing order of position, the positive letters of $\tbfi_0$ first come from the lift of $\mathbf w$, and the product of their reflections is $\sigma$.
  Their roots are the positive roots $\beta$ with $\sigma\beta<0$, so they swap values~$c<d$ with $\sigma(c)>\sigma(d)$, i.e.\ $c\le n<d$, and none of them has a simple root $\alpha_g$, $g\neq n$.
  So the column inversion letter of a root $a$ of $\tg$, at a position $p$ and level $x$, lies in the lift of $\mathbf l$.
  There $v_{>p}=\sigma v'$ with $v'$ in the parabolic subgroup generated by the $s_k$, $k\neq n$.
  Then $\sigma v'(\alpha_x)=\alpha_g$, $g=\tg_{\cc}(a)$, gives $v'(\alpha_x)=\sigma(\alpha_g)$, which is $\alpha_{g+n}$ if $g<n$ and $\alpha_{g-n}$ if $g>n$, so $x>n$ if $a<n$ and $x<n$ if $a>n$.
  For rows we argue in the same way with the negative letters in increasing order of position, which by symmetry have the levels of the positive letters in decreasing order.
  The thawed chamber of a root $a<n$ lies at the level $\tau(a)>n$, in the same half as its inversion letter.
  Hence its dipole entries are in rows and columns of chambers at levels $>n$, and symmetrically for $a>n$.
  All other nonzero entries have row and column chambers whose levels differ by at most one (\cref{lem: local entries}).
  Admissibility follows by considering the two cases: for a non-fixed chamber $\widetilde{k}$ and an exchangeable orbit $\{\widetilde{\ell},\tau(\widetilde{\ell})\}$ of size two, at most one of $b_{\widetilde{k}\widetilde{\ell}}$ and $b_{\widetilde{k}\tau(\widetilde{\ell})}$ is nonzero, while for a fixed chamber $\widetilde{k}$, $b_{\widetilde{k}\tau(\widetilde{\ell})}=b_{\widetilde{k}\widetilde{\ell}}$.
\end{proof}

\section{Proof of \texorpdfstring{\cref{thm: main}}{the main theorem}}
\label{sec: main}
\subsection{The cluster structure}
\begin{definition}
  \label{def: C gamma}
  The cluster structure $\mcC_\gamma$ on $\mcG=\Sp$ is the one defined by the BD word seed $\Sigma_\gamma(\bfi_0)$ (\cref{def: initial word,def: BD word seed,thm: fold}).
  We write $\oA_\gamma$ for its upper cluster algebra.
\end{definition}

\begin{proposition}
  \label{prop: adjacent}
  For every exchangeable chamber $k$ of $\Sigma_\gamma(\bfi_0)$, $\mu_{O_k}\big(\Sigma_{\tg}(\tbfi_0)\big)$ is a $\vt$-symmetric seed of $\mcC_{\tg}$ whose fold is $\mu_k\big(\Sigma_\gamma(\bfi_0)\big)$.
  In particular the cluster variables of the seeds adjacent to $\Sigma_\gamma(\bfi_0)$ are restrictions to~$\mcG$ of cluster variables of $\mcC_{\tg}$, hence regular functions on $\mcG$.
\end{proposition}
\begin{proof}
  The claim follows from \cref{prop: fold mutation}, by \cref{thm: symmetric,thm: fold,lem: halves,thm: GSV}.
\end{proof}

\begin{theorem}
  \label{thm: structure}
  The cluster structure $\mcC_\gamma$ has the properties (i), (iii) and (iv) of the introduction:
  \begin{enumerate}[(i)]
    \item $\mcC_\gamma$ has $2n-|\Gamma_{\rr}|-|\Gamma_{\cc}|$ frozen variables, and its exchange matrices have full rank.
      \setcounter{enumi}{2}
    \item The global toric action on $\mcC_\gamma$ is generated by the action of $\mcH_\gamma$.
    \item $\mcC_\gamma$ is compatible with $\{\cdot,\cdot\}_{\gamma,s}$ for every $s\in\mcS(\gamma)$.
  \end{enumerate}
\end{theorem}
\begin{proof}
  (i) The frozen chambers of $\bfi_0$ are the $2n$ end chambers other than the $|\Gamma_{\rr}|+|\Gamma_{\cc}|$ thawed chambers.

  (iv) Let $s\in\mcS(\gamma)$ and $\widetilde{s}$ as in \cref{sec: Sp}.
  By \cref{thm: GSV}(iv), $\{\cdot,\cdot\}_{\tg,\widetilde{s}}$ is log-canonical with respect to the extended cluster of the seed $\Sigma_{\tg}(\tbfi_0)$ of $\mcC_{\tg}$, with coefficient matrix satisfying \eqref{eq: compatibility} with all $\kappa_\ell\neq0$.
  By \cref{prop: restriction,thm: symmetric,thm: fold}, \cref{prop: fold compatible} applies, and gives compatibility and full rank.

  (iii) By \cref{thm: GSV}(iii), the cluster variables of $\Sigma_{\tg}(\tbfi_0)$ are homogeneous for the action of $\widetilde{\mcH}_{\tg}$.
  The weights span its character space by the argument of \cref{sec: toric A}, as the cluster variables generate $\bbC(\SL)$.
  Also $\tB^TP=0$ for their weight matrix $P$.
  The subtorus $\mcH_\gamma=\widetilde{\mcH}_{\tg}\cap(\mcH\times\mcH)$ (\cref{prop: lift}(c), applied to $\gamma_{\rr}$ and $\gamma_{\cc}$) preserves $\mcG$ and commutes with $\vt$.
  By \cref{prop: fold toric}, the weights of the extended cluster of $\Sigma_\gamma(\bfi_0)$ for $\mcH_\gamma$ span a space of dimension $\dim\mfh_\gamma=2n-|\Gamma_{\rr}|-|\Gamma_{\cc}|$, and $\overline{B}^T\overline{P}=0$ for their weight matrix $\overline{P}$, where $\overline{B}=\tB_\gamma(\bfi_0)$.
  As $\overline{B}$ has full rank, this is a global toric action of maximal rank.
\end{proof}

\subsection{Regularity and completeness}
We obtain Laurent expressions for regular functions by restriction from type A, and use factoriality to exclude denominators in the upper cluster algebra; this proves regularity and completeness.

\begin{lemma}
  \label{lem: primes}
  Let $R$ be a factorial domain with $R^\times=\bbC^*$, and $(\bx,\tB)$ a seed in its fraction field with $\tB$ of full rank, such that $\bx\subseteq R\subseteq\mcL(\bx)\cap\bigcap_{k\in \mcI_\ex}\mcL(\mu_k(\bx))$.
  Then the exchangeable elements of $\bx$ are pairwise non-associate primes of $R$, and $x'_k/x_k\notin R$ for every $k\in \mcI_\ex$.
\end{lemma}
\begin{proof}
  As $\bx\subseteq R\subseteq\mcL(\bx)$, we have $R[x_k^{-1}\mid k\in \mcI_\ex]=\mcL(\bx)$.
  The units of $\mcL(\bx)$ are $c\prod_{k\in \mcI_\ex}x_k^{a_k}$, $c\in\bbC^*$, by algebraic independence, as frozen variables are not inverted, and those of $R[x_k^{-1}\mid k\in \mcI_\ex]$ are $c\prod_jq_j^{b_j}$, where $q_1,\dots,q_s$ represent the primes dividing $\prod_{k\in \mcI_\ex}x_k$.
  Comparing these descriptions, we obtain $s=|\mcI_\ex|$ and mutually inverse integer matrices $F=(F_{kj})$ and $G=(G_{jk})$ such that $x_k=c_k\prod_jq_j^{F_{kj}}$ and $q_j=c'_j\prod_{k\in \mcI_\ex}x_k^{G_{jk}}$ for $k\in \mcI_\ex$ and $j\in[1,s]$; moreover $F_{kj}\ge0$.
  Write $x_kx'_k=M_++M_-$, with distinct monomials $M_\pm$ not involving $x_k$ ($\tB$ has no zero column).
  The exchange monomials have disjoint supports, so $M_++M_-$ is coprime to every monomial in the adjacent Laurent ring and, being a sum of distinct monomials, is not a unit there.
  If $G_{jk}<0$, then $q_j=(x'_k)^{-G_{jk}}(M_++M_-)^{G_{jk}}M'$ with $M'$ a Laurent monomial in $\mu_k(\bx)\setminus\{x'_k\}$; this cannot belong to $\mcL(\mu_k(\bx))$, a contradiction.
  So $F$ and $G$ have non-negative entries.
  For every column $j$ of $F$, choose $k$ with $F_{kj}>0$; then $0=(FG)_{k\ell}\ge F_{kj}G_{j\ell}\ge0$ for~$\ell\neq k$, so the only non-zero entry of the row $j$ of $G$ is in the column $k$.
  Hence the invertible matrix $G$ is monomial, and so is $F=G^{-1}$; their non-zero entries are positive integers inverse to each other, hence equal to $1$.
  So $F$ is a permutation matrix, that is, the $x_k$ are pairwise non-associate primes.

  If $x'_k=x_kg$ with $g\in R$, then $g=(x'_k)^2/(M_++M_-)\in\mcL(\mu_k(\bx))$, which is similarly a contradiction.
\end{proof}

Let $\bx$ be the extended cluster of $\Sigma_\gamma(\bfi_0)$, with index set $\mcI=C(\bfi_0)$.

\begin{lemma}
  \label{lem: Laurent}
  $\bbC[\mcG]\subseteq\mcL(\bx)\cap\bigcap_{k\in \mcI_\ex}\mcL(\mu_k(\bx))$.
\end{lemma}
\begin{proof}
  Let $f\in\bbC[\mcG]$; as $\mcG$ is a closed subvariety of the affine variety $\SL$, restriction $\bbC[\SL]\to\bbC[\mcG]$ is surjective, so $f=F|_\mcG$ with $F\in\bbC[\SL]$.
  By \cref{thm: GSV}(ii), $F$ is a Laurent polynomial in the extended clusters $\by$ of $\Sigma_{\tg}(\tbfi_0)$ and of $\mu_{O_k}(\Sigma_{\tg}(\tbfi_0))$ for every $k\in \mcI_\ex$, polynomial in the frozen variables.
  By \cref{prop: adjacent}, restriction to $\mcG$ maps the variables of $\by$ to those of $\bx$, resp.\ $\mu_k(\bx)$, exchangeable to exchangeable and frozen to frozen; so $f$ is a Laurent polynomial of the same kind.
\end{proof}

\begin{theorem}
  \label{thm: complete}
  The cluster structure $\mcC_\gamma$ has the property (ii) of the introduction: $\oA_\gamma=\bbC[\mcG]$; in particular $\mcC_\gamma$ is regular and complete.
\end{theorem}
\begin{proof}
  By \cref{lem: Laurent}, $\bbC[\mcG]$ is contained in the upper bound of $\Sigma_\gamma(\bfi_0)$, which equals $\oA_\gamma$ by full rank (\cref{thm: structure}).
  Conversely, let $f\in\oA_\gamma$.
  The ring $R=\bbC[\mcG]$ is factorial \cite[Proposition~1 and Corollary]{Pop74}, \cite[Theorem~3]{KP83}, and its units are the non-zero constants: for a unit $f$ and $U\in\mcG$, the polynomials $t\mapsto f(U\exp(tE_i))$ and $t\mapsto f(U\exp(tF_i))$ are nowhere vanishing, hence constant, and these one-parameter subgroups generate $\mcG$.
  So, by \cref{prop: adjacent,lem: Laurent}, \cref{lem: primes} applies to $\bx$ and $R$.
  Write $f=g/M$ in lowest terms, with $g\in R$ and $M$ a product of exchangeable variables of $\bx$, which is possible as $f\in\mcL(\bx)$.
  Suppose that $x_k$ divides $M$.
  As $f\in\mcL(\mu_k(\bx))$, $f=g'/M'$ with $g'\in R$ and $M'$ a monomial in $x'_k$ and the exchangeable $x_\ell$, $\ell\neq k$.
  Then $gM'=g'M$, and since $x_k$ divides $M$ but not $g$, the prime $x_k$ divides $x'_k$ or some $x_\ell$, which \cref{lem: primes} excludes.
  So $f\in R$.
  Every cluster variable lies in $\oA_\gamma$ by the Laurent phenomenon in geometric type, polynomial in the frozen variables \cite[Proposition~11.2]{FZ03}, hence is regular.
\end{proof}

\subsection{Arbitrary words}
\begin{theorem}
  \label{thm: all words}
  For every double reduced word $\bfi$ of type $C_n$, $\Sigma_\gamma(\bfi)$ is a seed of $\mcC_\gamma$, and the variables of its extended cluster are restrictions of cluster variables of $\mcC_{\tg}$.
  Commutations do not change $\Sigma_\gamma(\bfi)$, sign swaps and braid moves of length $3$ act on it by the mutation at the replaced chamber, and braid moves of length $4$ by $\mu_k\mu_\ell\mu_k$, where $k$ and $\ell$ are the replaced chambers at the levels $n$ and $n-1$.
\end{theorem}
\begin{proof}
  The action of the moves is \cref{thm: moves}, as $h_\gamma$ is birational (\cref{thm: fold}(c)) and $\alpha_n$ is the long simple root.
  So $\Sigma_\gamma(\bfi)$ is a seed of $\mcC_\gamma$, as every double reduced word is connected to $\bfi_0$ by moves.
  By \cref{thm: fold}, $\Sigma_\gamma(\bfi)$ is the fold of $\Sigma_{\tg}(\tbfi)$, which is a seed of $\mcC_{\tg}$ by \cref{cor: type A seeds}.
\end{proof}

\section{Compatible Poisson brackets}
\label{sec: brackets}
In this section $\mcG$ is an arbitrary simply connected simple complex algebraic group and $\gamma$ an arbitrary BD pair of $\mcG$; aperiodicity is not needed.
We classify the Poisson brackets compatible with a cluster structure on $\mcG$ under three hypotheses on one of its seeds, and then check them for the cluster structures $\mcC_\gamma$ on~$SL_N$ and $\Sp$ (\cref{cor: classification}).
By block-indecomposability, the coefficients $\kappa_\ell$ in \eqref{eq: compatibility} of a compatible bracket are a common multiple of those of $\{\cdot,\cdot\}_{\gamma,s}$, and the toric action of $\mcH_\gamma$ accounts for the remaining difference.

For $\xi=(\xi_{\rr},\xi_{\cc})\in\mfh_\gamma$ let $\partial_\xi:=\partial^\rr_{\xi_{\rr}}+\partial^\cc_{\xi_{\cc}}$, the infinitesimal action of $\mcH_\gamma$.
For $\Xi=\sum_u\xi_u\otimes\eta_u\in\bigwedge^2\mfh_\gamma$ let~$\toric_\Xi(f,g):=\sum_u\partial_{\xi_u}f\,\partial_{\eta_u}g$, and write $\Xi=\Xi_\mathrm{rr}+\Xi_\mathrm{rc}+\Xi_\mathrm{cc}$ with $\Xi_\mathrm{rr}\in\bigwedge^2\mfh_{\gamma_{\rr}}$, $\Xi_\mathrm{rc}\in\mfh_{\gamma_{\rr}}\otimes\mfh_{\gamma_{\cc}}$, and~$\Xi_\mathrm{cc}\in\bigwedge^2\mfh_{\gamma_{\cc}}$, where $\mfh_{\gamma_{\rr}}\otimes\mfh_{\gamma_{\cc}}$ is embedded by $\zeta\otimes\zeta'\mapsto(\zeta,0)\wedge(0,\zeta')$.

\begin{theorem}
  \label{thm: classification}
  Let $\mcC$ be a cluster structure on $\mcG$ with a seed $\Sigma=(\bx,\tB)$, with index set $\mcI=\mcI_\ex\sqcup\mcI_\fr$, such that:
  \begin{enumerate}[(a)]
    \item for some $s_0\in\mcS(\gamma)$, $\{\cdot,\cdot\}_{\gamma,s_0}$ is compatible with $\mcC$, and its coefficient matrix in $\bx$ satisfies \eqref{eq: compatibility} with all $\kappa_\ell\neq0$;
    \item every cluster variable of $\mcC$ is $\mcH_\gamma$-homogeneous, and the weights of the variables of $\bx$ span the character space of $\mcH_\gamma$ over $\bbC$, which has dimension $|\mcI_\fr|$;
    \item the principal part of $\tB$ is block-indecomposable.
  \end{enumerate}
  Then a Poisson bracket on $\bbC(\mcG)$ is compatible with $\mcC$ if and only if it is of the form
  \[
    q\{\cdot,\cdot\}_{\gamma,s}+\toric_\Xi,
    \qquad
    q\in\bbC,\quad s\in\mcS(\gamma),\quad \Xi\in\textstyle\bigwedge^2\mfh_\gamma.
  \]
  If $q\neq0$, the bracket equals $q\{\cdot,\cdot\}_{\gamma,s'}+\toric_{\Xi_\mathrm{rc}}$ with $s'=(s_{\rr}-q^{-1}\Xi_\mathrm{rr},\,s_{\cc}-q^{-1}\Xi_\mathrm{cc})\in\mcS(\gamma)$, with $\Xi_\mathrm{rr}$ and $\Xi_\mathrm{cc}$ regarded as elements of $\mfh\wedge\mfh$.
  So it is a multiple of a bracket $\{\cdot,\cdot\}_{\gamma,s''}$, $s''\in\mcS(\gamma)$, if and only if $\Xi_\mathrm{rc}=0$.
\end{theorem}
\begin{proof}
  Let $s_0$ be as in (a), and $\Omega$ the coefficient matrix of $\{\cdot,\cdot\}_{\gamma,s_0}$ in $\bx$.
  By (a), $\tB$ has full rank (\cref{sec: cluster conventions}).
  By (b), the action of $\mcH_\gamma$ is a global toric action; so its weight matrix $P$ in $\bx$, written in a basis of the character lattice of $\mcH_\gamma$, satisfies $\tB^TP=0$, and it has rank $\dim\mfh_\gamma=|\mcI_\fr|$.
  In particular the coefficient matrix $\Upsilon$ of a toric bracket $\toric_\Xi$ in $\bx$, whose columns are combinations of those of $P$, satisfies~$\tB^T\Upsilon=0$: toric brackets contribute nothing to the $\kappa_\ell$ in \eqref{eq: compatibility}.

  \emph{Changing the Cartan data.} For $s'_{\rr}\in\mcS(\gamma_{\rr})$ we have $s'_{\rr}-s_{\rr}\in\mfh_{\gamma_{\rr}}\wedge\mfh_{\gamma_{\rr}}$, as $\mcS(\gamma_{\rr})=s_{\rr}+\mfh_{\gamma_{\rr}}\wedge\mfh_{\gamma_{\rr}}$ (\cref{sec: BD}).
  For $f,g$ of bidegrees $(\lambda^\rr_f,\lambda^\cc_f)$, $(\lambda^\rr_g,\lambda^\cc_g)$, replacing $s_{\rr}$ by $s'_{\rr}$ changes $\{f,g\}_{\gamma,s}$ by $-(\lambda^\rr_f\otimes\lambda^\rr_g)(s'_{\rr}-s_{\rr})fg$, by \eqref{eq: bracket}, i.e.\ it adds $\toric_{s_{\rr}-s'_{\rr}}$, with $s_{\rr}-s'_{\rr}$ regarded as an element of $\bigwedge^2\mfh_{\gamma_{\rr}}\subseteq\bigwedge^2\mfh_\gamma$.
  Similarly, replacing $s_{\cc}$ by~$s'_{\cc}\in\mcS(\gamma_{\cc})$ changes it by $-(\lambda^\cc_f\otimes\lambda^\cc_g)(s'_{\cc}-s_{\cc})fg$, as $s_{\cc}$ enters $r_\cc$ through the factor swap, i.e.\ it adds $\toric_{s_{\cc}-s'_{\cc}}$.
  Thus $\{\cdot,\cdot\}_{\gamma,s'}=\{\cdot,\cdot\}_{\gamma,s}+\toric_{s-s'}$ for $s,s'\in\mcS(\gamma)$, with $s-s'=(s_{\rr}-s'_{\rr},s_{\cc}-s'_{\cc})\in\bigwedge^2\mfh_{\gamma_{\rr}}\oplus\bigwedge^2\mfh_{\gamma_{\cc}}\subseteq\bigwedge^2\mfh_\gamma$.

  \emph{The brackets are Poisson and compatible.} Let $s\in\mcS(\gamma)$.
  For $y_{\rr}\in\mcH_{\gamma_{\rr}}$, $\mathrm{Ad}_{y_{\rr}}$ fixes $\mfh$ and commutes with~$\gamma_{\rr}$ on $\mfn_\pm$, as $y_{\rr}^{\gamma_{\rr}(\alpha)}=y_{\rr}^\alpha$, so it fixes $r_\rr$.
  Hence left multiplication by $y_{\rr}$ is a Poisson automorphism of $\{\cdot,\cdot\}_{\gamma,s}$, and similarly for right multiplication by $\mcH_{\gamma_{\cc}}$.
  So the vector fields $\partial_\xi$, $\xi\in\mfh_\gamma$, preserve~$\{\cdot,\cdot\}_{\gamma,s}$, and they pairwise commute, as $\mcH_\gamma$ is abelian.
  Hence $\toric_\Xi$ is a Poisson bracket whose Schouten bracket with $\{\cdot,\cdot\}_{\gamma,s}$ vanishes, and $q\{\cdot,\cdot\}_{\gamma,s}+\toric_\Xi$ is Poisson.
  It is log-canonical in every cluster, as $\{\cdot,\cdot\}_{\gamma,s}=\{\cdot,\cdot\}_{\gamma,s_0}+\toric_{s_0-s}$, all cluster variables are $\mcH_\gamma$-homogeneous, and $\{\cdot,\cdot\}_{\gamma,s_0}$ is compatible.

  \emph{Every compatible bracket is of this form.} Let $\{\cdot,\cdot\}'$ be compatible, with coefficient matrix $\Omega'$ in $\bx$.
  For $\ell$ exchangeable and $k\neq\ell$, $\{\log(M_++M_-),\log x_k\}'$ is constant, where $x_\ell x'_\ell=M_++M_-$; since~$M_\pm$ are distinct monomials, $\{\log M_+,\log x_k\}'=\{\log M_-,\log x_k\}'$, i.e.\ $\sum_{m}b_{m\ell}\omega'_{mk}=0$.
  So $\Omega'$ satisfies \eqref{eq: compatibility} for some $\kappa'_\ell$.
  The matrix $\tB^T\Omega'\tB$ is skew-symmetric with $(\ell,m)$ entry $\kappa'_\ell b_{\ell m}$, so $\kappa'_\ell b_{\ell m}=-\kappa'_mb_{m\ell}$, and similarly for the $\kappa_\ell$.
  As the $\kappa_\ell$ are nonzero, $\kappa'_\ell/\kappa_\ell=\kappa'_m/\kappa_m$ whenever $b_{\ell m}\neq0$.
  By hypothesis (c), the graph on $\mcI_\ex$ with these edges is connected, so this ratio is a constant $q$, i.e.\ $\kappa'_\ell=q\kappa_\ell$ for all $\ell$.
  Then $\Upsilon:=\Omega'-q\Omega$ satisfies $\tB^T\Upsilon=0$.
  By full rank, $\ker\tB^T$ has dimension $|\mcI_\fr|=\dim\mfh_\gamma$, and it is spanned by the columns of $P$.
  Let $(\zeta_u)_{u=1}^{\dim\mfh_\gamma}$ be the basis of $\mfh_\gamma$ dual to the basis of the character lattice of $\mcH_\gamma$ in which $P$ is written, so that $\partial_{\zeta_u}x_k=P_{ku}x_k$.
  As $\Upsilon$ is skew-symmetric with columns in $\ker\tB^T$, and $P$ has full column rank, $\Upsilon=PZP^T$ for a unique skew-symmetric matrix $Z=(Z_{uv})$.
  Let~$\Xi:=\sum_{u,v}Z_{uv}\,\zeta_u\otimes\zeta_v\in\bigwedge^2\mfh_\gamma$; then $\toric_\Xi(x_k,x_\ell)=(PZP^T)_{k\ell}\,x_kx_\ell$, so $\Upsilon$ is the coefficient matrix of~$\toric_\Xi$ in $\bx$, and $\{\cdot,\cdot\}'=q\{\cdot,\cdot\}_{\gamma,s_0}+\toric_\Xi$.

  \emph{The last claims.} For a bracket $q\{\cdot,\cdot\}_{\gamma,s}+\toric_\Xi$ with $q\neq0$, put $s'_{\rr}:=s_{\rr}-q^{-1}\Xi_\mathrm{rr}$ and $s'_{\cc}:=s_{\cc}-q^{-1}\Xi_\mathrm{cc}$; then $s'=(s'_{\rr},s'_{\cc})\in\mcS(\gamma)$ and
  \[
    q\{\cdot,\cdot\}_{\gamma,s}+\toric_\Xi=q\{\cdot,\cdot\}_{\gamma,s'}+\toric_{\Xi_\mathrm{rc}}.
  \]
  So the bracket is a multiple of $\{\cdot,\cdot\}_{\gamma,s'}$ if $\Xi_\mathrm{rc}=0$.
  Conversely, if $q\{\cdot,\cdot\}_{\gamma,s'}+\toric_{\Xi_\mathrm{rc}}=q'\{\cdot,\cdot\}_{\gamma,s''}$ for some $s''\in\mcS(\gamma)$, comparing the $\kappa_\ell$ gives $q=q'$, as toric brackets contribute nothing to them; then we get~$\toric_{\Xi_\mathrm{rc}}=q(\{\cdot,\cdot\}_{\gamma,s''}-\{\cdot,\cdot\}_{\gamma,s'})=\toric_{\Xi'}$ with $\Xi'=q(s'-s'')\in\bigwedge^2\mfh_{\gamma_{\rr}}\oplus\bigwedge^2\mfh_{\gamma_{\cc}}$, and, by uniqueness,~$\Xi_\mathrm{rc}=0$.

\end{proof}

\begin{corollary}
  \label{cor: classification}
  Property (v) of the introduction holds for the cluster structures $\mcC_\gamma$ of \cref{thm: GSV} on~$SL_N$ and of \cref{def: C gamma} on $\Sp$: the Poisson brackets compatible with $\mcC_\gamma$ are exactly ${q\{\cdot,\cdot\}_{\gamma,s}+\toric_\Xi}$ with $q\in\bbC$, $s\in\mcS(\gamma)$ and $\Xi\in\bigwedge^2\mfh_\gamma$.
\end{corollary}
\begin{proof}
  By \cref{thm: classification}, it suffices to check (a)--(c) for one seed of $\mcC_\gamma$.
  For $\mcG=SL_N$ and an aperiodic BD pair $\gamma$ of type $A_{N-1}$, take $\Sigma=\Sigma_\gamma(\bfi_\mathrm{GSV})$: (a) holds for every $s_0\in\mcS(\gamma)$ by \cref{thm: GSV}(iv), (b) by \cref{thm: GSV}(i) and \cref{sec: toric A}, and (c) by \cref{lem: indecomposable}.
  For $\mcG=\Sp$ and an aperiodic BD pair $\gamma$ of type $C_n$, take $\Sigma=\Sigma_\gamma(\bfi_0)$: (a) holds by the proof of \cref{thm: structure}(iv), with all $\kappa_\ell\neq0$ by \cref{prop: fold compatible}, (b) by \cref{thm: structure}(i) and the proof of \cref{thm: structure}(iii), and (c) by \cref{lem: indecomposable}.
\end{proof}

\begin{remark}
  \label{rmk: v}
  Let $\mcC$ satisfy the hypotheses of \cref{thm: classification}.
  By \cref{thm: classification}, a compatible bracket with~$q\neq0$ is a multiple of a bracket in (iv) of the introduction exactly when $\Xi_\mathrm{rc}=0$.
  Mixed terms exist for every pair, as $\dim\mfh_{\gamma_\bullet}=|\Gamma|-|\Gamma_\bullet|\ge1$; they correspond to the term $V_{12}$ in \cite[Section~4]{GSV11}.
  For instance, for $\gamma=(0,0)$ on $SL_2$ and $\Xi=\alpha_1^\vee\otimes\alpha_1^\vee\in\mfh_{\gamma_{\rr}}\otimes\mfh_{\gamma_{\cc}}$, the bracket $\toric_\Xi$ gives $\{x_{11},x_{12}\}=-2x_{11}x_{12}$ and $\{x_{11},x_{22}\}=0$, so it is not a multiple of the standard bracket.

  The original \cite[Conjecture~3.2(v)]{GSV11} concerns Poisson--Lie brackets; in our conventions $\{\cdot,\cdot\}_{\gamma,s}$ is a Sklyanin bracket when $\gamma_{\cc}=\gamma_{\rr}^{-1}$ and $s_{\cc}=-s_{\rr}$.
  For $q\neq0$, a bracket of \cref{thm: classification} can be Poisson--Lie only if $\gamma_{\cc}=\gamma_{\rr}^{-1}$, and then it is Poisson--Lie exactly when it is a multiple of a Sklyanin bracket of the BD class of $\gamma_{\rr}$ \cite[Lemma~4.3(ii)]{GSV11}; in particular, multiplicativity excludes mixed terms.
  For $q=0$, some toric brackets are Poisson--Lie as well: if $\gamma_{\cc}=\gamma_{\rr}^{-1}$ and $0\neq\Xi_0\in\bigwedge^2\mfh_{\gamma_{\rr}}$, then $\Xi_0(f,g)^\rr-\Xi_0(f,g)^\cc$ is the difference of two Sklyanin brackets of the BD class of $\gamma_{\rr}$, hence multiplicative, and it is Poisson and compatible with $\mcC$ by \cref{thm: classification}, but it is not a multiple of such a bracket, whose $R$-matrix has non-zero root components.
  Indeed, $\mfh_{\gamma_{\cc}}=\mfh_{\gamma_{\rr}}$ as $\gamma_{\cc}=\gamma_{\rr}^{-1}$, and $\Xi_0(f,g)^\rr-\Xi_0(f,g)^\cc=\toric_\Xi$ for the tensor~$\Xi$ given by~$\Xi_\mathrm{rr}=\Xi_0$, $\Xi_\mathrm{cc}=-\Xi_0$ and zero mixed component $\Xi_\mathrm{rc}$; this is the case $q=0$ of \cref{thm: classification}.
  Such~$\Xi_0$ exist when $\dim\mfh_{\gamma_{\rr}}\ge2$, for instance for $\gamma=(0,0)$ on $SL_3$ or $Sp_4$.
  Thus, for $q\neq0$, the compatible Poisson--Lie brackets are as in \cite[Conjecture~3.2(v)]{GSV11}, and for $q=0$ there are also these toric brackets.
\end{remark}

\bibliographystyle{amsalpha}
\bibliography{references}

\end{document}